\documentclass[12pt,oneside,reqno]{amsart}
\usepackage{amssymb}
\usepackage{pifont}
\usepackage{amsaddr}

\usepackage{amsmath,amssymb,amsthm,amscd,a4wide,upref}
\usepackage[toc]{appendix}

\usepackage{hyperref}
\hypersetup{colorlinks,citecolor=blue,filecolor=black,linkcolor=blue,urlcolor=blue}
\usepackage{enumerate}
\usepackage{mathtools}

\begin{document}

\newcommand{\ad}{{\rm ad}}
\newcommand{\cri}{{\rm cri}}
\newcommand{\End}{{\rm{End}\ts}}
\newcommand{\Rep}{{\rm{Rep}\ts}}
\newcommand{\Hom}{{\rm{Hom}}}
\newcommand{\Mat}{{\rm{Mat}}}
\newcommand{\ch}{{\rm{ch}\ts}}
\newcommand{\chara}{{\rm{char}\ts}}
\newcommand{\diag}{{\rm diag}}
\newcommand{\non}{\nonumber}
\newcommand{\wt}{\widetilde}
\newcommand{\wh}{\widehat}
\newcommand{\ot}{\otimes}
\newcommand{\la}{\lambda}
\newcommand{\La}{\Lambda}
\newcommand{\De}{\Delta}
\newcommand{\al}{\alpha}
\newcommand{\be}{\beta}
\newcommand{\ga}{\gamma}
\newcommand{\Ga}{\Gamma}
\newcommand{\ep}{\epsilon}
\newcommand{\ka}{\kappa}
\newcommand{\vk}{\varkappa}
\newcommand{\si}{\sigma}
\newcommand{\vs}{\varsigma}
\newcommand{\vp}{\varphi}
\newcommand{\de}{\delta}
\newcommand{\ze}{\zeta}
\newcommand{\om}{\omega}
\newcommand{\Om}{\Omega}
\newcommand{\ee}{\epsilon^{}}
\newcommand{\su}{s^{}}
\newcommand{\hra}{\hookrightarrow}
\newcommand{\ve}{\varepsilon}
\newcommand{\ts}{\,}
\newcommand{\pr}{^{\tss\prime}}
\newcommand{\vac}{\mathbf{1}}
\newcommand{\di}{\partial}
\newcommand{\qin}{q^{-1}}
\newcommand{\tss}{\hspace{1pt}}
\newcommand{\Sr}{ {\rm S}}
\newcommand{\U}{ {\rm U}}
\newcommand{\BL}{ {\overline L}}
\newcommand{\BE}{ {\overline E}}
\newcommand{\BP}{ {\overline P}}
\newcommand{\AAb}{\mathbb{A}\tss}
\newcommand{\CC}{\mathbb{C}\tss}
\newcommand{\KK}{\mathbb{K}\tss}
\newcommand{\QQ}{\mathbb{Q}\tss}
\newcommand{\SSb}{\mathbb{S}\tss}
\newcommand{\TT}{\mathbb{T}\tss}
\newcommand{\ZZ}{\mathbb{Z}\tss}
\newcommand{\DY}{ {\rm DY}}
\newcommand{\DX}{ {\rm DX}}
\newcommand{\ZDX}{ {\rm ZDX}}
\newcommand{\X}{ {\rm X}}
\newcommand{\Y}{ {\rm Y}}
\newcommand{\Z}{{\rm Z}}
\newcommand{\Ac}{\mathcal{A}}
\newcommand{\Lc}{\mathcal{L}}
\newcommand{\Pc}{\mathcal{P}}
\newcommand{\Qc}{\mathcal{Q}}
\newcommand{\Rc}{\mathcal{R}}
\newcommand{\Sc}{\mathcal{S}}
\newcommand{\Tc}{\mathcal{T}}
\newcommand{\Bc}{\mathcal{B}}
\newcommand{\Ec}{\mathcal{E}}
\newcommand{\Fc}{\mathcal{F}}
\newcommand{\Gc}{\mathcal{G}}
\newcommand{\Hc}{\mathcal{H}}
\newcommand{\Uc}{\mathcal{U}}
\newcommand{\Vc}{\mathcal{V}}
\newcommand{\Wc}{\mathcal{W}}
\newcommand{\Yc}{\mathcal{Y}}
\newcommand{\Ar}{{\rm A}}
\newcommand{\Br}{{\rm B}}
\newcommand{\Ir}{{\rm I}}
\newcommand{\Fr}{{\rm F}}
\newcommand{\Jr}{{\rm J}}
\newcommand{\Or}{{\rm O}}
\newcommand{\GL}{{\rm GL}}
\newcommand{\Spr}{{\rm Sp}}
\newcommand{\Rr}{{\rm R}}
\newcommand{\Zr}{{\rm Z}}
\newcommand{\ZX}{{\rm ZX}}
\newcommand{\gl}{\mathfrak{gl}}
\newcommand{\middd}{{\rm mid}}
\newcommand{\ev}{{\rm ev}}
\newcommand{\Pf}{{\rm Pf}}
\newcommand{\Norm}{{\rm Norm\tss}}
\newcommand{\oa}{\mathfrak{o}}
\newcommand{\spa}{\mathfrak{sp}}
\newcommand{\osp}{\mathfrak{osp}}
\newcommand{\f}{\mathfrak{f}}
\newcommand{\g}{\mathfrak{g}}
\newcommand{\h}{\mathfrak h}
\newcommand{\n}{\mathfrak n}
\newcommand{\z}{\mathfrak{z}}
\newcommand{\Zgot}{\mathfrak{Z}}
\newcommand{\p}{\mathfrak{p}}
\newcommand{\sll}{\mathfrak{sl}}
\newcommand{\agot}{\mathfrak{a}}
\newcommand{\qdet}{ {\rm qdet}\ts}
\newcommand{\Ber}{ {\rm Ber}\ts}
\newcommand{\HC}{ {\mathcal HC}}
\newcommand{\cdet}{{\rm cdet}}
\newcommand{\rdet}{{\rm rdet}}
\newcommand{\tr}{ {\rm tr}}
\newcommand{\gr}{ {\rm gr}\ts}
\newcommand{\str}{ {\rm str}}
\newcommand{\loc}{{\rm loc}}
\newcommand{\Gr}{{\rm G}}
\newcommand{\sgn}{ {\rm sgn}\ts}
\newcommand{\sign}{{\rm sgn}}
\newcommand{\ba}{\bar{a}}
\newcommand{\bb}{\bar{b}}
\newcommand{\eb}{\bar{e}}
\newcommand{\bi}{\bar{\imath}}
\newcommand{\bj}{\bar{\jmath}}
\newcommand{\bk}{\bar{k}}
\newcommand{\bl}{\bar{l}}
\newcommand{\hb}{\mathbf{h}}
\newcommand{\Sym}{\mathfrak S}
\newcommand{\fand}{\quad\text{and}\quad}
\newcommand{\Fand}{\qquad\text{and}\qquad}
\newcommand{\For}{\qquad\text{or}\qquad}
\newcommand{\OR}{\qquad\text{or}\qquad}
\newcommand{\grpr}{{\rm gr}^{\tss\prime}\ts}
\newcommand{\degpr}{{\rm deg}^{\tss\prime}\tss}

\numberwithin{equation}{section}

\newtheorem{thm}{Theorem}[section]
\newtheorem{lem}[thm]{Lemma}
\newtheorem{prop}[thm]{Proposition}
\newtheorem{cor}[thm]{Corollary}
\newtheorem{conj}[thm]{Conjecture}
\newtheorem*{mthm}{Main Theorem}
\newtheorem*{mthma}{Theorem A}
\newtheorem*{mthmb}{Theorem B}
\newtheorem*{mthmc}{Theorem C}
\newtheorem*{mthmd}{Theorem D}

\theoremstyle{definition}
\newtheorem{defin}[thm]{Definition}

\theoremstyle{remark}
\newtheorem{remark}[thm]{Remark}
\newtheorem{example}[thm]{Example}

\newcommand{\bthm}{\begin{thm}}
\newcommand{\ethm}{\end{thm}}
\newcommand{\bpr}{\begin{prop}}
\newcommand{\epr}{\end{prop}}
\newcommand{\ble}{\begin{lem}}
\newcommand{\ele}{\end{lem}}
\newcommand{\bco}{\begin{cor}}
\newcommand{\eco}{\end{cor}}
\newcommand{\bde}{\begin{defin}}
\newcommand{\ede}{\end{defin}}
\newcommand{\bex}{\begin{example}}
\newcommand{\eex}{\end{example}}
\newcommand{\bre}{\begin{remark}}
\newcommand{\ere}{\end{remark}}
\newcommand{\bcj}{\begin{conj}}
\newcommand{\ecj}{\end{conj}}

\newcommand{\bal}{\begin{aligned}}
\newcommand{\eal}{\end{aligned}}
\newcommand{\beq}{\begin{equation}}
\newcommand{\eeq}{\end{equation}}
\newcommand{\ben}{\begin{equation*}}
\newcommand{\een}{\end{equation*}}

\newcommand{\bpf}{\begin{proof}}
\newcommand{\epf}{\end{proof}}

\def\beql#1{\begin{equation}\label{#1}}

\title[Yangian doubles and their representations]
{Yangian doubles of classical types and their vertex representations}

\author{Naihuan Jing${^{1,2}}$, Fan Yang$^{\dagger 1}$, Ming Liu${^{1,3}}$}
\address{$^1$School of Mathematics, South China University of Technology, Guangzhou, Guangdong 510640, China}
\email{jing@math.ncsu.edu, ming.l1984@gmail.com, 1329491781@qq.com}
\address{$^2$Department of Mathematics, North Carolina State University, Raleigh, NC 27695, USA}
\address{$^3$School of Mathematics and Statistics, University of Sydney, NSW 2006, Australia}
%\email{jing@math.ncsu.edu}
\thanks{$\dagger$Corresponding author: Fan Yang}
%\\Supported in part by National Natural Science Foundation of China grant no. 11531004 and Simons Foundation grant no. 523868.}
\date{} % Start Octorber 2017
\maketitle

\vspace{5 mm}

\begin{abstract}
The Yangian double $\DY_\hbar(\g_N)$ is introduced for the classical types of $\g_N=\oa_{2n+1}$, $\spa_{2n}$, $\oa_{2n}$.
Via the Gauss decomposition of the generator matrix, the Yangian double is given the Drinfeld presentation.
In addition, bosonization of level $1$ realizations for the Yangian double $\DY_\hbar(\g_N)$ of non-simply-laced types are explicitly constructed.
\end{abstract}

\vspace{5 mm}

\section{Introduction}
\label{sec:int}

The quantum double is defined by gluing together a quantum group with its
Hopf dual (in certain sense) \cite{D} and this
provides conceptual foundation for the universal R-matrix of the quantum
enveloping algebra ~\cite{cp:gq}. The very construction not only
solves the quantum Yang-Baxter equation on each irreducible representation of the quantum group (see \cite{C}),
but has also spurred important applications
in representation theory \cite{FR2, FM}, the $q$-conformal field theory \cite{FR1} and knot theory \cite{R, RT}. %They can be also used in

The quantum doubles corresponding to finite dimensional simple Lie algebras provide
a unified construction of solutions for the quantum Yang-Baxter equation \cite{Jb, KR, C}.
The quantum enveloping algebras associated to affine Lie algebras are introduced in connection with
the trigonometric quantum R-matrices \cite{Jb1}. %the quangives rise to the quantum affine algebras.
%enveloping algebras of Drinfeld and Jimbo.
Besides the quantum enveloping algebras, Drinfeld also introduced another class of
quantum groups called the Yangians $\Y(\g)$    %$\Y_\hbar(\g)$
corresponding to the rational R-matrices, and Yangians since then have quickly developed
into one of the most important algebraic structures with numerous applications, see \cite{m:yc} for more information.

The Yangian algebra corresponding to the rational $R$-matrix is in fact some deformation of
the current algebra $\mathfrak g[t]$
of the simple Lie algebra $\mathfrak g$. The general linear Yangian algebra $\Y(\mathfrak{gl}_n)$ was first given by
Tarasov \cite{Ta} for $\mathfrak{gl}_2$ and then introduced in general by Drinfeld \cite{Dr}.
The simply-laced types of the Yangian $\Y(\g)$
have recently been studied by Guay et al \cite{GRW2} and the PBW theorem for the simply-laced affine types was obtained using
the vertex representation. There is also a uniform proof of PBW bases for several type $A$ algebras \cite{T} (see
\cite{Z} for the super type $A$).
A complete proof the isomorphism between the $RTT$ and Drinfeld presentation
for the Yangian algebra $\Y(\g)$ in type $A$ was obtained by Brundan and Kleshchev in \cite{bk:pp} using the Gauss decomposition of the generator matrix $T(u)$.

The quantum double $\DY(\mathfrak g)$ and its universal R-matrix in general were constructed % and its universal R-matrix in general
by Khoroshkin and Tolstoy in \cite{KT}, where explicit formulas for $\DY(\mathfrak{sl}_2)$ were given.
The quantum double $\DY_\hbar(\mathfrak{sl}_n)$ %of the Yangian $\Y(\mathfrak{sl}_n)$
was studied by Iohara %\cite{io:br}
in the $RTT$ formulism with a central extension, and the author also stated its Drinfeld commutation relations. % were stated in \cite{io:br}.
Higher level bosonization of $\DY_\hbar(\mathfrak{sl}_n)$
was also constructed in \cite{DHHZ}.
%The Drinfeld realization for the Yangian doubles is particularly interesting as it is useful to classify finite dimensional representations via the Drinfeld polynomials and construct universal R-matrices in this case \cite{RT}.
The Yangian algebra $\Y(\g)$ in types B, C, and D was studied \cite{amr:otr} in the $RTT$ format
together with their Poincar\'e-Birkhoff-Witt bases.
Recently the identification of the Drinfeld realization and RTT presentation of the Yangian algebras $\Y(\g)$ for types B, C, and D
has been proved in \cite{ji:iso}, which
makes it possible to study the general non-simply-laced Drinfeld realization of the Yangian double algebra for all classical types.
In this paper, we will construct a central extension of the Yangian double $\DY_\hbar(\g_N)$ of types $B,C$ and $D$ by using the $RTT$ relation under a unitarity condition.
The Drinfeld generators are constructed by using the Gauss decomposition, and we will prove that these two constructions are isomorphic, moreover,
we also present a Poincar\'e-Birkhoff-Witt (PBW) basis using monomials in the Weyl like generators $t_{ij}^{(r)}$ and $c$. Finally, we also construct level one modules of the Yangian double $\DY_\hbar(\g_N)$ in terms of bosons. As a by-product a detailed proof of all Drinfeld generators
and explicit construction
for the Yangian double in type $A$ is also obtained.

Vertex representations of quantum affine algebras of simply-laced types were constructed by Frenkel-Jing in \cite{FJ}
and other non-simply laced types
were realized subsequently \cite{Be, JK, JKM} in untwisted types. The Yangian analog in type $A$ was given in \cite{io:br} and similar
properties are generalized in \cite{Ko}. Vertex operator representations and Drinfeld realization of quantum affine superalgebras in simply-laced types were recently given in \cite{XZ}.
Our construction of vertex representations for the central extension of the Yangian doubles provides level one modules
for other classical (non-simply laced) types in this
paper, and we would like to remark that our construction seems to be a little bit different from those in \cite{Be, JK, JKM}.

To explain our construction, let $\g$ be the simple Lie
algebra associated with the Cartan matrix $A=[a_{ij}]_{i,j=1}^n$. Let $\al_1,\dots,\al_n$
be the corresponding simple roots.
The {\em Drinfeld Yangian double} $\DY_\hbar^{D}(\g)$ is topologically generated by
elements $h_{i\tss r}$,~$\xi_{i\tss r}^{+}$,~$\xi_{i\tss r}^{-}$ and the central element $c$ with $i=1,\dots,n$ and $r\in \mathbb{Z}$
subject to the defining relations
\begin{align*}
[H_{i}^{\pm}(u),H_{j}^{\pm}(v)]=0,
\end{align*}
\begin{align*}
(u_{\pm}-v_{\mp}-B_{ij}\hbar)(u_{\mp}-v_{\pm}+B_{ij}\hbar)H_{i}^{\pm}(u)H_{j}^{\mp}(v)\\
=(u_{\pm}-v_{\mp}+B_{ij}\hbar)(u_{\mp}-v_{\pm}-B_{ij}\hbar)H_{j}^{\mp}(v)H_{i}^{\pm}(u),
\end{align*}
\begin{align*}
H_{i}^{\pm}(u)^{-1}E_{j}(v)H_{i}^{\pm}(u)=\frac{u_{\mp}-v+B_{ij}\hbar}{u_{\mp}-v-B_{ij}\hbar}E_{j}(v),\\
H_{i}^{\pm}(u)F_{j}(v)H_{i}^{\pm}(u)^{-1}=\frac{u_{\pm}-v+B_{ij}\hbar}{u_{\pm}-v-B_{ij}\hbar}F_{j}(v),
\end{align*}
\begin{align*}
(u-v+B_{ij}\hbar)E_{i}(u)E_{j}(v)=(u-v-B_{ij}\hbar)E_{j}(v)E_{i}(u),\\
(u-v-B_{ij}\hbar)F_{i}(u)F_{j}(v)=(u-v+B_{ij}\hbar)F_{j}(v)F_{i}(u),
\end{align*}
\begin{align*}
\sum_{\sigma\in\mathfrak{S}_{m}}[E_{i}(u_{\sigma(1)}),[E_{i}(u_{\sigma(2)})\cdots,[E_{i}(u_{\sigma(m)}),E_{j}(v)]\cdots]=0, \\
\sum_{\sigma\in\mathfrak{S}_{m}}[F_{i}(u_{\sigma(1)}),[F_{i}(u_{\sigma(2)})\cdots,[F_{i}(u_{\sigma(m)}),F_{j}(v)]\cdots]=0, \\
~~ i\neq j,m=1-a_{ij},
\end{align*}
\begin{align*}
[E_{i}(u),F_{j}(v)]=\frac{1}{\hbar}\delta_{ij}\{\delta(u_{-}-v_{+})H_{i}^{-}(v_{+})-\delta(u_{+}-v_{-})H_{i}^{+}(u_{+})\},
\end{align*}
where $B_{ij}=\frac{1}{2}(\al_i,\al_j)$, $u_{\pm}=u\pm \frac14 \hbar c$ and
\begin{align*}
H^{+}_i(u)&=1+\hbar\sum_{r=0}^{\infty}h_{i\tss r}\ts u^{-r-1},\quad H^{-}_i(u)=1-\hbar\sum_{r=1}^{\infty}h_{i\tss -r}\ts u^{r-1},\\
E_i(u)&=\sum_{r\in \mathbb{Z}}\xi^{+}_{i\tss r}\ts u^{-r-1},\quad \quad \quad F_i(u)=\sum_{r\in \mathbb{Z}}\xi^{-}_{i\tss r}\ts u^{-r-1}.
\end{align*}

If $\g=\g_N$
is the orthogonal Lie algebra $\oa_N$ (with $N=2n$ or $N=2n+1$)
or symplectic Lie algebra $\spa_N$
(with even $N=2n$) then the algebra $\DY_\hbar^{R}(\g_N)$
(the {\em Yangian double in the $RTT$ presentation})
can be defined via the rational $R$-matrix. The defining relations take the form
of the $RTT$ {\em relation}
\begin{align}\label{RTTint}
R(u-v)\ts T^{\pm}_1(u)\ts T^{\pm}_2(v)&=T^{\pm}_2(v)\ts T^{\pm}_1(u)\ts R(u-v),\\
R(u_{+}-v_{-})T_{1}^{+}(u)T_{2}^{-}(v)&=T_{2}^{-}(v)T_{1}^{+}(u)R(u_{-}-v_{+}),
\non
\end{align}
together with the {\em unitarity condition}
\beql{unitaryint}
T^{\pm}(u+\kappa\hbar)^{t}T^{\pm}(u)=1,
\eeq
with the notation explained below in Section~\ref{sec:nd}. Here
$T^{\pm}(u)$ is a square matrix of size $N$ whose $(i,j)$ entries are the formal series
\ben
t^{+}_{ij}(u)=\delta_{ij}+\hbar\sum_{r=1}^{\infty}t^{(r)}_{ij}u^{-r},\quad t^{-}_{ij}(u)=\delta_{ij}-\hbar\sum_{r=1}^{\infty}t^{(-r)}_{ij}u^{r-1},
\een
so that the algebra $\DY_\hbar^{R}(\g_N)$ is topologically generated by all coefficients
$t_{ij}^{(r)}$ and $c$ subject to the defining relations
\eqref{RTTint} and \eqref{unitaryint}. Our goal is to show the following.

\begin{mthm}
The algebra $\DY_\hbar^{R}(\g_N)$ is isomorphic to the algebra $\DY_\hbar^{D}(\g_N)$.
\end{mthm}

The paper is organized as follows. In section~\ref{sec:nd}, we give the definitions of the extended Yangian double $\DX_\hbar^{R}(\g_N)$ and Yangian double $\DY_\hbar^{R}(\g_N)$. We also explain the notations to be used in the rest of the paper. In section~\ref{sec:PBW}, we prove the Poincar\'e-Birkhoff-Witt theorem for the algebras $\DX_\hbar^{R}(\g_N)$ and $\DY_\hbar^{R}(\g_N)$. In section~\ref{sec:et}, we give the imbedding theorem which will play a key role in the proof of the main theorem. In sections~\ref{sec:dpey} and \ref{sec:isom}, we prove the main theorem by calculating the relations for the Gaussian generators, which then imply the Drinfeld commutation relations. The surjectivity and injectivity are also proved step by step. As a by-product
we also obtain a detailed proof of the isomorphism between the Drinfeld and $RTT$ presentations in type $A$, which was originally stated in
\cite{io:br}. Section~\ref{sec:boson} contains construction of the level $1$ modules for the Yangian double $\DY_\hbar^{D}(\g_N)$, and they
show that the central extension of the Yangian double algebra is indeed nontrivial.

\section{Basic definitions}
\label{sec:nd}

Let $\mathfrak g_N$ be the orthogonal/symplectic Lie algebra: $\mathfrak g_N=\oa_N$ with $N=2n+1, 2n$ (type $B, D$ respectively) or
$\spa_N$ with $N=2n$ (type $C$). As a subalgebra of $\gl_N$,
$\mathfrak g_N$ is spanned by the elements $F_{i\tss j}$:
\beql{fij}
F_{i\tss j}=E_{i\tss j}-E_{j\pr i\pr}\Fand F_{i\tss j}
=E_{i\tss j}-\ve_i\ts\ve_j\ts E_{j\pr i\pr},
\eeq
for $\oa_N$ and $\spa_N$ respectively, and $E_{i\tss j}$ are the unit matrices in $End(\mathbb C^N)$.
%with $1$ sitting at the $(i, j)$ position and zero elsewhere.
The index flip is defined by $i\pr=N-i+1$, and %$\ve_i=1$ for the orthogonal case and or $i=1,\dots,n$ and
$\ve_i=-1$ for $i=n+1,\dots,2n$
in the symplectic case and $\ve_i=1$ for other cases.

The elements $F_{ij}$ satisfy the relations
\begin{align}\label{fijfji}
F_{ij}+\theta_{ij}F_{j\pr i\pr}=0,
\end{align}
for any $1\leq i,j\leq N$ and
\begin{align}\label{fijfkl}
[F_{ij},F_{kl}]=\delta_{kj}F_{il}-\delta_{il}F_{kj}-\delta_{k i\pr}\theta_{ij}F_{j\pr l}+\delta_{l j\pr}\theta_{ij}F_{k i\pr},
\end{align}
where $\theta_{ij}=\ve_i\tss\ve_j$, so $\theta_{ij}=1$ in the orthogonal case, and
$\theta_{ij}=\ve_i\tss\ve_j$ in the symplectic case.
%To consider the three cases $B$, $C$ and $D$
%simultaneously we will use
%the notation $\g_N$ for any of the Lie algebras $\oa_N$
%or $\spa_N$.

%The simple roots of $\g_N$ can be uniformly defined as follows. Let $\ep_1,\dots,\ep_n$ be an orthonormal basis of the $n$-dimensional Euclidian space
%with the bilinear form $(.\ts,.)$, then the simple roots of
%$\g_N$ are $\al_1,\dots,\al_n$ with
%\beql{ali}
%\al_i=\ep_i-\ep_{i+1},\qquad i=1,\dots,n-1,
%\eeq
%and
%\beql{aln}
%\al_n=\begin{cases}
%\ep_n
%\qquad&\text{for}\quad \oa_{2n+1}\\
%2\tss \ep_n
%\qquad&\text{for}\quad \spa_{2n}\\
%\ep_{n-1}+\ep_n
%\qquad&\text{for}\quad \oa_{2n}
%\end{cases}.
%\eeq

Set $\Ac=\mathbb{C}[[\hbar]]$. To introduce the $R$-matrix presentation of the Yangian double, %we will use the standard
%tensor notation.
let $e_1,\dots,e_N$ be the canonical basis of $\CC^N$ and $e_{ij}$ the matrix units in
$\End\CC^N$. Following the notation of \cite{amr:otr}, for any element
%is identified with the algebra of $N\times N$ matrices with the matrix units
%$e_{ij}$ to denote the matrix unit of $\End\CC^{N}$. % with $i,j\in\{1,\dots,N\}$ form a basis of $\End\CC^N$.
%Consider the tensor product algebra of the form
%\beql{tenprkea}
%\End(\CC^{N})^{\ot m}\ot\Ac=\underbrace{\End\CC^{N}\ot\dots\ot\End\CC^{N}}_m{}\ot\Ac,
%\eeq
%where $\Ac$ is a unital associative algebra.
\beql{mata}
X=\sum_{i,j=1}^N e_{ij}\ot X_{ij}\in \End\CC^N\ot \Ac
\eeq
and $a\in\{1,\dots,m\}$, we associate $X_a$ the following element%the following element
%\eqref{mata} associated with the $a$-th copy of $\End\CC^{N}$, namely
\beql{xa}
X_a=\sum_{i,j=1}^N 1^{\ot(a-1)}\ot e_{ij}\ot 1^{\ot(m-a)}\ot X_{ij}
\in \End(\CC^{N})^{\ot m}\ot\Ac,
\eeq
where $1$
is the identity endomorphism. One also defines $C_{ab}\in \End(\CC^{N})^{\ot m}$ for an operator $C\in \End \CC^N\ot\End \CC^N$ in the same manner.
%In the same spirit, for any operator
%\ben
%C=\sum_{i,j,k,l=1}^N c^{}_{ijkl}\ts e_{ij}\ot e_{kl}\in
%\End \CC^N\ot\End \CC^N
%\een
%and for each pair of indices $1\leq a<b\leq m$,
%we denote by $C_{ab}$ the following endomorphism
%\beql{ars}
%C_{ab}=\sum_{i,j,k,l=1}^N c^{}_{ijkl}\ts
%1^{\ot(a-1)}\ot e_{ij}\ot 1^{\ot(b-a-1)}\ot e_{kl}\ot 1^{\ot(m-b)}\in\End(\CC^{N})^{\ot m}.
%\eeq
%In particular, one use $C_{ab}$ for the element $C_{ab}\ot1$ in the algebra
%\eqref{tenprkea}.

Set $\ve=1$ (or $-1$) according to the orthogonal (or symplectic) case. Let $\kappa= N/2-\ve$.
%\ben
%\kappa=\begin{cases} N/2-1&\qquad\text{in the orthogonal case,}\\[0.3em]
%N/2+1&\qquad\text{in the symplectic case}.
%\end{cases}
%\een
%As defined in ,
The $R$-{\em matrix} $R(u)$ associated with the defining representation of $\mathfrak g_N$ is the rational function in $u$
with values in the tensor product algebra
$\End\CC^N\ot\End\CC^N$ given by \cite{zz:rf}
\beql{zamolr}
R(u)=1-\frac{\hbar P}{u}+\frac{\hbar Q}{u-\kappa\hbar},
\eeq
where
\beql{p}
P=\sum_{i,j=1}^N e_{ij}\ot e_{ji}, \qquad Q=\sum_{i,j=1}^N \ve_i\tss\ve_j\ts e_{ij}
\ot e_{i'j'}
\eeq

%in the orthogonal and symplectic case, respectively. Note the relations
It is easy to see that
$
P^2=1,\ Q^2=N\ts Q
$
and
\ben
PQ=QP=\ve Q
%\begin{cases} \phantom{-}Q&\qquad\text{in the orthogonal case,}\\
%-Q&\qquad\text{in the symplectic case}.
%\end{cases}
\een
The $R$-{\em matrix} \eqref{zamolr} satisfies the Yang--Baxter equation
\beql{yberep}
R_{12}(u-v)\ts R_{13}(u)\ts R_{23}(v)
=R_{23}(v)\ts R_{13}(u)\ts R_{12}(u-v).
\eeq

We now define the extended Yangian double $\DX_\hbar^{R}(\mathfrak{g}_N)$ and the Yangian double $\DY_\hbar^{R}(\mathfrak{g}_N)$
corresponding to the R-matrix $R(u)$.
\bde The algebra
$\DX_\hbar^{R}(\mathfrak{g}_N)$ is an associative algebra over $\Ac$ topologically generated by $t_{ij}^{(r)}$ and $c$ where $1\leq i,j\leq N$, $r\in \mathbb{Z}^{\times}$
subject to the relations written in terms of
generating series:
\beq
[T^{\pm}(u),c]=0,
\eeq
\beq\label{YDrelation1}
R(u-v)T^{\pm}_{1}(u)T^{\pm}_2(v)=T^{\pm}_2(v)T^{\pm}_{1}(u)R(u-v),
\eeq
\beq\label{relationpm}
R(u_{+}-v_{-})T_{1}^{+}(u)T_{2}^{-}(v)=T_{2}^{-}(v)T_{1}^{+}(u)R(u_{-}-v_{+}),
\eeq
where $u_{\pm}=u\pm \frac{1}{4}\hbar c$, the generating series $T^{\pm}(u)=\sum_{i,j=1}^{N}t^{\pm}_{ij}(u)\otimes E_{ij}$ are defined by
%\beq
%T^{\pm}(u)=\sum_{i,j=1}^{N}t^{\pm}_{ij}(u)\otimes E_{ij},
%\eeq
%%with
\beql{tiju}
t^{+}_{ij}(u)=\delta_{ij}+\hbar\sum_{r=1}^{\infty}t^{(r)}_{ij}u^{-r},\quad t^{-}_{ij}(u)=\delta_{ij}-\hbar\sum_{r=1}^{\infty}t^{(-r)}_{ij}u^{r-1}.
\eeq
and the $R$-matrices are expanded as power series in $u^{-1}$, for instance, $\frac1{u_+-v_-}=\sum_{k=0}^{\infty}u^{-1}(\frac{v-\hbar c/2}u)^k$.
Note that the matrix elements $t^{\pm}_{ii}(u)$ and $T^{\pm}(u)$ are defined over $\mathbb C[[\hbar]]$ via the $\hbar$-adic topology. In fact, they are invertible elements
as their constant terms are invertible as formal power series in $\hbar$.
%Then the generating relations can be written as the following form:
\ede

In terms of the R-matrix %\eqref{tiju}
the defining relations \eqref{YDrelation1} and \eqref{relationpm}
can be written explicitly as
\begin{align}\label{defrel++}
[\tss t^{\pm}_{ij}(u),t^{\pm}_{kl}(v)]&=\frac{\hbar}{u-v}
\Big(t^{\pm}_{kj}(u)\ts t^{\pm}_{il}(v)-t^{\pm}_{kj}(v)\ts t^{\pm}_{il}(u)\Big)\\
{}&-\frac{\hbar}{u-v-\kappa\hbar}
\Big(\de_{k i\pr}\sum_{p=1}^N\theta_{ip}\ts t^{\pm}_{pj}(u)\ts t^{\pm}_{p'l}(v)-
\de_{l j\pr}\sum_{p=1}^N\theta_{jp}\ts t^{\pm}_{k\tss p'}(v)\ts t^{\pm}_{ip}(u)\Big),
\non
\end{align}
\begin{align}\label{defrel+-}
&[\tss t^{+}_{ij}(u),t^{-}_{kl}(v)]=\frac{\hbar}{u_{+}-v_{-}}
t^{+}_{kj}(u)\ts t^{-}_{il}(v)-\frac{\hbar}{u_{-}-v_{+}}t^{-}_{kj}(v)\ts t^{+}_{il}(u)\\
{}&-\frac{\hbar}{u_{+}-v_{-}-\kappa\hbar}
\de_{k i\pr}\sum_{p=1}^N\theta_{ip}\ts t^{+}_{pj}(u)\ts t^{-}_{p'l}(v)+\frac{\hbar}{u_{-}-v_{+}-\kappa\hbar}
\de_{l j\pr}\sum_{p=1}^N\theta_{jp}\ts t^{-}_{k\tss p'}(v)\ts t^{+}_{ip}(u),
\non
\end{align}
where the factors $\frac1{u-v}, \frac1{u_+-v_-}$ etc are understood as infinite series expanded in powers $u^{-1}$, e.g.
$\frac1{u_+-v_--\kappa\hbar}=\sum_{k=0}^{\infty}(v-\hbar c/2+\kappa\hbar)^ku^{-k-1}$. We remark that the formal series identities
\eqref{defrel++}-\eqref{defrel+-} are understood as follows. First, if $A(u)=\sum_{r=0}^{\infty}a_ru^{-r}\in
\mathfrak[[u^{-1}]]$ and
$B(v)=\sum_{r=0}^{\infty}b_rv^{r}\in \mathfrak a[[v]]$ over an algebra $\mathfrak a$, then $A(u)B(v)=\sum_{r, s=0}^{\infty}a_rb_su^{-r}v^{s}\in \mathfrak a[[u^{-1}, v]]$ as a formal power series.
%which is also viewed as an element of $\sum_{r=0}^{\infty} a_rB(v)u^{-r}\in \mathbb C[[v, u^{-1}]]$.
Moreover, suppose $f(u)\in \mathfrak a[[u, u^{-1}]]$ over an algebra $\mathfrak a$ and that $S$ is a multiplicative nonzero subset of the ring $\mathfrak a[[u, u^{-1}]]$, the localization $S^{-1}\mathfrak a[[u, u^{-1}]]$ is defined in the usual manner. In particular, the function
$$\frac1{u-v-c}=\frac{u^{-1}v^{-1}}{v^{-1}-u^{-1}-cu^{-1}v^{-1}}$$
is viewed inside the extension $\widehat{\mathfrak a}[[u^{-1}, v^{-1}]]:=S^{-1}\mathfrak a[[u^{-1}, v^{-1}]]$ for $S=\{(u^{-1}-v^{-1}+cu^{-1}v^{-1})^i|i\in\mathbb N\}$. For example, \eqref{defrel+-} is understood as an identity over
the extension $\widehat{\DX}^R(\mathfrak g_N)[[u^{-1}, v]]$. Also, we may need to work over
$\mathbb{C}[[\hbar]]$ with the help of the $\hbar$-adic topology as remarked at the end of \eqref{tiju}.
In the sequel similar identities will always be treated likewise in an appropriate way.

Let $f(u)\in \Ac[[u^{-1}]]$ and $g(u)\in \Ac[[u]]$ be any formal series with constant term $1$. It is easy to see that the map
\begin{align}\label{auto}
\mu_{f,g}:\quad T^{+}(u)&\mapsto f(u)T^{+}(u), \quad T^{-}(u)\mapsto g(u)T^{-}(u)  \\
c&\mapsto c \non
\end{align}
defines an algebra homomorphism of $\DX_\hbar^{R}(\mathfrak{g}_N)$.

The Yangian double $\DY_\hbar^{R}(\mathfrak{g}_N)$ is defined as the subalgebra of $\DX_\hbar^{R}(\mathfrak{g}_N)$ consisting of the elements stable under all homomorphisms of the form \eqref{auto}.

\bpr
There exists formal series $z^{\pm}_N(u)\in \DX_\hbar^R(\mathfrak{g}_N)[[u^{\mp 1}]]$ (over $\mathbb C[[\hbar]]$) such that
\beql{zcenter}
T^{\pm}(u)T^{\pm}(v)^{t}\mid_{v=u+\kappa\hbar}=T^{\pm}(v)^{t}T^{\pm}(u)\mid_{v=u+\kappa\hbar}=z^{\pm}_N(u)I,
\eeq
%where
\beql{zn}
z^{+}_N(u)=1+\hbar\sum_{ r=1}^{\infty} z^{(r)}_N\ts u^{-r}, \quad z^{-}_N(u)=1-\hbar\sum_{r=1}^{\infty} z^{(-r)}_N\ts u^{r-1},
\eeq
and the symbol $t$ denotes the matrix transposition defined by
\beql{transp}
(X^{\tss t})_{ij}=\begin{cases} X_{j'i'}\qquad&\text{in the orthogonal case,}\\
\ve_i\tss\ve_j\ts X_{j'i'}\qquad&\text{in the symplectic case.}
\end{cases}
\eeq
We will write \eqref{zcenter} as $T^{\pm}(u)T^{\pm}(u+\kappa\hbar)^{t}=T^{\pm}(u+\kappa\hbar)^{t}T^{\pm}(u)=z^{\pm}_N(u)I$ for simplicity.
\epr

\bpf
The result is proved by applying the property of $Q$, see \cite[Sec.~2]{amr:otr} for more details.
\epf

Let  $\ZDX(\g_N)$ be the subalgebra of $\DX_\hbar^R(\mathfrak{g}_N)$ topologically generated by all the elements $z_N^{(r)}$ with $r\in \mathbb{Z}^{\times}$.
%We have the following results similar to those in \cite{amr:otr}.
Similar to \cite[Sec.~3]{amr:otr}, the following result is seen by normalizing the
generating series $t_{ij}^{\pm}(u)$ using \eqref{zcenter}.

\bthm\label{defDY}
One has the tensor product decomposition
\begin{align}\label{tensordecom}
\DX_\hbar^R(\mathfrak{g}_N)=\ZDX(\g_N)\widetilde{\otimes} \DY_\hbar^R(\mathfrak{g}_N),
\end{align}
where $\ZDX(\g_N)\widetilde{\otimes} \DY_\hbar^R(\mathfrak{g}_N)$ $=\displaystyle\lim_{\longleftarrow\atop n}(\ZDX(\g_N)\otimes_{\Ac} \DY_\hbar^R(\mathfrak{g}_N)/(\hbar^n))$ is a topological tensor product (see \cite[p.~390]{CKa} for the definition) over $\Ac$. Moreover, the Yangian double $\DY_\hbar^R(\mathfrak{g}_N)$ is isomorphic to the quotient of $\DX_\hbar^R(\mathfrak{g}_N)$ by the ideal generated by the coefficients $z^{(r)}_N$ of the series $z^{\pm}_N(u)$.
\ethm

Next we introduce some notations for later discussion. By the Gauss decomposition
the matrix $T^{\pm}(u)$ can be factored as:
\beql{gd}
T^{\pm}(u)=F^{\pm}(u)\ts H^{\pm}(u)\ts E^{\pm}(u),
\eeq
where $F^{\pm}(u)$, $H^{\pm}(u)$ and $E^{\pm}(u)$ are uniquely determined matrices of the form
\ben
F^{\pm}(u)=\begin{bmatrix}
1&0&\dots&0\ts\\
f^{\pm}_{21}(u)&1&\dots&0\\
\vdots&\vdots&\ddots&\vdots\\
f^{\pm}_{N1}(u)&f^{\pm}_{N2}(u)&\dots&1
\end{bmatrix},
\qquad
E^{\pm}(u)=\begin{bmatrix}
\ts1&e^{\pm}_{12}(u)&\dots&e^{\pm}_{1N}(u)\ts\\
\ts0&1&\dots&e^{\pm}_{2N}(u)\\
\vdots&\vdots&\ddots&\vdots\\
0&0&\dots&1
\end{bmatrix},
\een
and $H^{\pm}(u)=\diag\ts\big[k^{\pm}_1(u),\dots,k^{\pm}_N(u)\big]$. The coefficients of the matrix elements are referred
as {\em Gaussian generators}, explicitly
%the Introduce the coefficients of the series $e^{\pm}_{ij}(u),f^{\pm}_{ji}(u)$ and $k^{\pm}_{i}(u)$ by
\ben
e^{+}_{ij}(u)=\hbar\sum_{r=1}^{\infty} e_{ij}^{(r)}\tss u^{-r},\qquad
f^{+}_{ji}(u)=\hbar\sum_{r=1}^{\infty} f_{ji}^{(r)}\tss u^{-r},\qquad
k^{+}_i(u)=1+\hbar\sum_{r=1}^{\infty} k_i^{(r)}\tss u^{-r},\qquad
\een
\ben
e^{-}_{ij}(u)=-\hbar\sum_{r=1}^{\infty} e_{ij}^{(-r)}\tss u^{r-1},\qquad
f^{-}_{ji}(u)=-\hbar\sum_{r=1}^{\infty} f_{ji}^{(-r)}\tss u^{r-1},\qquad
k^{-}_i(u)=1-\hbar\sum_{r=1}^{\infty} k_i^{(-r)}\tss u^{r-1}.
\een
By definition, the algebra $\DX_\hbar^R(\mathfrak{g}_N)$
is generated by the elements $e^{(r)}_{ij}$, $f^{(r)}_{ji}$ and $k^{(\pm r)}_{i}$. %Thus, we will refer them as the {\em Gaussian generators}.
We further define the series $H^{\pm}_i(u)$ with coefficients in $\DX_\hbar^R(\mathfrak{g}_N)$ by
\ben
H^{\pm}_i(u)=k^{\pm}_i\big(u-(i-1)\hbar/2\big)^{-1}\ts k^{\pm}_{i+1}\big(u-(i-1)\hbar/2\big)
\een
for $i=1,\dots,n-1$, and
\ben
H^{\pm}_n(u)=\begin{cases}
k^{\pm}_n\big(u-(n-1)\hbar/2\big)^{-1}\ts k^{\pm}_{n+1}\big(u-(n-1)\hbar/2\big)
\qquad&\text{for}\quad \oa_{2n+1} \\[0.6em]
2\tss k^{\pm}_n\big(u-n\hbar/2\big)^{-1}\ts k^{\pm}_{n+1}\big(u-n\hbar/2\big)
\qquad&\text{for}\quad \spa_{2n}\\[0.6em]
k^{\pm}_{n-1}\big(u-(n-2)\hbar/2\big)^{-1}\ts k^{\pm}_{n+1}\big(u-(n-2)\hbar/2\big)
\qquad&\text{for}\quad \oa_{2n}.
\end{cases}
\een
Furthermore, we set
\ben
X_i^+(u)=e_{i}^{+}\big(u_{+}\big)-e_{i}^{-}\big(u_{-}\big),\qquad
X_i^-(u)=f_{i}^{+}\big(u_{-}\big)-f_{i}^{-}\big(u_{+}\big)
\een
where
\ben
e_i^{\pm}(u)=e_{i\ts i+1}^{\pm}(u),\qquad
f_i^{\pm}(u)=f_{i+1\ts i}^{\pm}(u)
\een
for $i=1,\dots,n-1$,
\ben
e_n^{\pm}(u)=\begin{cases}
e_{n\ts n+1}^{\pm}(u)\\
%\qquad&\text{for}\quad \oa_{2n+1} \quad \mbox{or}\quad \spa_{2n}\\[0.3em]
%e_{n\ts n+1}^{\pm}(u)
%\qquad&\text{for}\quad \spa_{2n}\\[0.3em]
e_{n-1\ts n+1}^{\pm}(u)\\
%\qquad&\text{for}\quad \oa_{2n}
\end{cases},\quad
f_n^{\pm}(u)=\begin{cases}
f_{n+1\ts n}^{\pm}(u)
\qquad&\text{for}\quad \oa_{2n+1}/\spa_{2n}\\[0.3em]
%f_{n+1\ts n}^{\pm}(u)
%\qquad&\text{for}\quad \spa_{2n}\\[0.3em]
f_{n+1\ts n-1}^{\pm}(u)
\qquad&\text{for}\quad \oa_{2n}
\end{cases}.
\een
%and
%\ben
%f_n^{\pm}(u)=\begin{cases}
%f_{n+1\ts n}^{\pm}(u)
%\qquad&\text{for}\quad \oa_{2n+1} \quad \mbox{or} \quad \spa_{2n}\\[0.3em]
%%f_{n+1\ts n}^{\pm}(u)
%%\qquad&\text{for}\quad \spa_{2n}\\[0.3em]
%f_{n+1\ts n-1}^{\pm}(u)
%\qquad&\text{for}\quad \oa_{2n}
%\end{cases}
%\een

Finally, we set
\ben
E_{i}(u)=\frac{1}{\hbar}X_i^+\big(u-(i-1)\hbar/2\big),\qquad
F_{i}(u)=\frac{1}{\hbar}X_i^-\big(u-(i-1)\hbar/2\big)
\een
for $i=1,\dots,n-1$,
\ben
E_{n}(u)=\begin{cases}
\frac{1}{\hbar}X_n^+\big(u-(n-1)\hbar/2\big)\\
%\qquad&\text{for}\quad \oa_{2n+1}\\[0.3em]
\frac{1}{\hbar}X_n^+\big(u-n\hbar/2\big)\\
%\qquad&\text{for}\quad \spa_{2n}\\[0.3em]
\frac{1}{\hbar}X_n^+\big(u-(n-2)\hbar/2\big)
%\qquad&\text{for}\quad \oa_{2n}
\end{cases},
F_{n}(u)=\begin{cases}
\frac{1}{\hbar}X_n^-\big(u-(n-1)\hbar/2\big)
\qquad&\text{for}\quad \oa_{2n+1}\\[0.3em]
\frac{1}{\hbar}X_n^-\big(u-n\hbar/2\big)
\qquad&\text{for}\quad \spa_{2n}\\[0.3em]
\frac{1}{\hbar}X_n^-\big(u-(n-2)\hbar/2\big)
\qquad&\text{for}\quad \oa_{2n}
\end{cases}.
\een
%and
%\ben
%F_{n}(u)=\begin{cases}
%X_n^-\big(u-(n-1)/2\big)
%\qquad&\text{for}\quad \oa_{2n+1}\\[0.3em]
%X_n^-\big(u-n/2\big)
%\qquad&\text{for}\quad \spa_{2n}\\[0.3em]
%X_n^-\big(u-(n-2)/2\big)
%\qquad&\text{for}\quad \oa_{2n}
%\end{cases}
%\een

Now we have
explicit formulas for the series $z^{\pm}_N(u)$
in terms of the Gaussian generators $k^{\pm}_i(u)$ (see \cite{ji:iso}).

\bthm\label{thm:Center}
In $\DX_\hbar^R(\mathfrak{g}_N)$ one has respectively for $\g_N=\oa_{2n+1}$ and $\g_N=\oa_{2n}/\spa_{2n}$:
\ben
z^{\pm}_N(u)=\begin{cases}\prod_{i=1}^n k^{\pm}_i(u+\ka\hbar-i\hbar)^{-1}\tss \prod_{i=1}^n k^{\pm}_i(u+\ka\hbar-i\hbar+\hbar)\cdot
k^{\pm}_{n+1}(u)\tss k^{\pm}_{n+1}\big(u-\hbar/2\big),\\
\prod_{i=1}^{n-1} k^{\pm}_i(u+\ka\hbar-i\hbar)^{-1}\tss \prod_{i=1}^n k^{\pm}_i(u+\ka\hbar-i\hbar+\hbar)\cdot
k^{\pm}_{n+1}(u).\end{cases}
\een
%for $\g_N=\oa_{2n+1}$,
%\ben
%z^{\pm}_N(u)=\prod_{i=1}^{n-1} k^{\pm}_i(u+\ka-i)^{-1}\tss \prod_{i=1}^n k^{\pm}_i(u+\ka-i+1)\cdot
%k^{\pm}_{n+1}(u)
%\een
%for $\g_N=\oa_{2n}$ and $\g_N=\spa_{2n}$.
\ethm

\section{PBW theorems}
\label{sec:PBW}

In this section, we will prove the Poincar\'e-Birkhoff-Witt theorem for the algebras $\DX_\hbar^R(\mathfrak{g}_N)$ and $\DY_\hbar^R(\mathfrak{g}_N)$. We usually omit the superscript $R$ and denote $\DY_\hbar^R(\mathfrak{g}_N)$~(resp. $\DX_\hbar^R(\mathfrak{g}_N)$) simply by $\DY_\hbar(\g_{N})$~(resp. $\DX_\hbar(\g_{N})$).
%Firstly, we need a lemma which is implied by the proof of Theorem~\ref{defDY}.
The following lemma is easily obtained following the method of \cite{amr:otr}.

\ble\label{ZDXiso}
The elements $z^{(r)}_N$ (resp. $z^{(-r)}_N$) $(r\geq 1)$ are algebraically independent over $\Ac$, and the subalgebra $\ZDX^{+}(\g_N)$ (resp. $\ZDX^{-}(\g_N)$) topologically generated by $z^{(r)}_N$ (resp. $z^{(-r)}_N$) $(r\geq 1)$ is isomorphic to the algebra of polynomials in countably many variables over $\Ac$.
%Similarly, the subalgebra $\ZDX^{-}(\g_N)$ generated by $z^{(-r)}_N~(r\geq 1)$ is also isomorphic to the algebra of polynomials in countably many variables.
\ele

By Theorem~\ref{defDY}, the Yangian double $\DY_\hbar(\g_{N})$ is generated by the elements $\tau_{ij}^{(r)}$ and $c$, $1\leq i,j\leq N$, $r\in\mathbb{Z}^{\times}$
subject to the relations
\begin{align}\label{defrelDY++}
[\tss \tau^{\pm}_{ij}(u),\tau^{\pm}_{kl}(v)]&=\frac{\hbar}{u-v}
\Big(\tau^{\pm}_{kj}(u)\ts \tau^{\pm}_{il}(v)-\tau^{\pm}_{kj}(v)\ts \tau^{\pm}_{il}(u)\Big)\\
{}&-\frac{\hbar}{u-v-\kappa\hbar}
\Big(\de_{k i\pr}\sum_{p=1}^N\theta_{ip}\ts \tau^{\pm}_{pj}(u)\ts \tau^{\pm}_{p'l}(v)-
\de_{l j\pr}\sum_{p=1}^N\theta_{jp}\ts \tau^{\pm}_{k\tss p'}(v)\ts \tau^{\pm}_{ip}(u)\Big),
\non
\end{align}
\begin{align}\label{defrelDY+-}
&[\tss \tau^{+}_{ij}(u),\tau^{-}_{kl}(v)]=\frac{\hbar}{u_{+}-v_{-}}
\tau^{+}_{kj}(u)\ts \tau^{-}_{il}(v)-\frac{\hbar}{u_{-}-v_{+}}\tau^{-}_{kj}(v)\ts \tau^{+}_{il}(u)\\
{}&-\frac{\hbar}{u_{+}-v_{-}-\kappa\hbar}
\de_{k i\pr}\sum_{p=1}^N\theta_{ip}\ts \tau^{+}_{pj}(u)\ts \tau^{-}_{p'l}(v)+\frac{\hbar}{u_{-}-v_{+}-\kappa\hbar}
\de_{l j\pr}\sum_{p=1}^N\theta_{jp}\ts \tau^{-}_{k\tss p'}(v)\ts \tau^{+}_{ip}(u),
\non
\end{align}
and
\begin{align}\label{unitaryDY}
\sum_{p=1}^N\theta_{kp}\ts \tau^{\pm}_{p'\tss k'}(u+\kappa\hbar)\tau^{\pm}_{pl}(u)=\delta_{kl}.
\end{align}
where
$\tau^{+}_{ij}(u)=\delta_{ij}+\hbar\sum_{r=1}^{\infty}\tau^{(r)}_{ij}u^{-r}$ and $\tau^{-}_{ij}(u)=\delta_{ij}-\hbar\sum_{r=1}^{\infty}\tau^{(-r)}_{ij}u^{r-1}.$

Define a natural ascending filtration on the Yangian double $\DX_\hbar(\g_{N})$ by setting $\deg t_{ij}^{(r)}=r-1,$ $\deg t_{ij}^{(-r)}=-r,$ for all $r\geqslant 1$ and $\deg c=\deg \hbar=0$. Denote by $\bar t_{ij}^{\ts(\pm r)}$ the image of $t_{ij}^{(\pm r)}$ in the $(r-1)$ (or $(-r)$)-th component of the associated graded algebra $\gr\DX_\hbar(\g_{N})$. %, $\bar t_{ij}^{\ts(-r)}$ the image of $t_{ij}^{(-r)}$ in the $(-r)$-th component of the associated graded algebra $\gr\DX_\hbar(\g_{N})$.
The filtration on the Yangian double $\DY_\hbar(\g_{N})$ is induced by that of $\DX_\hbar(\g_{N})$. Denote by $\bar \tau_{ij}^{\ts(\pm r)}$ %and  $\bar \tau_{ij}^{\ts(-r)}$
the image of $\tau_{ij}^{\ts(\pm r)}$ in the $(r-1)$ (or $(-r)$)-th component of the associated graded algebra $\gr\DY_\hbar(\g_{N})$. One then has the following result for the graded algebras. %immediately from the defining relations.
%We have the following result which follows from the defining relations of $\DY_\hbar(\g_{N})$.

\bpr\label{gradedalg}
The mapping
\beq
F_{ij}\tss x^{r-1}\mapsto \bar \tau_{ij}^{\ts(r)},\quad F_{ij}\tss x^{-r}\mapsto \bar \tau_{ij}^{\ts(-r)},\quad K\mapsto \bar c,\quad \hbar\mapsto\bar{\hbar}
\eeq
defines a homomorphism
\ben
\U(\g_N[x,x^{-1}]\oplus \mathbb{C}K)[[\hbar]]\rightarrow\gr\DY_\hbar(\g_{N}),
\een
where $K$ is the central element.
\epr

\bpf
It follows from \eqref{unitaryDY}  that $\bar \tau_{ij}^{\ts(r)}+\theta_{ij}\bar \tau_{j\pr i\pr}^{\ts(r)}=0$ for any $1\leq i,j\leq N$ and $r\in\mathbb{Z}^{\times}$. Using the expansion
\begin{align*}
\frac{1}{u-v}=u^{-1}+u^{-2}v+u^{-3}v^{2}+\cdots,
\end{align*}
and taking the coefficient at $u^{-r}v^{s}~(r,s\geq 1)$ and keeping the highest degree terms on both sides of the relation \eqref{defrelDY+-}
gives that
\beq
[\bar \tau_{ij}^{\ts(r)},\bar \tau_{kl}^{\ts(-s)}]=\begin{cases}
\delta_{kj}\bar \tau_{il}^{\ts(r-s-1)}-\delta_{il}\bar \tau_{kj}^{\ts(r-s-1)}-\delta_{k i\pr}\theta_{ij}\bar \tau_{j\pr l}^{\ts(r-s-1)}+\delta_{l j\pr}\theta_{ij}\bar \tau_{k i\pr}^{\ts(r-s-1)}, \quad r\leq s\\
\delta_{kj}\bar \tau_{il}^{\ts(r-s)}-\delta_{il}\bar \tau_{kj}^{\ts(r-s)}-\delta_{k i\pr}\theta_{ij}\bar \tau_{j\pr l}^{\ts(r-s)}+\delta_{l j\pr}\theta_{ij}\bar \tau_{k i\pr}^{\ts(r-s)}, \quad r>s
\end{cases}
\eeq
Similarly, the coefficients at $u^{-r}v^{-s}~(r,s\geq 1)$ and $u^{r}v^{s}~(r,s\geq 1)$ of \eqref{defrelDY++} imply that
\begin{align}
[\bar \tau_{ij}^{\ts(r)},\bar \tau_{kl}^{\ts(s)}]&=\delta_{kj}\bar \tau_{il}^{\ts(r+s-2)}-\delta_{il}\bar \tau_{kj}^{\ts(r+s-2)}-\delta_{k i\pr}\theta_{ij}\bar \tau_{j\pr l}^{\ts(r+s-2)}+\delta_{l j\pr}\theta_{ij}\bar \tau_{k i\pr}^{\ts(r+s-2)},\\
%Take the coefficients at  on both sides of the relation \eqref{defrelDY++}, keeping the highest degree terms, we have
[\bar \tau_{ij}^{\ts(-r)},\bar \tau_{kl}^{\ts(-s)}]&=\delta_{kj}\bar \tau_{il}^{\ts(-r-s)}-\delta_{il}\bar \tau_{kj}^{\ts(-r-s)}-\delta_{k i\pr}\theta_{ij}\bar \tau_{j\pr l}^{\ts(-r-s)}+\delta_{l j\pr}\theta_{ij}\bar \tau_{k i\pr}^{\ts(-r-s)}.
\end{align}
The relations \eqref{fijfkl} then give the final result.
\epf

Let $\rho:$
$F_{ij}\mapsto e_{ij}-\theta_{ij}e_{j\pr i\pr}$ be the vector representation of $\g_N$ on $\mathbb{C}^{N}$.
For any $c\in \mathbb{C}$ let $\rho_{c}$ be the evaluation representation of $\g_{N}[x,x^{-1}]$ defined by
$$\rho_{c}:F_{ij}x^{s}\mapsto c^{s}\rho(F_{ij}),\quad s\in \mathbb{Z}.$$
For any $c_{1},\cdots,c_{l}\in \mathbb{C}$ define the tensor product %of the $\rho_{c_i}$'s %evaluation representations of $\g_{N}[x,x^{-1}]$:
$\rho_{c_{1},\cdots,c_{l}}=\rho_{c_{1}}\otimes\cdots\otimes\rho_{c_{l}}.$ By the similar method as in \cite[Sec.~3]{amr:otr}, we can prove the following property about the kernel of
$\rho_{c_{1},\cdots,c_{l}}$.

\ble\label{kernel} As the complex parameters $c_{1},\cdots,c_{l}$ and integer $l\geq 0$ vary, the intersection
of all $\mathrm{ker}\rho_{c_{1},\cdots,c_{l}}$ in $\U(\g_N[x,x^{-1}])$ %of the kernels of all representations $\rho_{c_{1},\cdots,c_{l}}$
is trivial.
\ele

Now we can %in a position to
prove the Poincar\'e-Birkhoff-Witt theorem for the algebra $\DX_\hbar(\g_{N})$.
We order the generators of $\DX_\hbar(\g_{N})$ as follows. Define
$t_{ij}^{+}(u)\prec t_{kl}^{-}(u)$ for all $i,j,k,l$ and set $t_{ij}^{(r)}\prec t_{kl}^{(s)}$
(resp. $t_{ij}^{(-r)}\prec t_{kl}^{(-s)}$) when $(i,j,r)\prec (k,l,s)$ (resp. $(i,j,-r)\prec (k,l,-s)$)
 in the lexicographical order.
 %set $t_{ij}^{(-r)}\prec t_{kl}^{(-s)}$ when $(i,j,-r)\prec (k,l,-s)$ in the lexicographical order, where $r,s\geq 1$.
 This clearly defines a well-defined total ordering among the generators $t_{ij}^{\pm}, c$ of $\DX_\hbar(\g_{N})$~(with the central element $c$ included in the ordering in an arbitrary way).

\bthm\label{PBW}
The ordered monomials in the generators $t_{ij}^{(r)},t_{ij}^{(-r)}$ and $c$ form a topological basis of the algebra $\DX_\hbar(\g_{N})$.
\ethm

\bpf
It follows by an easy induction from the defining relations of $\DX_\hbar(\g_{N})$ that the span of the ordered monomials in the generators $t_{ij}^{(r)},t_{ij}^{(-r)}$ and $c$ is dense in $\DX_\hbar(\g_{N})$ in the $\hbar$-adit topology.

To show that the ordered monomials are linearly independent, we consider the level zero subalgebra $\DX_{0}(\g_{N})$ first, which is the quotient of $\DX_\hbar(\g_{N})$ by the ideal generated by $c$. For each nonzero $a\in \mathbb{C}$, we have the evaluation modules of $\DX_{0}(\g_{N})$ defined by
$$\pi_{a}:~\DX_{0}(\g_{N})\rightarrow \mathrm{End}\mathbb{C}^{N}[[\hbar]],\quad t_{ij}^{(r)}\mapsto a^{r-1}\rho(F_{ij}),~~t_{ij}^{(-r)}\mapsto -a^{-r}\rho(F_{ij}),~~\hbar\mapsto\hbar$$
for all $r\geq 1$. If there is a nontrivial linear combination of ordered monomials equal to zero, say $t_{i_{1}j_{1}}^{(s_{1})}\cdots t_{i_{p}j_{p}}^{(s_{p})}t_{k_{1}l_{1}}^{(-t_{1})}\cdots t_{k_{q}l_{q}}^{(-t_{q})}$ with certain complex coefficients $A_{i_{1}j_{1}\cdots i_{p}j_{p}k_{1}l_{1}\cdots k_{q}l_{q}}^{(s_{1}\cdots s_{p}-t_{1}\cdots -t_{q})}$ where the indices $s_{1},\cdots,s_{p},t_{1},\cdots,t_{q}\geq 1$ and the number $p,q\geq 0$ may vary. Consider the image of this combination under the representation $\pi_{a_{1}}\otimes\cdots\otimes\pi_{a_{l}}$, which depends on $a_{1},\cdots,a_{l}$ polynomially. Take the terms of this polynomial with maximal total degree in $a_{1},\cdots,a_{l}$. Let $A$ the sum of these terms. Therefore $A\in (\mathrm{End}\mathbb{C}^{N})^{\otimes l}[[\hbar]]$ coincides with the image of the sum $B$
$$\sum A_{i_{1}j_{1}\cdots i_{p}j_{p}k_{1}l_{1}\cdots k_{q}l_{q}}^{(s_{1}\cdots s_{p}-t_{1}\cdots -t_{q})}(F_{i_{1}j_{1}}x^{s_{1}-1})\cdots(F_{i_{p}j_{p}}x^{s_{p}-1})(F_{k_{1}l_{1}}x^{-t_{1}})\cdots(F_{k_{q}l_{q}}x^{-t_{q}})\in \U(\g_N[x,x^{-1}])[[\hbar]]$$
under the representation $\rho_{a_{1},\cdots,a_{l}}$. By the Poincar\'e-Birkhoff-Witt theorem for $\U(\g_N[x,x^{-1}])$, $B\neq 0$. Due to Lemma~\ref{kernel}, there exists $a_{1},\cdots,a_{l}\in \mathbb{C}$ such that $\rho_{a_{1},\cdots,a_{l}}(B)\neq 0$, i.e $A\neq 0$. But this linear combination equal to zero implies $A=0$ for all $a_{1},\cdots,a_{l}\in \mathbb{C}$. This is a contradiction. Hence, the ordered monomials are linearly independent for the level zero algebra $\DX_{0}(\g_{N})$.

Then the similar argument of \cite[Theorem 2.2]{ji:cen} implies that ordered monomials are linearly independent in $\DX_\hbar(\g_{N})$.
%we can apply the same argument as in \cite[Theorem 2.2]{ji:cen}.
\epf

The following version of the Poincar\'e-Birkhoff-Witt theorem for the algebra $\DX_\hbar(\g_{N})$ is also useful.

\bco\label{PBW2}
Given any total ordering on the set of elements $t_{ij}^{(r)},t_{ij}^{(-r)},z_{N}^{(r)},z_{N}^{(-r)}$ and $c$ with
$$i>j\pr,\quad r\geq 1,\quad in~the~orthogonal~case,$$
and
$$i\geq j\pr,\quad r\geq 1,\quad in~the~symplectic~case,$$
the ordered monomials in these elements form a topological basis of $\DX_\hbar(\g_{N})$.
\eco

%\bpf
%It follows from Theorem~\ref{PBW} and the identity
%$$\sum_{i=1}^{N}\theta_{ki}t^{\pm}_{i\pr k\pr}(u+\kappa)t^{\pm}_{ik}(u)=z^{\pm}_{N}(u).$$
%\epf

\bco\label{PBWDY}
Given any total ordering on the set of elements $\tau_{ij}^{(r)},\tau_{ij}^{(-r)}$ and $c$ with
$$i>j\pr,\quad r\geq 1,\quad in~the~orthogonal~case,$$
and
$$i\geq j\pr,\quad r\geq 1,\quad in~the~symplectic~case,$$
the ordered monomials in these elements form a topological basis of $\DY_\hbar(\g_{N})$.
\eco

\bpf
It follows from Corollary~\ref{PBW2} and the tensor product decomposition
\begin{align*}
\DX_\hbar(\mathfrak{g}_N)=\ZDX(\g_N)\widetilde{\otimes}\DY_\hbar(\mathfrak{g}_N).
\end{align*}
\epf

\bco\label{gradediso}
The mapping defined in Proposition~\ref{gradedalg} is an algebra isomorphism.
\eco

\bpf
The injectivity follows from Corollary~\ref{PBWDY}.
\epf

\bco\label{gradealgDX}
The mapping
\beql{isomgrpol}
F_{ij}\tss x^{r-1}\mapsto \bar t_{ij}^{\ts(r)}-\frac{1}{2}\delta_{ij}\bar z_{N}^{(r)}, \ze_r\mapsto \bar z_{N}^{(r)}, F_{ij}\tss x^{-r}\mapsto \bar t_{ij}^{\ts(-r)}-\frac{1}{2}\delta_{ij}\bar z_{N}^{(-r)}, \vs_r\mapsto \bar z_{N}^{(-r)}, K\mapsto \bar c, \hbar\mapsto\bar{\hbar}
\eeq
defines an algebra isomorphism $\psi$:
\ben
\U(\g_N[x,x^{-1}]\oplus \mathbb{C}K)[[\hbar]]\ot\Ac[\ze_1,\ze_2,\dots]\ot\Ac[\vs_1,\vs_2,\dots]\rightarrow \gr\DX_\hbar(\g_{N}),
\een
where $\Ac[\ze_1,\ze_2,\dots]$ is the algebra of polynomials in variables $\ze_i$ over $\Ac$, $\Ac[\vs_1,\vs_2,\dots]$ is the algebra of polynomials in variables $\vs_i$ over $\Ac$, $K$ is the central element.
\eco

\bpf
It follows from Lemma~\ref{ZDXiso}, Corollary~\ref{gradediso} and the tensor product decomposition
\begin{align*}
\DX_\hbar(\mathfrak{g}_N)=\ZDX(\g_N)\widetilde{\otimes}\DY_\hbar(\mathfrak{g}_N).
\end{align*}
\epf

\section{Embedding theorems}
\label{sec:et}

We first recall the quasideterminants of a matrix.
Let $A=[a_{ij}]$ be an $N\times N$ matrix, with the row vectors $r_i$ and column vectors $c_j$, over a ring with $1$. Denote by $A^{ij}$ the submatrix
of $A$ by deleting the $i$-th row
and $j$-th column.
If the submatrix
$A^{ij}$ is invertible, the $ij$-{\em th quasideterminant of} $A$ \cite{gr:dm} is defined as
\ben
|A|_{ij}=a_{ij}-r^{\tss j}_i(A^{ij})^{-1}\ts c^{\tss i}_j,
\een
where $r^{\tss j}_i$ is the row vector obtained from $r_i$ by deleting the element $a_{ij}$, and $c^{\tss i}_j$
is the column vector obtained from $c_j$ by deleting the element $a_{ij}$.
For example,
the four quasideterminants of a $2\times 2$ matrix $A$ are
\ben
\bal
&|A|_{11}=a^{}_{11}-a^{}_{12}\ts a_{22}^{-1}\ts a^{}_{21},\qquad
|A|_{12}=a^{}_{12}-a^{}_{11}\ts a_{21}^{-1}\ts a^{}_{22},\\
&|A|_{21}=a^{}_{21}-a^{}_{22}\ts a_{12}^{-1}\ts a^{}_{11},\qquad
|A|_{22}=a^{}_{22}-a^{}_{21}\ts a_{11}^{-1}\ts a^{}_{12}.
\eal
\een
We also denote the quasideterminant $|A|_{ij}$
by boxing the entry $a_{ij}$,
\ben
|A|_{ij}=\left|\begin{matrix}a_{11}&\dots&a_{1j}&\dots&a_{1N}\\
                                   &\dots&      &\dots&      \\
                             a_{i1}&\dots&\boxed{a_{ij}}&\dots&a_{iN}\\
                                   &\dots&      &\dots&      \\
                             a_{N1}&\dots&a_{Nj}&\dots&a_{NN}
                \end{matrix}\right|.
\een

Now suppose that $n\geqslant 1$ in type $B$, and $n\geqslant 2$ in types $C$ and $D$.
Consider the algebra $\DX_\hbar(\g_{N-2})$ and assume the indices of
the generators $t_{ij}^{(r)}$ satisfy
$2\leqslant i,j\leqslant 2\pr$ and $r=1,2,\dots,-1,-2,\dots$ (here
$i\pr=N-i+1$). Next we will give the first main result.

\bthm\label{thm:embed}
The following mapping gives an injective homomorphism from\\
$\DX_\hbar^{\tilde{R}}(\g_{N-2})$ into $\DX_\hbar^R(\g_N)$
\beql{embedgen}
t^{\pm}_{ij}(u)\mapsto s^{\pm}_{ij}(u)= \left|\begin{matrix}
t_{11}^{\pm}(u)&t_{1j}^{\pm}(u)\\
t_{i1}^{\pm}(u)&\boxed{t_{ij}^{\pm}(u)}
\end{matrix}\right|=t_{ij}^{\pm}(u)-t_{i1}^{\pm}(u)t_{11}^{\pm}(u)^{-1}t_{1j}^{\pm}(u),\qquad 2\leqslant i,j\leqslant 2\pr,
\eeq
where
$$\tilde{R}(u)=1-\frac{\hbar P}{u}+\frac{\hbar Q}{u-\kappa\hbar+\hbar}.$$
\ethm

Introduce power series
$\tau^{\pm\tss a_1 a_2}_{\tss b_1 b_2}(u)$ in $u^{\mp 1}$
with coefficients in $\DX_\hbar(\g_{N})$
by the formulas:
\beql{quamintau}
R_{12}(\hbar)\ts T^{\pm}_1(u)\ts T^{\pm}_2(u-\hbar)=\sum_{a_i,b_i}e_{a_1b_1}\ot e_{a_2b_2}
\ot \tau^{\pm\tss a_1 a_2}_{\tss b_1 b_2}(u).
\eeq
We have the following properties about the series $\tau^{\pm\tss a_1 a_2}_{\tss b_1 b_2}(u)$; see \cite{ji:iso}.

\ble\label{lem:skewsymm}\quad
(i)\quad
If $a_1\ne a'_2$ then $\tau^{\pm\tss a_1 a_2}_{\tss b_1 b_2}(u)
=-\tau^{\pm\tss a_2 a_1}_{\tss b_1 b_2}(u)$.

\medskip
\noindent
(ii)\quad
If $b_1\ne b'_2$ then $\tau^{\pm\tss a_1 a_2}_{\tss b_1 b_2}(u)
=-\tau^{\pm\tss a_1 a_2}_{\tss b_2 b_1}(u)$.
\ele

\bre\label{rem:unnecc}
The skew-symmetry properties still hold without the assumptions
in the symplectic case.
\qed
\ere

\ble\label{lem:toneone}
For any $2\leqslant i,j\leqslant 2\pr$ we have
\beql{wttfo}
s^{\pm}_{ij}(u)=t^{\pm}_{11}(u+\hbar)^{-1}\tss \tau^{\pm\tss 1 i}_{\tss 1 j}(u+\hbar).
\eeq
Moreover,
\beql{commtoo}
\big[t^{\pm}_{11}(u), \tau^{\pm\tss 1 i}_{\tss 1 j}(v)\big]=0,
\eeq
\beql{tpmspm}
\big[t^{\pm}_{11}(u), s^{\pm}_{\tss i j}(v)\big]=0,
\eeq
\beql{commtoopm}
\frac{u_{+}-v_{-}-\hbar}{u_{+}-v_{-}}t_{11}^{+}(u)\tau^{-\tss 1 i}_{\tss 1 j}(v)=\frac{u_{-}-v_{+}-\hbar}{u_{-}-v_{+}}\tau^{-\tss 1 i}_{\tss 1 j}(v)t_{11}^{+}(u),
\eeq
\beql{commtoomp}
\frac{u_{-}-v_{+}+\hbar}{u_{-}-v_{+}+2\hbar}t_{11}^{-}(u)\tau^{+\tss 1 i}_{\tss 1 j}(v)
=\frac{u_{+}-v_{-}+\hbar}{u_{+}-v_{-}+2\hbar}\tau^{+\tss 1 i}_{\tss 1 j}(v)t_{11}^{-}(u),
\eeq
\beql{tpmsmp}
[t^{+}_{11}(u),s^{-}_{ij}(v)]=0,
\eeq
\beql{tpsmmp}
\frac{(u_{-}-v_{+})^{2}}{(u_{-}-v_{+})^{2}-\hbar^{2}}t^{-}_{11}(u)s^{+}_{ij}(v)=\frac{(u_{+}-v_{-})^{2}}{(u_{+}-v_{-})^{2}-\hbar^{2}}s^{+}_{ij}(v)t^{-}_{11}(u).
\eeq
\ele

\bpf
\eqref{embedgen} implies that
\ben
t_{11}^{\pm}(u+\hbar)\ts s_{ij}^{\pm}(u)=t_{11}^{\pm}(u+\hbar)\tss
t_{ij}^{\pm}(u)-t_{11}^{\pm}(u+\hbar)\tss t_{i1}^{\pm}(u)\tss t_{11}^{\pm}(u)^{-1}\tss t_{1j}^{\pm}(u).
\een
However, $t_{11}^{\pm}(u+\hbar)\tss t_{i1}^{\pm}(u)=t_{i1}^{\pm}(u+\hbar)\tss t_{11}^{\pm}(u)$ by \eqref{defrel++},
hence
\ben
t_{11}^{\pm}(u+\hbar)\ts s_{ij}^{\pm}(u)=t_{11}^{\pm}(u+\hbar)\tss
t_{ij}^{\pm}(u)-t_{i1}^{\pm}(u+\hbar)\tss t_{1j}^{\pm}(u).
\een
The definition \eqref{quamintau} implies that this
equals $\tau^{\pm\tss 1 i}_{\tss 1 j}(u+\hbar)$ hence \eqref{wttfo} follows. Relation \eqref{commtoo} and \eqref{tpmspm} is already proved in \cite{ji:iso}.
Relation \eqref{commtoopm} follows easily from the defining relations of $\DX_\hbar(\g_N)$.
Consider the following equation:
\beq
\begin{aligned}
&R_{01}(u_{+}-v_{-})R_{02}(u_{+}-v_{-}+\hbar)T_0^{+}(u)R_{12}(\hbar)T_1^{-}(v)T_2^{-}(v-\hbar)\\
&=R_{12}(\hbar)T_1^{-}(v)T_2^{-}(v-\hbar)T_0^{+}(u)R_{02}(u_{-}-v_{+}+\hbar)R_{01}(u_{-}-v_{+})
\end{aligned}
\eeq
Furthermore, we have
\beq
\begin{aligned}
&\langle 1,1,i| R_{01}(u_{+}-v_{-})R_{02}(u_{+}-v_{-}+\hbar)T_0^{+}(u)R_{12}(\hbar)T_1^{-}(v)T_2^{-}(v-\hbar)|1,1,j\rangle\\
&=\langle 1,1,i|R_{12}(\hbar)T_1^{-}(v)T_2^{-}(v-\hbar)T_0^{+}(u)R_{02}(u_{-}-v_{+}+\hbar)R_{01}(u_{-}-v_{+})|1,1,j\rangle
\end{aligned}
\eeq
which is just \eqref{commtoopm}. We can prove \eqref{commtoomp} similarly by the following relation:
\begin{multline}
\frac{(u_{-}-v_{+})^{2}}{(u_{-}-v_{+})^{2}-\hbar^{2}}R(u_{-}-v_{+})T_1^{-}(u)T_2^{+}(v)\\
=\frac{(u_{+}-v_{-})^{2}}{(u_{+}-v_{-})^{2}-\hbar^{2}}T_2^{+}(v)T_1^{-}(u)R(u_{+}-v_{-}).
\end{multline}
\eqref{tpmsmp} follows from \eqref{wttfo} and \eqref{commtoopm}. The proof of \eqref{tpsmmp} is similar.
\epf

The following lemma is needed for the first main result.

\ble\label{lem:sylv}
Set
\ben
\Gamma^{\pm}(u)=\sum_{ij=2}^{2'} E_{ij}\otimes \tau^{\pm\tss 1 i}_{\tss 1 j}(u).
\een
The following relations hold:
\beq
[\Gamma^{\pm}(u),c]=0,
\eeq
\beq
R^{[1]}(u-v)\Gamma_1^{\pm}(u)\Gamma_2^{\pm}(v)=\Gamma_2^{\pm}(v)\Gamma_1^{\pm}(u)R^{[1]}(u-v),
\eeq
\beq\label{RTpTm}
f(u_{+}-v_{-})R^{[1]}(u_{+}-v_{-})\Gamma_{1}^{+}(u)\Gamma_{2}^{-}(v)=\Gamma_{2}^{-}(v)\Gamma_{1}^{+}(u)R^{[1]}(u_{-}-v_{+})f(u_{-}-v_{+}),
\eeq
here, $f(u)=\frac{u-2\hbar}{u-\hbar}$, $R^{[1]}(u)=1-\frac{\hbar P}{u}+\frac{\hbar Q}{u-\kappa\hbar+\hbar}$.
\ele

\bpf
We prove the relation \eqref{RTpTm} to show the method. %Consider the tensor product algebra \eqref{tenprkea} with $m=4$.
Repeatedly using the Yang-Baxter equation \eqref{yberep} and
the $RTT$ relation \eqref{YDrelation1}, \eqref{relationpm} implies immediately that
\begin{multline}
\label{longrel}
R_{23}(a_{+}-\hbar)\tss R_{13}(a_{+})\tss R_{24}(a_{+})\tss R_{14}(a_{+}+\hbar)\tss
R_{12}(\hbar)\tss T^{+}_1(u)\tss T^{+}_2(u-\hbar)
R_{34}(\hbar)\tss T^{-}_3(v)\tss T^{-}_4(v-\hbar)\\[0.4em]
{}=R_{34}(\hbar)\tss T^{-}_3(v)\tss T^{-}_4(v-\hbar)\tss R_{12}(\hbar)\tss T^{+}_1(u)\tss T^{+}_2(u-\hbar)\tss
R_{14}(a_{-}+\hbar)\tss R_{24}(a_{-})\tss R_{13}(a_{-})\tss R_{23}(a_{-}-\hbar),
\end{multline}
where $a_+=u_{+}-v_{-}$, $a_-=u_{-}-v_{+}$. Let $V$ be the subspace of $(\CC^N)^{\ot 4}$ spanned by the
basis vectors of the form $e_1\ot e_j\ot e_1\ot e_l$ with $j,l\in\{2,\dots,2\pr\}$. Apply the same calculation as in \cite{ji:iso}, we can conclude that the restriction of the operator on the right hand side of
\eqref{longrel} to the subspace $V$ coincides with the operator
\beql{oprhs}
\frac{a_{-}-2\hbar}{a_{-}-\hbar}\ts R_{34}(\hbar)\tss T_3^{-}(v)\tss T_4^{-}(v-\hbar)\tss
R_{12}(\hbar)\tss T_1^{+}(u)\tss T_2^{+}(u-\hbar)\tss \Big(1-\frac{\hbar P_{24}}{a_{-}}+\frac{\hbar Q_{24}}{a_{-}-\ka\hbar+\hbar}\Big).
\eeq
Similarly, the restriction of the operator on the left hand side of
\eqref{longrel} to the subspace $V$ coincides with the operator
\beql{oplhs}
\frac{a_{+}-2\hbar}{a_{+}-\hbar}\ts \Big(1-\frac{\hbar P_{24}}{a_{+}}+\frac{\hbar Q_{24}}{a_{+}-\ka\hbar+\hbar}\Big)\ts
R_{12}(\hbar)\tss T_1^{+}(u)\tss T_2^{+}(u-\hbar)\tss
R_{34}(\hbar)\tss T_3^{-}(v)\tss T_4^{-}(v-\hbar).
\eeq
Moreover, we have
\begin{multline}
R_{12}(\hbar)\tss T^{+}_1(u)\tss T^{+}_2(u-\hbar)\tss
R_{34}(\hbar)\tss T^{-}_3(v)\tss T^{-}_4(v-\hbar)\tss(e_1\ot e_j\ot e_1\ot e_l)\\
{}\equiv\sum_{c,d=1}^{1\pr}\tau^{+\tss 1 c}_{\tss 1 j}(u)\tss \tau^{-\tss 1 d}_{\tss 1 l}(v)
\tss(e_1\ot e_c\ot e_1\ot e_d),
\label{tautwo}
\end{multline}
where we only keep the basis vectors which can give a nonzero contribution
to the coefficient of $e_1\ot e_i\ot e_1\ot e_k$ with $i,k\in\{2,\dots,2\pr\}$ after the subsequent
application of the operator $\Big(1-\frac{\hbar P_{24}}{a_{+}}+\frac{\hbar Q_{24}}{a_{+}-\ka\hbar+\hbar}$\Big). By Lemma~\ref{lem:skewsymm}\tss$(i)$, $\tau^{+\tss 1 1}_{\tss 1 j}(u)=\tau^{-\tss 1 1}_{\tss 1 l}(v)=0,$ hence the values $c=1$ and $d=1$ can be excluded from the range of the summation indices. The definition of $P,Q$ implies that the values
$c=1\pr$ and $d=1\pr$ can also be excluded. Therefore, we can write
an operator equality
\ben
1-\frac{\hbar P_{24}}{a_{+}}+\frac{\hbar Q_{24}}{a_{+}-\ka\hbar+\hbar}=R_{24}(a_{+}),
\een
where $R(u)$ is the $R$-matrix
associated with the algebra $\DX_\hbar(\g_{N-2})$.
Similarly, we can prove the same property for the operator on the right hand side of
\eqref{longrel} with the use of Lemma~\ref{lem:skewsymm}\tss$(ii)$. So we have
\begin{multline}
\frac{a_{+}-2\hbar}{a_{+}-\hbar}\ts R_{24}(a_{+})\ts
R_{12}(\hbar)\tss T_1^{+}(u)\tss T_2^{+}(u-\hbar)\tss
R_{34}(\hbar)\tss T_3^{-}(v)\tss T_4^{-}(v-\hbar)\\
=\frac{a_{-}-2\hbar}{a_{-}-\hbar}\ts
R_{34}(\hbar)\tss T_3^{-}(v)\tss T_4^{-}(v-\hbar)\tss
R_{12}(\hbar)\tss T_1^{+}(u)\tss T_2^{+}(u-\hbar)\tss
R_{24}(a_{-}).
\end{multline}
By equating the matrix elements, this completes the proof of the equation \eqref{RTpTm}.
\epf

Using Lemma~\ref{lem:toneone} and \ref{lem:sylv}, we can easily verify that the mapping \eqref{embedgen} is an algebra homomorphism. Next we show that the homomorphism \eqref{embedgen} is injective. Define a filtration $(\DX_\hbar^R(\g_{N})_{m})_{m\in \mathbb{Z}}$ on the extended Yangian double $\DX_\hbar^R(\g_{N})$ by setting $\deg t_{ij}^{(r)}=r-1,$ $\deg t_{ij}^{(-r)}=-r$ for all $r\geqslant 1$ and $\deg c=0$. Let $a$ be any nonzero element in $\DX_\hbar^R(\g_{N})$, Theorem~\ref{PBW} implies that $a\notin\cap_{m\in \mathbb{Z}}\DX_\hbar^R(\g_{N})_{m}$. Hence, $\cap_{m\in \mathbb{Z}}\DX_\hbar^R(\g_{N})_{m}=\{0\}$, i.e, the filtration on $\DX_\hbar^R(\g_{N})$ is separated. Denote by $\bar t_{ij}^{\ts(r)}$ the image of $t_{ij}^{(r)}$ in the $(r-1)$-th component of the associated graded algebra $\gr\DX_\hbar^R(\g_{N})$, $\bar t_{ij}^{\ts(-r)}$ the image of $t_{ij}^{(-r)}$
in the $(-r)$-th component of the associated graded algebra $\gr\DX_\hbar^R(\g_{N})$. The map \eqref{embedgen} induces a
homomorphism of the associated graded algebras $\gr\DX_\hbar^{\tilde{R}}(\g_{N-2})\to \gr\DX_\hbar^R(\g_{N})$. By Corollary~\ref{gradealgDX}, the mapping
\beql{isomapping}
\bar t_{ij}^{\ts(r)}\mapsto F_{ij}\tss x^{r-1}+\frac12\ts\de_{ij}\ts\ze_r,\quad \bar t_{ij}^{\ts(-r)}\mapsto F_{ij}\tss x^{-r}+\frac12\ts\de_{ij}\ts\vs_r,\quad \bar c\mapsto K,\quad \bar{\hbar}\mapsto\hbar
\eeq
defines an isomorphism
\ben
\gr\DX_\hbar^R(\g_{N})\cong \U(\g_N[x,x^{-1}]\oplus \mathbb{C}K)[[\hbar]]\ot\Ac[\ze_1,\ze_2,\dots]\ot\Ac[\vs_1,\vs_2,\dots].
\een
The variables $\ze_r,\vs_r$ correspond to the images of the elements $z^{(r)}_N,~z^{(-r)}_N$
defined in \eqref{zn},
\beql{znim}
\bar z^{\ts(r)}_N\mapsto \ze_r,\quad \bar z^{\ts(-r)}_N\mapsto \vs_r.
\eeq
Hence the homomorphism $\gr\DX_\hbar^{\tilde{R}}(\g_{N-2})\to \gr\DX_\hbar^R(\g_{N})$
is injective. Since the filtration on $\DX_\hbar^{\tilde{R}}(\g_{N-2})$ is separated, the homomorphism \eqref{embedgen} is also injective.

Fix a positive integer $m$ such that
$m\leqslant n$ for type $B$ and $m\leqslant n-1$ for types $C$ and $D$.
Assume that the generators $t_{ij}^{(r)}$ of the algebra $\DX_\hbar(\g_{N-2m})$ are
labelled by the indices
$m+1\leqslant i,j\leqslant (m+1)\pr$ and $r\in \mathbb{Z}^{\times}$. We will give a generalization of
Theorem~\ref{thm:embed}.

\bthm\label{thm:red}
The mapping $\psi_{m}$:
\beql{redu}
t^{\pm}_{ij}(u)\mapsto \left|\begin{matrix}
t^{\pm}_{11}(u)&\dots&t^{\pm}_{1m}(u)&t^{\pm}_{1j}(u)\\
\dots&\dots&\dots&\dots\\
t^{\pm}_{m1}(u)&\dots&t^{\pm}_{mm}(u)&t^{\pm}_{mj}(u)\\
t^{\pm}_{i1}(u)&\dots&t^{\pm}_{im}(u)&\boxed{t^{\pm}_{ij}(u)}
\end{matrix}\right|,\qquad m+1\leqslant i,j\leqslant (m+1)\pr,
\eeq
is an injective homomorphism $\DX_\hbar^{\tilde{R}}(\g_{N-2m})\rightarrow \DX_\hbar^R(\g_N)$.
Here $\DX_\hbar^{\tilde{R}}(\g_{N-2m})$ is the algebra corresponding to $\tilde{R}(u)=1-\frac{\hbar P}{u}+\frac{\hbar Q}{u-\kappa\hbar+m\hbar}$.
\ethm

\bpf
It follows from Theorem~\ref{thm:embed} and the well-known identity:
\beq
\left|\begin{matrix}
t^{\pm}_{11}(u)&\dots&t^{\pm}_{1m}(u)&t^{\pm}_{1j}(u)\\
\dots&\dots&\dots&\dots\\
t^{\pm}_{m1}(u)&\dots&t^{\pm}_{mm}(u)&t^{\pm}_{mj}(u)\\
t^{\pm}_{i1}(u)&\dots&t^{\pm}_{im}(u)&\boxed{t^{\pm}_{ij}(u)}
\end{matrix}\right|=
\left|\begin{matrix}
s^{\pm}_{22}(u)&\dots&s^{\pm}_{2m}(u)&s^{\pm}_{2j}(u)\\
\dots&\dots&\dots&\dots\\
s^{\pm}_{m2}(u)&\dots&s^{\pm}_{mm}(u)&s^{\pm}_{mj}(u)\\
s^{\pm}_{i2}(u)&\dots&s^{\pm}_{im}(u)&\boxed{s^{\pm}_{ij}(u)}
\end{matrix}\right|.
\eeq
\epf

Now we set
\ben
F^{[m]\pm}(u)=\begin{bmatrix}
1&0&\dots&0\ts\\
f^{\pm}_{m+2\ts m+1}(u)&1&\dots&0\\
\vdots&\ddots&\ddots&\vdots\\
f^{\pm}_{(m+1)'\tss m+1}(u)&\dots&f^{\pm}_{(m+1)'\ts (m+2)'}(u)&1
\end{bmatrix},
\een
\ben
E^{[m]\pm}(u)=\begin{bmatrix} 1&e^{\pm}_{m+1\tss m+2}(u)&\ldots&e^{\pm}_{m+1\tss (m+1)'}(u)\\
                        0&1&\ddots &\vdots\\
                         \vdots&\vdots&\ddots&e^{\pm}_{(m+2)'\tss(m+1)'}(u)\\
                         0&0&\ldots&1\\
           \end{bmatrix}
\een
and $H^{[m]\pm}(u)=\diag\ts\big[k^{\pm}_{m+1}(u),\dots,k^{\pm}_{(m+1)'}(u)\big]$. Moreover, we define matrices $T^{[m]\pm}(u)$ as follows:
\ben
T^{[m]\pm}(u)=F^{[m]\pm}(u)\tss H^{[m]\pm}(u)\tss E^{[m]\pm}(u).
\een
And the entries of $T^{[m]\pm}(u)$ will be denoted by $t^{[m]\pm}_{ij}(u)$, here $m+1\leqslant i,j\leqslant (m+1)'$. Then we have a useful property about the series $t^{[m]\pm}_{ij}(u)$; see \cite[Section~4]{ji:iso}.

\bpr\label{prop:gauss-consist}
The series $t^{[m]\pm}_{ij}(u)$ coincides with the image of the generator series $t^{\pm}_{ij}(u)$ of the extended Yangian double $\DX_\hbar^{\tilde{R}}(\g_{N-2m})$ under the embedding
\eqref{redu},
\begin{align*}
t^{[m]\pm}_{ij}(u)=\psi_{m}(t^{\pm}_{ij}(u)),\quad m+1\leq i,j \leq (m+1)'.
\end{align*}
\epr

The following corollary is immediate from Proposition~\ref{prop:gauss-consist}.

\bco\label{cor:guass-embed}
The subalgebra $\DX_\hbar^{[m]}(\g_N)$ generated by the coefficients of all series $t^{[m]\pm}_{ij}(u)$ with $m+1\leqslant i,j\leqslant (m+1)'$ and $c$ is isomorphic to the extended Yangian
double $\DX_\hbar^{\tilde{R}}(\g_{N-2m})$. Here $\DX_\hbar^{\tilde{R}}(\g_{N-2m})$ is the algebra corresponding to $\tilde{R}(u)=1-\frac{\hbar P}{u}+\frac{\hbar Q}{u-\kappa\hbar+m\hbar}$. Furthermore, we have the relations
\begin{align}
\tilde{R}_{12}(u-v)T_{1}^{[m]\pm}(u)T_{2}^{[m]\pm}(v)=T_{2}^{[m]\pm}(v)T_{1}^{[m]\pm}(u)\tilde{R}_{12}(u-v),
\end{align}
\begin{align}
\tilde{R}_{12}(u_{+}-v_{-})T_{1}^{[m]+}(u)T_{2}^{[m]-}(v)
=T_{2}^{[m]-}(v)T_{1}^{[m]+}(u)\tilde{R}_{12}(u_{-}-v_{+}).
\end{align}
\eco

Suppose that $\{a_{1},\cdots,a_{k}\}$ and $\{b_{1},\cdots,b_{k}\}$ are subsets of $\{1,\cdots,N\}$ satisfying $a_{1}\leq \cdots\leq a_{k}$ and
$b_{1}\leq \cdots\leq b_{k}$. Assume that $a_{i}\neq a\pr_{j}$ and
$b_{i}\neq b\pr_{j}$ for all $i,j$. Introduce the corresponding $type~ A~ quantum~ minors$ by the following identities:
\begin{align*}
R_{k-1,k}(R_{k-2,k}R_{k-2,k-1})\cdots(R_{1,k}\cdots R_{1,2})T_{1}^{\pm}(u)\cdots T_{k}^{\pm}(u-k\hbar+\hbar)\\
=\sum_{a_{i},b_{i}}e_{a_{1},b_{1}}\otimes \cdots \otimes e_{a_{k},b_{k}}t^{\pm\ts a_{1}\cdots a_{k}}_{\ts b_{1}\cdots b_{k}}(u),
\end{align*}
where $R_{i,j}=R_{ij}((j-i)\hbar)$. Apply $R$-matrix calculations which are quite analogous to the Yangian case, we get the following relations:
\ben
\big[t^{\pm}_{a_{i}b_{j}}(u),t^{\pm\ts a_{1}\cdots a_{k}}_{\ts b_{1}\cdots b_{k}}(v)\big]=0,
\een
\ben
\frac{u_{+}-v_{-}-\hbar}{u_{+}-v_{-}}t^{+}_{a_{i}b_{j}}(u)t^{-\ts a_{1}\cdots a_{k}}_{\ts b_{1}\cdots b_{k}}(v)=
\frac{u_{-}-v_{+}-\hbar}{u_{-}-v_{+}}t^{-\ts a_{1}\cdots a_{k}}_{\ts b_{1}\cdots b_{k}}(v)t^{+}_{a_{i}b_{j}}(u),
\een
\begin{align*}
\prod_{i=1}^{k-1}\frac{(u_{-}-v_{+}+i\hbar)^{2}}{(u_{-}-v_{+}+i\hbar)^{2}-\hbar^{2}}\cdot\frac{u_{-}-v_{+}}{u_{-}-v_{+}+\hbar}t^{-}_{a_{i}b_{j}}(u)t^{+\ts a_{1}\cdots a_{k}}_{\ts b_{1}\cdots b_{k}}(v)\\
=\prod_{i=1}^{k-1}\frac{(u_{+}-v_{-}+i\hbar)^{2}}{(u_{+}-v_{-}+i\hbar)^{2}-\hbar^{2}}\cdot\frac{u_{+}-v_{-}}{u_{+}-v_{-}+\hbar}t^{+\ts a_{1}\cdots a_{k}}_{\ts b_{1}\cdots b_{k}}(v)t^{-}_{a_{i}b_{j}}(u),
\end{align*}
where $1\leq i,j\leq k$.

\bpr\label{cor:commu}
We have the relations
\ben
\big[t^{\pm}_{ab}(u),t^{[m]\pm}_{ij}(v)\big]=0,
\een
\ben
\big[t^{+}_{ab}(u),t^{[m]-}_{ij}(v)\big]=0,
\een
\ben
\frac{(u_{-}-v_{+})^{2}}{(u_{-}-v_{+})^{2}-\hbar^{2}}t^{-}_{ab}(u)t^{[m]+}_{ij}(v)=\frac{(u_{+}-v_{-})^{2}}{(u_{+}-v_{-})^{2}-\hbar^{2}}t^{[m]+}_{ij}(v)t^{-}_{ab}(u),
\een
for all $1\leqslant a,b\leqslant m$ and $m+1\leqslant i,j\leqslant (m+1)\pr$.
\epr

\bpf
All relations are verified with the use of the relations between the quasideterminants
and quantum minors:
\beq
\left|\begin{matrix}
t^{\pm}_{11}(u)&\dots&t^{\pm}_{1m}(u)&t^{\pm}_{1j}(u)\\
\dots&\dots&\dots&\dots\\
t^{\pm}_{m1}(u)&\dots&t^{\pm}_{mm}(u)&t^{\pm}_{mj}(u)\\
t^{\pm}_{i1}(u)&\dots&t^{\pm}_{im}(u)&\boxed{t^{\pm}_{ij}(u)}
\end{matrix}\right|=
t^{\pm\ts 1\cdots m}_{\ts 1\cdots m}(u+m\hbar)^{-1}t^{\pm\ts 1\cdots m\ts i}_{\ts 1\cdots m\ts j}(u+m\hbar).
\eeq
\epf

\section{Drinfeld realization of extended Yangian doubles}
\label{sec:dpey}

Drinfeld realization provides a new set of generators for the extended
Yangian double $\DX_\hbar(\g_N)$ that facilitates the study of finite dimensional representations \cite{io:br}.
%analogous to that of the Yangian double $\DY_\hbar(\gl_N)$
 We first recall the homomorphism between Yangian double in low ranks \cite{amr:otr}.
%for between classical Lie algebras
%in low ranks lead to corresponding  isomorphisms; see
%and \cite{jl:ib}. We begin by reviewing them in the context of Drinfeld presentations.

\subsection{Low rank isomorphisms}
\label{subsec:lri}

Introduce the normalized R-matrix
\begin{align*}
\bar{R}(u)=f(u)(I-\frac{\hbar}{u}P).
\end{align*}
Recall that the Yangian double $\DY_\hbar^{\bar{R}}(\gl_N)$ is defined as an associative algebra over $\Ac$ generated by $T_{ij}^{(r)}$ and $C$ where $1\leq i,j\leq N$, $r\in \mathbb{Z}^{\times}$
subject to the relations written in terms of
generating series:
\beq
[T^{\circ\pm}(u),C]=0,
\eeq
\beq
\bar{R}(u-v)T^{\circ\pm}_{1}(u)T^{\circ\pm}_2(v)=T^{\circ\pm}_2(v)T^{\circ\pm}_{1}(u)\bar{R}(u-v),
\eeq
\beq
\bar{R}(u_{+}-v_{-})T_{1}^{\circ +}(u)T_{2}^{\circ -}(v)=T_{2}^{\circ -}(v)T_{1}^{\circ +}(u)\bar{R}(u_{-}-v_{+}),
\eeq
where $u_{\pm}=u\pm \frac{1}{4}\hbar C$, the generating series $T^{\circ\pm}(u)=\sum_{i,j=1}^{N}T^{\pm}_{ij}(u)\otimes E_{ij}$ are defined by
\beq
T^{+}_{ij}(u)=\delta_{ij}+\hbar\sum_{r=1}^{\infty}T^{(r)}_{ij}u^{-r},\quad T^{-}_{ij}(u)=\delta_{ij}-\hbar\sum_{r=1}^{\infty}T^{(-r)}_{ij}u^{r-1}.
\eeq
We use $T^{\circ}$ here to distinguish the objects related to $\DY_\hbar^{\bar{R}}(\gl_N)$ from those related to the algebra $\DX_\hbar(\g_N)$. When $f(u)=1$, we denote $\DY_\hbar^{\bar{R}}(\gl_N)$ by $\DY_\hbar(\gl_N)$ for simplicity.
To distinguish two sets of Gaussian generators, we use capital letters to
denote the generating series %for follow the notation of \cite[Sec.~3]{io:br}
in the Gauss decomposition
of the generator matrix of the Yangian double $\DY_\hbar^{\bar{R}}(\gl_N)$, i.e.
 $K^{\pm}_i(u)$, $E^{\pm}_{ij}(u)$ and $F^{\pm}_{ji}(u)$
for the entries occurring in the $\DY_\hbar^{\bar{R}}(\gl_N)$ counterpart
of \eqref{gd}, while keeping the lower case letters for the
generating series for the corresponding factors in the Gauss decomposition
for the extended Yangian double $\DX_\hbar(\g_N)$.

\ble\label{lem:lowrank C}
The mapping
\begin{align*}
t^{\pm}_{ij}(u)&\mapsto T^{\pm}_{ij}(\frac{u}{2}),\quad i,j\in\{1,2\}\\
c&\mapsto 2C
\end{align*}
defines an isomorphism $\phi:\DX_\hbar(\spa_2)\to \DY_\hbar^{\bar{R}}(\gl_2)$, where $f(u)=\frac{2u-\hbar}{2u-2\hbar}$.
Moreover, the images of the Gaussian generators is given by
\begin{alignat*}{2}
k^{\pm}_{1}(u)&\mapsto K^{\pm}_{1}(u/2),\qquad\quad  & e^{\pm}_{12}(u)&\mapsto E^{\pm}_{12}(u/2),
\\[0.2em]
k^{\pm}_{2}(u)&\mapsto K^{\pm}_{2}(u/2),\qquad\quad & f^{\pm}_{21}(u)&\mapsto F^{\pm}_{21}(u/2).
\end{alignat*}
\ele

\bpf
By the similar method as \cite[Sec.~4.1]{amr:otr}, we can prove that $\phi$ is an isomorphism. Note that $k^{\pm}_{1}(u)=t^{\pm}_{11}(u)$ and $e^{\pm}_{12}(u)=t^{\pm}_{11}(u)^{-1}t^{\pm}_{12}(u)$. Hence $\phi(k^{\pm}_{1}(u))=T^{\pm}_{11}(u/2)=K^{\pm}_{1}(u/2)$. And $\phi(e^{\pm}_{12}(u))=T^{\pm}_{11}(u/2)^{-1}T^{\pm}_{12}(u/2)=E^{\pm}_{12}(u/2)$.
\epf

\ble\label{lem:lowrank B}
The mapping
\begin{align*}
T^{\pm}_{ij}(u)&\mapsto \frac{1+P}{2}T^{\circ\pm}_{1}(2u)T^{\circ\pm}_{2}(2u+\hbar),\quad i,j\in\{1,2\}\\
c&\mapsto \frac{1}{2}C
\end{align*}
defines an injective homomorphism $\phi:\DX_\hbar(\oa_3)\to \DY_\hbar^{\bar{R}}(\gl_2)$, where $f(u)$ is a formal power series in $u^{-1}$ satisfying the following equation:
\beq
f(u-\hbar)f^{2}(u)f(u+\hbar)=\frac{u+\hbar}{u-\hbar}.
\eeq
Moreover, the images of the Gaussian generators is given by
\begin{alignat*}{3}
k^{\pm}_{1}(u)&\mapsto K^{\pm}_{1}(2u)K^{\pm}_{1}(2u+\hbar),\qquad
&e^{\pm}_{12}(u)&\mapsto \sqrt{2}\tss E^{\pm}_{12}(2u+\hbar),\qquad
&f^{\pm}_{21}(u)&\mapsto \sqrt{2}\tss F^{\pm}_{21}(2u+\hbar),\\[0.3em]
k^{\pm}_{2}(u)&\mapsto K^{\pm}_{1}(2u)K^{\pm}_{2}(2u+\hbar),\qquad
&e^{\pm}_{23}(u)&\mapsto -\sqrt{2}\tss E^{\pm}_{12}(2u),\qquad
&f^{\pm}_{32}(u)&\mapsto -\sqrt{2}\tss F^{\pm}_{21}(2u),\\[0.3em]
k^{\pm}_3(u)&\mapsto K^{\pm}_{2}(2u)K^{\pm}_{2}(2u+\hbar),\qquad
&e^{\pm}_{13}(u)&\mapsto -E^{\pm}_{12}(2u+\hbar)^2,\qquad
&f^{\pm}_{31}(u)&\mapsto -F^{\pm}_{21}(2u+\hbar)^2.
\end{alignat*}
\ele

\bpf
Define a filtration on $\DY_\hbar^{\bar{R}}(\gl_2)$ by setting $\deg T_{ij}^{(r)}=r-1,$ $\deg T_{ij}^{(-r)}=-r$ for all $r\geqslant 1$ and $\deg C=0$.
By the similar method as \cite[Sec.~4.2]{amr:otr}, we can prove that $\phi$ is a homomorphism and $\phi$ induces an injective homomorphism of the associated graded algebras $\gr\phi:\gr\DX_\hbar(\oa_3)\to \gr\DY_\hbar^{\bar{R}}(\gl_2)$. Since the filtration on $\DX_\hbar(\oa_3)$ is separated, $\phi$ is also injective. The images of the Gaussian generators is obtained by tedious calculation. For instance,
\begin{align*}
\phi(e^{\pm}_{12}(u))&=K^{\pm}_{1}(2u+\hbar)^{-1}K^{\pm}_{1}(2u)^{-1}(\frac{1}{\sqrt{2}}T^{\pm}_{11}(2u)T^{\pm}_{12}(2u+\hbar)+\frac{1}{\sqrt{2}}T^{\pm}_{12}(2u)T^{\pm}_{11}(2u+\hbar))\\
&=\frac{1}{\sqrt{2}}E^{\pm}_{12}(2u+\hbar)+\frac{1}{\sqrt{2}}K^{\pm}_{1}(2u+\hbar)^{-1}E^{\pm}_{12}(2u)K^{\pm}_{1}(2u+\hbar)\\
&=\frac{1}{\sqrt{2}}E^{\pm}_{12}(2u+\hbar)+\frac{1}{\sqrt{2}}E^{\pm}_{12}(2u+\hbar)=\sqrt{2}E^{\pm}_{12}(2u+\hbar)
\end{align*}
\epf

\ble\label{lem:lowrank D}
The mapping
\begin{align*}
T^{\pm}(u)&\mapsto T^{\circ\pm}_{1}(u)T^{\circ\ts'\pm}_{2}(u),\\
c&\mapsto C
\end{align*}
defines an injective homomorphism $\phi:\DX_\hbar(\oa_4)\to \DY_\hbar^{\bar{R}}(\gl_2)\otimes\DY_\hbar^{\bar{R}}(\gl_2)$, where $f(u)$ is a formal power series in $u^{-1}$ satisfying the following equation:
\beq
f^{2}(u)=\frac{u}{u-\hbar}.
\eeq
We use notations from \cite[Sec.~4.3]{amr:otr} here. Moreover, the images of the Gaussian generators is given by
\begin{alignat*}{2}
k^{\pm}_{1}(u)&\mapsto K^{\pm}_{1}(u)\otimes K^{\pm}_{1}(u),\qquad\quad k^{\pm}_{2'}(u)&&\mapsto K^{\pm}_{2}(u)\otimes K^{\pm}_{1}(u),
\\[0.2em]
k^{\pm}_{2}(u)&\mapsto K^{\pm}_{1}(u)\otimes K^{\pm}_{2}(u),\qquad\quad
k^{\pm}_{1'}(u)&&\mapsto K^{\pm}_{2}(u)\otimes K^{\pm}_{2}(u),
\end{alignat*}
together with
\begin{alignat*}{2}
e^{\pm}_{12}(u)&\mapsto 1\otimes E^{\pm}_{12}(u),\qquad\quad & e^{\pm}_{12'}(u)&\mapsto E^{\pm}_{12}(u)\otimes 1,
\\[0.2em]
e^{\pm}_{11'}(u)&\mapsto -E^{\pm}_{12}(u)\otimes E^{\pm}_{12}(u),\qquad\quad &e^{\pm}_{22'}(u)&\mapsto 0,\\[0.2em]
e^{\pm}_{21'}(u)&\mapsto -E^{\pm}_{12}(u)\otimes 1,\qquad\quad &e^{\pm}_{2'1'}(u)&\mapsto -1\otimes E^{\pm}_{12}(u),
\end{alignat*}
and
\begin{alignat*}{2}
f^{\pm}_{21}(u)&\mapsto 1\otimes F^{\pm}_{21}(u),\qquad\quad & f^{\pm}_{2'1}(u)&\mapsto F^{\pm}_{21}(u)\otimes 1,
\\[0.2em]
f^{\pm}_{1'1}(u)&\mapsto -F^{\pm}_{21}(u)\otimes F^{\pm}_{21}(u),\qquad\quad &f^{\pm}_{2'2}(u)&\mapsto 0,\\[0.2em]
f^{\pm}_{1'2}(u)&\mapsto -F^{\pm}_{21}(u)\otimes 1,\qquad\quad &f^{\pm}_{1'2'}(u)&\mapsto -1\otimes F^{\pm}_{21}(u).
\end{alignat*}
\ele

\bpf
The lemma is proved similarly as type B.
\epf

\subsection{Type {\sl A} relations}
\label{subsec:tar}

Since the extended Yangian $\DX_\hbar(\g_{N})$ contains a subalgebra isomorphic
to the Yangian $\DY_\hbar(\gl_n)$, some relations between the Gaussian generators of
$\DX_\hbar(\g_{N})$ can be obtained from those of $\DY_\hbar(\gl_n)$.
We record them in the next proposition.

\bpr\label{pro:Arelations}
The following relations hold in
$\DX_\hbar(\g_{N})$, with the conditions on the indices $1\leq i,j\leq n-1$
and $1\leq l,m \leq n$.

$$k_{l}^{+}(u)k_{m}^{-}(v)=k_{m}^{-}(v)k_{l}^{+}(u)\quad (l<m),$$
$$\frac{u_{+}-v_{-}-\hbar}{u_{+}-v_{-}}k_{l}^{+}(u)k_{l}^{-}(v)=\frac{u_{-}-v_{+}-\hbar}{u_{-}-v_{+}}k_{l}^{-}(v)k_{l}^{+}(u),$$
\begin{equation}
\non
\frac{(v_{+}-u_{-})^{2}}{(v_{+}-u_{-})^2-\hbar^{2}}k_{m}^{+}(u)k_{l}^{-}(v)
=\frac{(v_{-}-u_{+})^{2}}{(v_{-}-u_{+})^2-\hbar^{2}}k_{l}^{-}(v)k_{m}^{+}(u)\quad (l<m),
\end{equation}
$$\frac{u_{\mp}-v_{\pm}-\hbar}{u_{\mp}-v_{\pm}}k_{l}^{\pm}(u)k_{l+1}^{\mp}(v)k_{l}^{\mp}(v)^{-1}=\frac{u_{\pm}-v_{\mp}-\hbar}{u_{\pm}-v_{\mp}}k_{l+1}^{\mp}(v)k_{l}^{\mp}(v)^{-1}k_{l}^{\pm}(u)\quad (l\leq n-1),$$
$$\frac{u_{\mp}-v_{\pm}+\hbar}{u_{\mp}-v_{\pm}}k_{l+1}^{\pm}(u)k_{l+1}^{\mp}(v)k_{l}^{\mp}(v)^{-1}=\frac{u_{\pm}-v_{\mp}+\hbar}{u_{\pm}-v_{\mp}}k_{l+1}^{\mp}(v)k_{l}^{\mp}(v)^{-1}k_{l+1}^{\pm}(u)\quad (l\leq n-1),$$
$$k_{i}^{\pm}(u)^{-1}X_{i}^{+}(v)k_{i}^{\pm}(u)=\frac{u_{\mp}-v-\hbar}{u_{\mp}-v}X_{i}^{+}(v),k_{i}^{\pm}(u)X_{i}^{-}(v)k_{i}^{\pm}(u)^{-1}=\frac{u_{\pm}-v-\hbar}{u_{\pm}-v}X_{i}^{-}(v),$$
$$k_{i+1}^{\pm}(u)^{-1}X_{i}^{+}(v)k_{i+1}^{\pm}(u)=\frac{u_{\mp}-v+\hbar}{u_{\mp}-v}X_{i}^{+}(v),k_{i+1}^{\pm}(u)X_{i}^{-}(v)k_{i+1}^{\pm}(u)^{-1}=\frac{u_{\pm}-v-\hbar}{u_{\pm}-v}X_{i}^{-}(v),$$
$$k_{l}^{\pm}(u)^{-1}X_{i}^{+}(v)k_{l}^{\pm}(u)=X_{i}^{+}(v),k_{l}^{\pm}(u)X_{i}^{-}(v)k_{l}^{\pm}(u)^{-1}=X_{i}^{-}(v)
\quad (l\neq i,i+1),$$
$$(u-v\pm \hbar)X_{i}^{\pm}(u)X_{i}^{\pm}(v)=(u-v\mp
\hbar)X_{i}^{\pm}(v)X_{i}^{\pm}(u),$$
$$(u-v-\hbar)X_{i}^{+}(u)X_{i+1}^{+}(v)=(u-v)X_{i+1}^{+}(v)X_{i}^{+}(u)\quad (i\leq n-2),$$
$$(u-v)X_{i}^{-}(u)X_{i+1}^{-}(v)=(u-v-\hbar)X_{i+1}^{-}(v)X_{i}^{-}(u)\quad (i\leq n-2),$$
\begin{align*}
&X_{i}^{\pm}(u_{1})X_{i}^{\pm}(u_{2})X_{j}^{\pm}(v)-2X_{i}^{\pm}(u_{1})X_{j}^{\pm}(v)X_{i}^{\pm}(u_{2})+X_{j}^{\pm}(v)X_{i}^{\pm}(u_{1})X_{i}^{\pm}(u_{2})+\{u_{1}\leftrightarrow
u_{2}\}\\&=0\quad if~ |i-j|=1,
\end{align*}
$$X_{i}^{\pm}(u)X_{j}^{\pm}(v)=X_{j}^{\pm}(v)X_{i}^{\pm}(u)\quad if~
|i-j|>1,$$
\begin{align*}
&[X_{i}^{+}(u),X_{j}^{-}(v)]=-\delta_{ij}\hbar\{\delta(u_{+}-v_{-})k_{i+1}^{+}(u_{+})k_{i}^{+}(u_{+})^{-1}-\delta(u_{-}-v_{+})k_{i+1}^{-}(v_{+})k_{i}^{-}(v_{+})^{-1}\}\\
&where ~\delta(u-v)=\Sigma_{k\in\mathbb{Z}}u^{-k-1}v^{k}.\nonumber
\end{align*}
\epr

\bpf
The coefficients of the series $t^{\pm}_{ij}(u)$ with $i,j\in\{1,\dots,n\}$ and $c$
generate a subalgebra of $\DX_\hbar(\g_{N})$ isomorphic to the Yangian $\DY_\hbar(\gl_n)$.
Hence, the elements $k^{\pm}_{i}(u)$ with $1\leq i\leq n$ and $X_{j}^{\pm}(v)$ with $1\leq j\leq n-1$
satisfy the same relations as the corresponding
elements of $\DY_\hbar(\gl_n)$; see \cite[Section~3]{io:br}.
\epf

We also have the following proposition, which can be proved by the similar method as in %The proof is similar as described in
\cite[Section~5]{ji:iso}.

\bpr\label{prop:Areldual}
All relations given in Proposition~\ref{pro:Arelations}
still hold if we replace the indices of all series
by $i\mapsto (n-i+1)'$ for $1\leqslant i\leqslant n$.
\epr

\subsection{Relations for Gaussian generators}
\label{subsec:rgd}

\bpr\label{pro:rellowrankB}
We have the relations in $\DX_\hbar(\oa_{2n+1})$, for $\epsilon=\pm$.
\begin{align}
k_{n}^{\pm}(u)^{-\epsilon}X_{n}^{\epsilon}(v)k_{n}^{\pm}(u)^{\epsilon}&=\frac{u_{\mp\epsilon}-v-\hbar}{u_{\mp\epsilon}-v}X_{n}^{\epsilon}(v), \\[0.3em]
%k_{n}^{\pm}(u)^{-1}X_{n}^{+}(v)k_{n}^{\pm}(u)&=\frac{u_{\mp}-v-1}{u_{\mp}-v}X_{n}^{+}(v), \\[0.3em]
%k_{n}^{\pm}(u)X_{n}^{-}(v)k_{n}^{\pm}(u)^{-1}&=\frac{u_{\pm}-v-1}{u_{\pm}-v}X_{n}^{-}(v), \\[0.3em]
k_{n+1}^{\pm}(u)^{-\epsilon}X_{n}^{\epsilon}(v)k_{n+1}^{\pm}(u)^{\epsilon}&
=\frac{(u_{\mp\epsilon}-v-\hbar)(u_{\mp\epsilon}-v+\frac{1}{2}\hbar)}{(u_{\mp\epsilon}-v)(u_{\mp\epsilon}-v-\frac{1}{2}\hbar)}X_{n}^{\epsilon}(v), \\[0.3em]
%k_{n+1}^{\pm}(u)X_{n}^{-}(v)k_{n+1}^{\pm}(u)^{-1}&=\frac{(u_{\pm}-v-1)(u_{\pm}-v+\frac{1}{2})}{(u_{\pm}-v)(u_{\pm}-v-\frac{1}{2})}X_{n}^{-}(v), \\[0.3em]
(u-v\pm \frac{1}{2}\hbar)X_{n}^{\pm}(u)X_{n}^{\pm}(v)&=(u-v\mp\frac{1}{2}\hbar)X_{n}^{\pm}(v)X_{n}^{\pm}(u), \\[0.3em]
%(u-v-\frac{1}{2})X_{n}^{-}(u)X_{n}^{-}(v)&=(u-v+\frac{1}{2})X_{n}^{-}(v)X_{n}^{-}(u), \\[0.3em]
[X_{n}^{+}(u),X_{n}^{-}(v)]&=2\hbar\{\delta(u_{-}-v_{+})k_{n+1}^{-}(v_{+})k_{n}^{-}(v_{+})^{-1}\\[0.3em]
&-\delta(u_{+}-v_{-})k_{n+1}^{+}(u_{+})k_{n}^{+}(u_{+})^{-1}\}.\nonumber
\end{align}
\epr

\bpf
By Corollary~\ref{cor:guass-embed}, the subalgebra $\DX_\hbar^{[n-1]}(\oa_N)$ of
$\DX_\hbar(\oa_N)$ is isomorphic to $\DX_\hbar(\oa_3)$.  Hence the relations are implied by
Lemma~\ref{lem:lowrank B}, and the Drinfeld presentation of the Yangian $\DY_\hbar(\gl_2)$ (see \cite[Sec.~3.1]{m:yc}), note that the $R$-matrix of $\DY_\hbar(\gl_2)$ should be modified and the relations are calculated accordingly. For instance, the first relation in the proposition can be verified by Lemma~\ref{lem:lowrank B}.
Since
\begin{align*}
K_{1}^{\pm}(2u+\hbar)^{-1}K_{1}^{\pm}(2u)^{-1}(\sqrt{2}E_{1}^{\mp}(2v+\hbar))K_{1}^{\pm}(2u)K_{1}^{\pm}(2u+\hbar)\\
=\frac{u-v-\hbar\mp \frac{\hbar C}{4}}{u-v\mp \frac{\hbar C}{4}}\sqrt{2}E_{1}^{\mp}(2v+\hbar)+\frac{\hbar}{u-v\mp \frac{\hbar C}{4}}\sqrt{2}E_{1}^{\pm}(2u+\hbar),
\end{align*}
then we have
\begin{align*}
k_{n}^{\pm}(u)^{-1}e_{n}^{\mp}(v)k_{n}^{\pm}(u)=\frac{u_{\mp}-v_{\pm}-\hbar}{u_{\mp}-v_{\pm}}e_{n}^{\mp}(v)+\frac{\hbar}{u_{\mp}-v_{\pm}}e_{n}^{\pm}(u).
\end{align*}
Similarly, we have
\ben
k_{n}^{\pm}(u)^{-1}e_{n}^{\pm}(v)k_{n}^{\pm}(u)=\frac{u-v-\hbar}{u-v}
e_{n}^{\pm}(v)+\frac{\hbar}{u-v}e_{n}^{\pm}(u).
\een
Since $X_{n}^{+}(v)=e_{n}^{+}(v_{+})-e_{n}^{-}(v_{-})$, then we have
\ben
k_{n}^{\pm}(u)^{-1}X_{n}^{+}(v)k_{n}^{\pm}(u)=\frac{u_{\mp}-v-\hbar}{u_{\mp}-v}X_{n}^{+}(v).
\een
\epf
The following two propositions are proved as in type B.
\bpr\label{pro:rellowrankC} For $\epsilon=\pm$ or $\pm 1$, the relations in $\DX_\hbar(\spa_{2n})$ are
\begin{align}
%k_{n}^{\pm}(u)^{-1}X_{n}^{+}(v)k_{n}^{\pm}(u)&=\frac{u_{\mp}-v-2}{u_{\mp}-v}X_{n}^{+}(v), \\[0.3em]
%k_{n}^{\pm}(u)X_{n}^{-}(v)k_{n}^{\pm}(u)^{-1}&=\frac{u_{\pm}-v-2}{u_{\pm}-v}X_{n}^{-}(v), \\[0.3em]
k_{n}^{\pm}(u)^{-\epsilon}X_{n}^{\epsilon}(v)k_{n}^{\pm}(u)^{\epsilon}&=\frac{u_{\mp\epsilon}-v-2\hbar}{u_{\mp\epsilon}-v}X_{n}^{\epsilon}(v), \\[0.3em]
k_{n+1}^{\pm}(u)^{-\epsilon}X_{n}^{\epsilon}(v)k_{n+1}^{\pm}(u)^{\epsilon}&=
\frac{u_{\mp\epsilon}-v+2\hbar}{u_{\mp\epsilon}-v}X_{n}^{\epsilon}(v), \\[0.3em]
%k_{n+1}^{\pm}(u)^{-1}X_{n}^{+}(v)k_{n+1}^{\pm}(u)&=\frac{u_{\mp}-v+2}{u_{\mp}-v}X_{n}^{+}(v), \\[0.3em]
%k_{n+1}^{\pm}(u)X_{n}^{-}(v)k_{n+1}^{\pm}(u)^{-1}&=\frac{u_{\pm}-v+2}{u_{\pm}-v}X_{n}^{-}(v), \\[0.3em]
(u-v\pm 2\hbar)X_{n}^{\pm}(u)X_{n}^{\pm}(v)&=(u-v\mp 2\hbar)X_{n}^{\pm}(v)X_{n}^{\pm}(u), \\[0.3em]
%(u-v-2)X_{n}^{-}(u)X_{n}^{-}(v)&=(u-v+2)X_{n}^{-}(v)X_{n}^{-}(u), \\[0.3em]
[X_{n}^{+}(u),X_{n}^{-}(v)]&=2\hbar\{\delta(u_{-}-v_{+})k_{n+1}^{-}(v_{+})k_{n}^{-}(v_{+})^{-1}\\[0.3em]
&-\delta(u_{+}-v_{-})k_{n+1}^{+}(u_{+})k_{n}^{+}(u_{+})^{-1}\}.\nonumber
\end{align}
\epr

%\bpf
%This proposition can be proved in a similar way as type B.
%\epf

\bpr\label{pro:rellowrankD}
We have the relations in $\DX_\hbar(\oa_{2n})$,
\begin{align}
%k_{n-1}^{\pm}(u)^{-1}X_{n}^{+}(v)k_{n-1}^{\pm}(u)&=\frac{u_{\mp}-v-1}{u_{\mp}-v}X_{n}^{+}(v), \\[0.3em]
%k_{n-1}^{\pm}(u)X_{n}^{-}(v)k_{n-1}^{\pm}(u)^{-1}&=\frac{u_{\pm}-v-1}{u_{\pm}-v}X_{n}^{-}(v), \\[0.3em]
k_{n-1}^{\pm}(u)^{-\epsilon}X_{n}^{\epsilon}(v)k_{n-1}^{\pm}(u)^{\epsilon}
&=\frac{u_{\mp\epsilon}-v-\hbar}{u_{\mp\epsilon}-v}X_{n}^{\epsilon}(v), \\[0.3em]
k_{n}^{\pm}(u)^{-\epsilon}X_{n}^{\epsilon}(v)k_{n}^{\pm}(u)^{\epsilon}
&=\frac{u_{\mp\epsilon}-v-\hbar}{u_{\mp\epsilon}-v}X_{n}^{\epsilon}(v), \\[0.3em]
%k_{n}^{\pm}(u)^{-1}X_{n}^{+}(v)k_{n}^{\pm}(u)&=\frac{u_{\mp}-v-1}{u_{\mp}-v}X_{n}^{+}(v), \\[0.3em]
%k_{n}^{\pm}(u)X_{n}^{-}(v)k_{n}^{\pm}(u)^{-1}&=\frac{u_{\pm}-v-1}{u_{\pm}-v}X_{n}^{-}(v),  \\[0.3em]
k_{n+1}^{\pm}(u)^{-\epsilon}X_{n}^{\epsilon}(v)k_{n+1}^{\pm}(u)^{\epsilon}
&=\frac{u_{\mp\epsilon}-v+\hbar}{u_{\mp\epsilon}-v}X_{n}^{\epsilon}(v),
%k_{n+1}^{\pm}(u)^{-1}X_{n}^{+}(v)k_{n+1}^{\pm}(u)&=\frac{u_{\mp}-v+1}{u_{\mp}-v}X_{n}^{+}(v), \\[0.3em]
%k_{n+1}^{\pm}(u)X_{n}^{-}(v)k_{n+1}^{\pm}(u)^{-1}&=\frac{u_{\pm}-v+1}{u_{\pm}-v}X_{n}^{-}(v),
\end{align}
\begin{align}
(u-v\pm \hbar)X_{n}^{\pm}(u)X_{n}^{\pm}(v)&=(u-v\mp \hbar)X_{n}^{\pm}(v)X_{n}^{\pm}(u), \\[0.3em]
%(u-v-1)X_{n}^{-}(u)X_{n}^{-}(v)&=(u-v+1)X_{n}^{-}(v)X_{n}^{-}(u), \\[0.3em]
[X_{n-1}^{\pm}(u),X_{n}^{\mp}(v)]&=[X_{n}^{\pm}(u),X_{n-1}^{\pm}(v)]=0,
\end{align}
\begin{align}
%X_{n-1}^{\pm}(u)X_{n}^{\pm}(v)&=X_{n}^{\pm}(v)X_{n-1}^{\pm}(u), \\[0.3em]
%X_{n-1}^{-}(u)X_{n}^{-}(v)&=X_{n}^{-}(v)X_{n-1}^{-}(u), \\[0.3em]
[X_{n}^{+}(u),X_{n}^{-}(v)]&=\hbar\{\delta(u_{-}-v_{+})k_{n+1}^{-}(v_{+})k_{n-1}^{-}(v_{+})^{-1}\\[0.3em]
&-\delta(u_{+}-v_{-})k_{n+1}^{+}(u_{+})k_{n-1}^{+}(u_{+})^{-1}\}.\nonumber
\end{align}
\epr

%\bpf
%This proposition can be proved in a similar way as type B.
%\epf

\ble\label{lem:emjmtkl} For nonnegative integer $m<[(N+1)/2]$ and
%For the parameter $m$ chosen in Section~\ref{sec:gd},
indices
$j,k,l$ such that $m+1\leqslant j,k,l\leqslant (m+1)'$ and $j\neq l'$,
the following relations hold in the extended Yangian double
$\DX_\hbar(\g_N)$,
\beql{emjmtkl}
\big[e^{\pm}_{m\tss j}(u), {t^{\ts[m]\pm}_{k\tss l}(v)}\big]
=\frac{\hbar}{u-v}{t^{\ts[m]\pm}_{k\tss j}(v)}\big(e^{\pm}_{m\tss l}(v)-e^{\pm}_{m\tss l}(u)\big),
\eeq
\beql{fjmmtkl}
\big[f^{\pm}_{j\tss m}(u), {t^{\ts[m]\pm}_{k\tss l}(v)}\big]
=\frac{\hbar}{u-v}\big(f^{\pm}_{k\tss m}(u)-f^{\pm}_{k\tss m}(v)\big)\ts{t^{\ts[m]\pm}_{j\tss l}(v)},
\eeq
\beql{emjmtklpm}
\big[e^{\pm}_{m\tss j}(u), {t^{\ts[m]\mp}_{k\tss l}(v)}\big]
=\frac{\hbar}{u_{\mp}-v_{\pm}}t_{kj}^{[m]\mp}(v)(e^{\mp}_{ml}(v)-e^{\pm}_{ml}(u)),
\eeq
\beql{fjmmtklpm}
\big[f^{\pm}_{j\tss m}(u), {t^{\ts[m]\mp}_{k\tss l}(v)}\big]
=\frac{\hbar}{u_{\pm}-v_{\mp}}(f^{\pm}_{km}(u)-f^{\mp}_{km}(v))t_{jl}^{[m]\mp}(v).
\eeq
\ele

\bpf
We only need to verify the relations for $m=1$; the general case
follows by applying the homomorphism $\psi_m$, see
Proposition~\ref{prop:gauss-consist}.
Let us verify \eqref{fjmmtklpm} as an example, the other relations can be checked similarly.
Since
\ben
k_1^{-}(v)=t_{11}^{-}(v),\qquad f_{k\tss 1}^{-}(v)=t_{k1}^{-}(v)t_{11}^{-}(v)^{-1},\qquad
e_{1\tss l}^{-}(v)=t_{11}^{-}(v)^{-1}t_{1\tss l}^{-}(v),
\een
we get
\ben
t^{\ts[1]-}_{kl}(v)=t_{kl}^{-}(v)-t_{k1}^{-}(v)t_{11}^{-}(v)^{-1}t_{1\tss l}^{-}(v)=t_{kl}^{-}(v)-f_{k1}^{-}(v)k_{1}^{-}(v)e_{1l}^{-}(v).
\een
By defining relations
\eqref{defrel+-}, we have
\ben
[t_{j\tss 1}^{+}(u), t_{k\tss l}^{-}(v)]=\frac{\hbar}{u_{+}-v_{-}}t_{k1}^{+}(u)t_{jl}^{-}(v)-\frac{\hbar}{u_{-}-v_{+}}t_{k1}^{-}(v)t_{jl}^{+}(u).
\een
Therefore,
\begin{multline}
\non
\big[t_{j\tss 1}^{+}(u), t^{\ts[1]-}_{k\tss l}(v)\big]
+\big[t_{j\tss 1}^{+}(u), f_{k\tss 1}^{-}(v)k_{1}^{-}(v)e_{1\tss l}^{-}(v)\big]\\
=\frac{\hbar}{u_{+}-v_{-}}(t^{+}_{k\tss 1}(u)t^{\ts[1]-}_{j\tss l}(v)+t^{+}_{k\tss 1}(u)f^{-}_{j\tss 1}(u)k^{-}_{1}(v)e^{-}_{1\tss l}(v))-
\frac{\hbar}{u_{-}-v_{+}}t^{-}_{k\tss 1}(v)t^{+}_{j\tss l}(u).
\end{multline}
Next we transform the second commutator
on the left hand side as
\begin{align*}
[t_{j\tss 1}^{+}(u), f_{k\tss 1}^{-}(v)k_{1}^{-}(v)e_{1\tss l}^{-}(v)]
&=[t_{j\tss 1}^{+}(u),f_{k\tss 1}^{-}(v)]t_{1\tss l}^{-}(v)+f_{k\tss 1}^{-}(v)[t_{j\tss 1}^{+}(u),t_{1\tss l}^{-}(v)]\\
&=[t_{j\tss 1}^{+}(u),t_{k\tss 1}^{-}(v)t_{1\tss 1}^{-}(v)^{-1}]t_{1\tss l}^{-}(v)+f_{k\tss 1}^{-}(v)[t_{j\tss 1}^{+}(u),t_{1\tss l}^{-}(v)]
\end{align*}
which coincides with
\begin{align*}
&[t_{j\tss 1}^{+}(u),f_{k\tss 1}^{-}(v)]e_{1\tss l}^{-}(v)
-f^{-}_{k\tss 1}(v)[t^{+}_{j\tss 1}(u),t^{-}_{11}(v)]e^{-}_{1\tss l}(v)
+f^{-}_{k\tss 1}(v)[t_{j\tss 1}^{+}(u),t^{-}_{1\tss l}(v)]\\
&=\frac{\hbar}{u_{+}-v_{-}}f^{+}_{k\tss 1}(u)k^{+}_{1}(u)t^{-}_{j\tss 1}(v)e^{-}_{1\tss l}(v)-\frac{\hbar}{u_{-}-v_{+}}f^{-}_{k\tss 1}(v)k^{-}_{1}(v)t^{+}_{j\tss 1}(u)e^{-}_{1\tss l}(v)\\
&-\frac{\hbar}{u_{+}-v_{-}}f^{-}_{k\tss 1}(v)k^{+}_{1}(u)t^{-}_{j\tss 1}(v)e^{-}_{1\tss l}(v)+\frac{\hbar}{u_{-}-v_{+}}f^{-}_{k\tss 1}(v)k^{-}_{1}(v)t^{+}_{j\tss 1}(u)e^{-}_{1\tss l}(v)\\
&+\frac{\hbar}{u_{+}-v_{-}}f^{-}_{k\tss 1}(v)k^{+}_{1}(u)t^{-}_{j\tss l}(v)-\frac{\hbar}{u_{-}-v_{+}}f^{-}_{k\tss 1}(v)k^{-}_{1}(v)t^{+}_{j\tss l}(u).
\end{align*}
Cancelling the common terms, we get that
\ben
\big[t^{+}_{j\tss 1}(u),t^{\ts[1]-}_{k\tss l}(v)\big]=\frac{\hbar}{u_{+}-v_{-}}\big(f^{+}_{k\tss 1}(u)-f^{-}_{k\tss 1}(v)\big)k^{+}_{1}(u)t^{\ts[1]-}_{j\tss l}(v).
\een
This yields
\ben
\big[f^{+}_{j\tss 1}(u),t^{\ts[1]-}_{k\tss l}(v)\big]=\frac{\hbar}{u_{+}-v_{-}}\big(f^{+}_{k\tss 1}(u)-f^{-}_{k\tss 1}(v)\big)t^{\ts[1]-}_{j\tss l}(v).
\een
Similarly, we can prove
$$\big[f^{-}_{j\tss 1}(u),t^{\ts[1]+}_{k\tss l}(v)\big]=\frac{\hbar}{u_{-}-v_{+}}\big(f^{-}_{k\tss 1}(u)-f^{+}_{k\tss 1}(v)\big)t^{\ts[1]+}_{j\tss l}(v).$$
So \eqref{fjmmtklpm} with $m=1$ holds.
\epf

\bpr\label{Ben-1en}
For $\DX_\hbar(\oa_{2n+1})$, we have
\begin{align}
[e^{\pm}_{n-1}(u),e^{\pm}_{n}(v)]=\frac{\hbar}{u-v}(e^{\pm}_{n-1,n+1}(v)-e^{\pm}_{n-1,n+1}(u)-e^{\pm}_{n-1}(v)e^{\pm}_{n}(v)+e^{\pm}_{n-1}(u)e^{\pm}_{n}(v)),\\
\label{Be+n-1e-n}
[e^{\pm}_{n-1}(u),e^{\mp}_{n}(v)]=\frac{\hbar}{u_{\mp}-v_{\pm}}(e^{\mp}_{n-1,n+1}(v)-e^{\pm}_{n-1,n+1}(u)-e^{\mp}_{n-1}(v)e^{\mp}_{n}(v)+e^{\pm}_{n-1}(u)e^{\mp}_{n}(v)).
\end{align}
\epr

\bpf
By Lemma~\ref{lem:emjmtkl}, we get
\begin{align*}
[e^{\pm}_{n-1}(u),t^{\ts[n-1]\pm}_{n\tss n+1}(v)]&=\frac{\hbar}{u-v}t^{\ts[n-1]\pm}_{n\tss n}(v)(e^{\pm}_{n-1,n+1}(v)-e^{\pm}_{n-1,n+1}(u))\\
&=[e^{\pm}_{n-1}(u),k^{\pm}_{n}(v)e^{\pm}_{n,n+1}(v)]\\
&=[e^{\pm}_{n-1}(u),k^{\pm}_{n}(v)]e^{\pm}_{n,n+1}(v)+k^{\pm}_{n}(v)[e^{\pm}_{n-1}(u),e^{\pm}_{n,n+1}(v)]\\
&=\frac{\hbar}{u-v}k^{\pm}_{n}(v)(e^{\pm}_{n-1}(v)-e^{\pm}_{n-1}(u))e^{\pm}_{n}(v)+k^{\pm}_{n}(v)[e^{\pm}_{n-1}(u),e^{\pm}_{n}(v)]\\
&=\frac{\hbar}{u-v}k^{\pm}_{n}(v)(e^{\pm}_{n-1,n+1}(v)-e^{\pm}_{n-1,n+1}(u)).
\end{align*}
Hence, $[e^{\pm}_{n-1}(u),e^{\pm}_{n}(v)]=\frac{\hbar}{u-v}(e^{\pm}_{n-1,n+1}(v)-e^{\pm}_{n-1,n+1}(u)-e^{\pm}_{n-1}(v)e^{\pm}_{n}(v)+e^{\pm}_{n-1}(u)e^{\pm}_{n}(v))$.
We can prove \eqref{Be+n-1e-n} similarly.
\epf

\bco\label{cor:Bxn-1xn}
For $\DX_\hbar(\oa_{2n+1})$, we have $(u-v-\hbar)X^{+}_{n-1}(u)X^{+}_{n}(v)=(u-v)X^{+}_{n}(v)X^{+}_{n-1}(u)$.
\eco

\bpf
One can prove the above corollary by Proposition~\ref{Ben-1en}.
\epf

\bpr\label{prop:C4ee}
For $\DX_\hbar(\spa_4)$, we have
\begin{align}\label{C4e1+e2+}
[e^{\pm}_{12}(u),e^{\pm}_{23}(v)]&=\frac{2\hbar}{u-v}(e^{\pm}_{13}(v)-e^{\pm}_{13}(u)-e^{\pm}_{12}(v)e^{\pm}_{23}(v)+e^{\pm}_{12}(u)e^{\pm}_{23}(v)),\\
\label{C4e1+e2-}
[e^{\pm}_{12}(u),e^{\mp}_{23}(v)]&=\frac{2\hbar}{u_{\mp}-v_{\pm}}(e^{\mp}_{13}(v)-e^{\pm}_{13}(u)-e^{\mp}_{12}(v)e^{\mp}_{23}(v)+e^{\pm}_{12}(u)e^{\mp}_{23}(v)).
\end{align}
\epr

\bpf
The proof of \eqref{C4e1+e2+} can be found in \cite{ji:iso}. Next we will prove \eqref{C4e1+e2-}.
From the defining relations, we have
\begin{align}\label{Ct12t23}
[t^{+}_{12}(u),t^{-}_{23}(v)]=&\frac{\hbar}{u_{+}-v_{-}}t^{+}_{22}(u)t^{-}_{13}(v)-\frac{\hbar}{u_{-}-v_{+}}t^{-}_{22}(v)t^{+}_{13}(u)\\
&+\frac{\hbar}{u_{-}-v_{+}-3\hbar}(t^{-}_{24}(v)t^{+}_{11}(u)+t^{-}_{23}(v)t^{+}_{12}(u)-t^{-}_{22}(v)t^{+}_{13}(u)-t^{-}_{21}(v)t^{+}_{14}(u)).\nonumber
\end{align}

The left hand side of \eqref{Ct12t23} can be written as
\begin{align*}
&k^{+}_{1}(u)[e^{+}_{12}(u),k^{-}_{2}(v)e^{-}_{23}(v)]+[k^{+}_{1}(u),k^{-}_{2}(v)e^{-}_{23}(v)]e^{+}_{12}(u)\\
&+[k^{+}_{1}(u)e^{+}_{12}(u),f^{-}_{21}(v)k^{-}_{1}(v)]e^{-}_{13}(v)+f^{-}_{21}(v)k^{-}_{1}(v)[k^{+}_{1}(u)e^{+}_{12}(u),e^{-}_{13}(v)].
\end{align*}
It follows from the defining relations that
\begin{align*}
&[k^{+}_{1}(u)e^{+}_{12}(u),f^{-}_{21}(v)k^{-}_{1}(v)]=[t^{+}_{12}(u),t^{-}_{21}(v)]=\frac{\hbar}{u_{+}-v_{-}}t^{+}_{22}(u)t^{-}_{11}(v)-\frac{\hbar}{u_{-}-v_{+}}t^{-}_{22}(v)t^{+}_{11}(u)\\
&=\frac{\hbar}{u_{+}-v_{-}}(k^{+}_{2}(u)+f^{+}_{21}(u)k^{+}_{1}(u)e^{+}_{12}(u))k^{-}_{1}(v)-\frac{\hbar}{u_{-}-v_{+}}(k^{-}_{2}(v)+f^{-}_{21}(v)k^{-}_{1}(v)e^{-}_{12}(v))k^{+}_{1}(u).
\end{align*}
To compute the Lie bracket $[k^{+}_{1}(u)e^{+}_{12}(u),e^{-}_{13}(v)]=[t^{+}_{12}(u),t^{-}_{11}(v)^{-1}t^{-}_{13}(v)]$,
%we should calculate the Lie bracket $[t^{+}_{12}(u),t^{-}_{11}(v)^{-1}]$ and $[t^{+}_{12}(u),t^{-}_{13}(v)]$ firstly.
We have $[t^{+}_{12}(u),t^{-}_{11}(v)^{-1}]=-t^{-}_{11}(v)^{-1}[t^{+}_{12}(u),t^{-}_{11}(v)]t^{-}_{11}(v)^{-1}$, which can be calculated easily by the defining relations.
Also we have $[e^{+}_{12}(u),k^{-}_{2}(v)]=\frac{\hbar}{u_{-}-v_{+}}k^{-}_{2}(v)(e^{-}_{1}(v)-e^{+}_{1}(u))$ by type A relations.
Hence, the left hand side of \eqref{Ct12t23} takes the form
\begin{align*}
&k^{+}_{1}(u)k^{-}_{2}(v)[e^{+}_{12}(u),e^{-}_{23}(v)]+\frac{\hbar}{u_{-}-v_{+}}k^{+}_{1}(u)k^{-}_{2}(v)(e^{-}_{1}(v)e^{-}_{2}(v)-e^{+}_{1}(u)e^{-}_{2}(v))\\
&+\frac{\hbar}{u_{+}-v_{-}}(k^{+}_{2}(u)k^{-}_{1}(v)+f^{+}_{21}(u)k^{+}_{1}(u)e^{+}_{12}(u)k^{-}_{1}(v))e^{-}_{13}(v)-\frac{\hbar}{u_{-}-v_{+}}(k^{-}_{2}(v)k^{+}_{1}(u)e^{-}_{13}(v)\\
&+f^{-}_{21}(v)k^{-}_{1}(v)e^{-}_{12}(v)k^{+}_{1}(u)e^{-}_{13}(v))+\frac{\hbar}{u_{-}-v_{+}-3\hbar}f^{-}_{21}(v)k^{-}_{1}(v)(e^{-}_{14}(v)k^{+}_{1}(u)+e^{-}_{13}(v)k^{+}_{1}(u)e^{+}_{12}(u)\\
&-e^{-}_{12}(v)k^{+}_{1}(u)e^{+}_{13}(u)-k^{+}_{1}(u)e^{+}_{14}(u))+\frac{\hbar}{u_{-}-v_{+}}(f^{-}_{21}(v)k^{-}_{1}(v)e^{-}_{12}(v)k^{+}_{1}(u)e^{-}_{13}(v)\\
&-f^{-}_{21}(v)k^{-}_{1}(v)e^{-}_{12}(v)k^{+}_{1}(u)e^{+}_{13}(u)).
\end{align*}
The right hand side of \eqref{Ct12t23} takes the form
\begin{align*}
&\frac{\hbar}{u_{+}-v_{-}}(f^{+}_{21}(u)k^{+}_{1}(u)e^{+}_{12}(u)+k^{+}_{2}(u))k^{-}_{1}(v)e^{-}_{13}(v)-\frac{\hbar}{u_{-}-v_{+}}(f^{-}_{21}(v)k^{-}_{1}(v)e^{-}_{12}(v)+k^{-}_{2}(v))\\
&k^{+}_{1}(u)e^{+}_{13}(u)+\frac{\hbar}{u_{-}-v_{+}-3\hbar}[(f^{-}_{21}(v)k^{-}_{1}(v)e^{-}_{14}(v)+k^{-}_{2}(v)e^{-}_{24}(v))k^{+}_{1}(u)+(f^{-}_{21}(v)k^{-}_{1}(v)e^{-}_{13}(v)\\
&+k^{-}_{2}(v)e^{-}_{23}(v))k^{+}_{1}(u)e^{+}_{12}(u)-(f^{-}_{21}(v)k^{-}_{1}(v)e^{-}_{12}(v)+k^{-}_{2}(v))k^{+}_{1}(u)e^{+}_{13}(u)-f^{-}_{21}(v)k^{-}_{1}(v)k^{+}_{1}(u)e^{+}_{14}(u)].
\end{align*}

By cancelling the common terms, we have
\begin{align*}
&[e^{+}_{12}(u),e^{-}_{23}(v)]=\frac{\hbar}{u_{-}-v_{+}}(e^{-}_{13}(v)-e^{+}_{13}(u)+e^{+}_{12}(u)e^{-}_{23}(v)-e^{-}_{12}(v)e^{-}_{23}(v))+\frac{\hbar}{u_{-}-v_{+}-3\hbar}k^{+}_{1}(u)^{-1}\\
&(e^{-}_{24}(v)k^{+}_{1}(u)+e^{-}_{23}(v)k^{+}_{1}(u)e^{+}_{12}(u)-k^{+}_{1}(u)e^{+}_{13}(u))=\frac{\hbar}{u_{-}-v_{+}}(e^{-}_{13}(v)-e^{+}_{13}(u)+e^{+}_{12}(u)e^{-}_{23}(v)\\
&-e^{-}_{12}(v)e^{-}_{23}(v))+\frac{\hbar}{u_{-}-v_{+}-2\hbar}(e^{-}_{24}(v)+e^{-}_{2}(v)e^{+}_{1}(u)-e^{+}_{13}(u))=\frac{2\hbar}{u_{-}-v_{+}}(e^{-}_{13}(v)-e^{+}_{13}(u)\\
&-e^{-}_{1}(v)e^{-}_{2}(v)+e^{+}_{1}(u)e^{-}_{2}(v))+\frac{\hbar}{u_{-}-v_{+}-\hbar}(e^{-}_{24}(v)+e^{-}_{1}(v)e^{-}_{2}(v)-e^{-}_{13}(v)).
\end{align*}
Therefore, $e^{-}_{24}(v)+e^{-}_{1}(v)e^{-}_{2}(v)-e^{-}_{13}(v)=0$. Then $$[e^{+}_{12}(u),e^{-}_{23}(v)]=\frac{2\hbar}{u_{-}-v_{+}}(e^{-}_{13}(v)-e^{+}_{13}(u)-e^{-}_{12}(v)e^{-}_{23}(v)+e^{+}_{12}(u)e^{-}_{23}(v)).$$ Similarly, we can prove
$$[e^{-}_{12}(u),e^{+}_{23}(v)]=\frac{2\hbar}{u_{+}-v_{-}}(e^{+}_{13}(v)-e^{-}_{13}(u)-e^{+}_{12}(v)e^{+}_{23}(v)+e^{-}_{12}(u)e^{+}_{23}(v)).$$
This finishes the proof of \eqref{C4e1+e2-}.
\epf

The following result follows easily from Proposition~\ref{prop:C4ee}.
\bco\label{cor:Cxn-1xn}
For $\DX_\hbar(\spa_{2n})$, we have
$$[e^{\pm}_{n-1}(u),e^{\pm}_{n}(v)]=\frac{2\hbar}{u-v}(e^{\pm}_{n-1,n+1}(v)-e^{\pm}_{n-1,n+1}(u)-e^{\pm}_{n-1}(v)e^{\pm}_{n}(v)+e^{\pm}_{n-1}(u)e^{\pm}_{n}(v)),$$
$$[e^{\pm}_{n-1}(u),e^{\mp}_{n}(v)]=\frac{2\hbar}{u_{\mp}-v_{\pm}}(e^{\mp}_{n-1,n+1}(v)-e^{\pm}_{n-1,n+1}(u)-e^{\mp}_{n-1}(v)e^{\mp}_{n}(v)+e^{\pm}_{n-1}(u)e^{\mp}_{n}(v)),$$
$$(u-v-2\hbar)X^{+}_{n-1}(u)X^{+}_{n}(v)=(u-v)X^{+}_{n}(v)X^{+}_{n-1}(u).$$
\eco

%\bpf
%One can prove the above corollary by Proposition~\ref{prop:C4ee}.
%\epf

\bpr\label{prop:BCen-1fn}
For $\DX_\hbar(\oa_{2n+1}),\DX_\hbar(\spa_{2n})$, we have $[e^{\pm}_{n-1}(u),f^{\mp}_{n}(v)]=[e^{\pm}_{n-1}(u),f^{\pm}_{n}(v)]$\\
$=[e^{\pm}_{n}(u),f^{\mp}_{n-1}(v)]=[e^{\pm}_{n}(u),f^{\pm}_{n-1}(v)]=0.$
\epr

\bpf
As $t^{\ts[n-1]\mp}_{n+1\tss n}(v)=f^{\mp}_{n}(v)k^{\mp}_{n}(v)$, it follows from Lemma~\ref{lem:emjmtkl} that
\beq
[e^{\pm}_{n-1}(u),f^{\mp}_{n}(v)k^{\mp}_{n}(v)]=\frac{\hbar}{u_{\mp}-v_{\pm}}f^{\mp}_{n}(v)k^{\mp}_{n}(v)(e^{\mp}_{n-1}(v)-e^{\pm}_{n-1}(u)).
\eeq
On the other hand, by type $A$ relations we have
%\beq
%[e^{\pm}_{n-1}(u),f^{\mp}_{n}(v)k^{\mp}_{n}(v)]=[e^{\pm}_{n-1}(u),f^{\mp}_{n}(v)]k^{\mp}_{n}(v)+f^{\mp}_{n}(v)[e^{\pm}_{n-1}(u),k^{\mp}_{n}(v)],
%\eeq
%whereas
\beq
[e^{\pm}_{n-1}(u),k^{\mp}_{n}(v)]=\frac{\hbar}{u_{\mp}-v_{\pm}}k^{\mp}_{n}(v)(e^{\mp}_{n-1}(v)-e^{\pm}_{n-1}(u)),
\eeq
which forces that $[e^{\pm}_{n-1}(u),f^{\mp}_{n}(v)]=0$. The other relations are verified similarly.
\epf

\bco\label{cor:BCxnxn-1}
For $\DX_\hbar(\oa_{2n+1}),\DX_\hbar(\spa_{2n})$, we have $[X^{+}_{n}(u),X^{-}_{n-1}(v)]=[X^{+}_{n-1}(u),X^{-}_{n}(v)]$\\
$=0.$
\eco

\bpf As $X^+_i(u)=e^+_i(u_+)-e^-_i(u_-)$ and $X^-_i(u)=-f^-_i(u_+)+f^+_i(u_-)$, the relations follows immediately from Prop.~\ref{prop:BCen-1fn}.
%One can prove the above corollary by Proposition~\ref{prop:BCen-1fn}.
\epf

Now we recall a symmetry property \cite{ji:iso} for the matrix elements
in the Gauss decomposition, which will be useful
for later discussion.

\bpr\label{prop:eiei'} The following symmetry relations hold in $\DX_\hbar(\g_{N})$,
\beq
e^{\pm}_{(i+1)'\ts i'}(u)=-e^{\pm}_{i}(u+\ka\hbar-i\hbar)\Fand
f^{\pm}_{i'\ts (i+1)'}(u)=-f^{\pm}_{i}(u+\ka\hbar-i\hbar)
\eeq
for $i=1,\dots,n-1$.
In addition, for $\g_N=\spa_{2n}$ one has that
\ben
e^{\pm}_{n+1}(u)=-e^{\pm}_{n-1}\big(u+2\hbar\big)\Fand f^{\pm}_{n+1}(u)=-f^{\pm}_{n-1}\big(u+2\hbar\big),
\een
and for $\g_N=\oa_{2n}$ one has that
\ben
e^{\pm}_{n+1\ts n+2}(u)=-e^{\pm}_{n-1}(u)\Fand f^{\pm}_{n+2\ts n+1}(u)=-f^{\pm}_{n-1}(u).
\een
\epr

\bpr\label{prop:BCDkn+1xn-1} In $\DX_\hbar(\g_{N})$, for $\epsilon=\pm$ one has the following
\begin{align*}
[k^{\pm}_{n+1}(u),X^{\epsilon}_{n-1}(v)]&=0 &\text{for } \g_{N}=\oa_{2n+1},\\
k^{\pm}_{n+1}(u)^{-1}X^{\epsilon}_{n-1}(v)k^{\pm}_{n+1}(u)&=\frac{u_{\mp\epsilon}-v+\hbar}{u_{\mp\epsilon}-v+2\hbar}X^{\epsilon}_{n-1}(v)
&\text{for } \g_{N}=\spa_{2n},\\
%k^{\pm}_{n+1}(u)^{-1}X^{+}_{n-1}(v)k^{\pm}_{n+1}(u)=\frac{u_{\mp}-v+1}{u_{\mp}-v+2}X^{+}_{n-1}(v),
k^{\pm}_{n+1}(u)^{-1}X^{\epsilon}_{n-1}(v)k^{\pm}_{n+1}(u)&=\frac{u_{\mp\epsilon}-v-\hbar}{u_{\mp\epsilon}-v}X^{\epsilon}_{n-1}(v)
&\text{for } \g_{N}=\oa_{2n}.
\end{align*}
%\begin{align*}
%[k^{\pm}_{n+1}(u),X^{+}_{n-1}(v)]=0\quad and \quad [k^{\pm}_{n+1}(u),X^{-}_{n-1}(v)]=0
%\end{align*}
%\begin{align*}
%k^{\pm}_{n+1}(u)^{-1}X^{+}_{n-1}(v)k^{\pm}_{n+1}(u)=\frac{u_{\mp}-v+1}{u_{\mp}-v+2}X^{+}_{n-1}(v),\\ k^{\pm}_{n+1}(u)^{-1}X^{-}_{n-1}(v)k^{\pm}_{n+1}(u)=\frac{u_{\pm}-v+1}{u_{\pm}-v+2}X^{-}_{n-1}(v)
%\end{align*}
%for $\g_{N}=\spa_{2n}$,
%\begin{align*}
%k^{\pm}_{n+1}(u)^{-1}X^{+}_{n-1}(v)k^{\pm}_{n+1}(u)=\frac{u_{\mp}-v-1}{u_{\mp}-v}X^{+}_{n-1}(v),\\ k^{\pm}_{n+1}(u)^{-1}X^{-}_{n-1}(v)k^{\pm}_{n+1}(u)=\frac{u_{\pm}-v-1}{u_{\pm}-v}X^{-}_{n-1}(v)
%\end{align*}
%for $\g_{N}=\oa_{2n}$.
\epr

\bpf
First consider $\g_N=\oa_{2n+1}$.
Corollary~\ref{cor:commu} says that
$k^{\pm}_{n+1}(u)$ commutes with the elements of the subalgebra
generated by the $t^{\pm}_{ij}(u)$ with $1\leqslant i,j \leqslant n$, thus the relations follow. Now let $\g_N=\spa_{2n}$ or $\oa_{2n}$.
Corollary~\ref{cor:guass-embed} implies that the subalgebra
$\DX_\hbar^{[n-2]}(\g_N)$ of $\DX_\hbar(\g_N)$ is isomorphic to $\DX_\hbar(\g_4)$.
Now we apply Proposition~\ref{prop:Areldual} to the subalgebra $\DX_\hbar^{[n-2]}(\g_N)$, then we have
\ben
k^{\pm}_{n+1}(u)^{-1}e^{\pm}_{n+1,n+2}(v)k^{\pm}_{n+1}(u)=\frac{u-v-\hbar}{u-v}e^{\pm}_{n+1,n+2}(v)+\frac{\hbar}{u-v}e^{\pm}_{n+1,n+2}(u),
\een
\ben
k^{\pm}_{n+1}(u)^{-1}e^{\mp}_{n+1,n+2}(v)k^{\pm}_{n+1}(u)=\frac{u_{\mp}-v_{\pm}-\hbar}{u_{\mp}-v_{\pm}}e^{\mp}_{n+1,n+2}(v)+\frac{\hbar}{u_{\mp}-v_{\pm}}e^{\pm}_{n+1,n+2}(u),
\een
\ben
k^{\pm}_{n+1}(u)f^{\pm}_{n+2,n+1}(v)k^{\pm}_{n+1}(u)^{-1}=\frac{u-v-\hbar}{u-v}f^{\pm}_{n+2,n+1}(v)+\frac{\hbar}{u-v}f^{\pm}_{n+2,n+1}(u),
\een
\ben
k^{\pm}_{n+1}(u)f^{\mp}_{n+2,n+1}(v)k^{\pm}_{n+1}(u)^{-1}=\frac{u_{\pm}-v_{\mp}-\hbar}{u_{\pm}-v_{\mp}}f^{\mp}_{n+2,n+1}(v)+\frac{\hbar}{u_{\pm}-v_{\mp}}f^{\pm}_{n+2,n+1}(u).
\een
Finally the relations follow by applying
Proposition~\ref{prop:eiei'}.
\epf

\bpr\label{pro:relxnxi}
For $i=1,\dots,n$ and $\epsilon=\pm$ we have in the algebra $\DX_\hbar(\g_{N})$
\begin{align}\label{Bkix+n}
k^{\pm}_{i}(u)^{-1}X^{\epsilon}_{n}(v)k^{\pm}_{i}(u)&=(\delta_{in}\frac{u_{\mp\epsilon}-v-\hbar}{u_{\mp\epsilon}-v}+\delta_{i<n})X^{\epsilon}_{n}(v)
&\text{for }\g_{N}=\oa_{2n+1},\\ \label{Ckix-n}
%k^{\pm}_{i}(u)^{-1}X^{-}_{n}(v)k^{\pm}_{i}(u)=(\delta_{in}\frac{u_{\pm}-v-1}{u_{\pm}-v}+\delta_{i<n})X^{-}_{n}(v)
k^{\pm}_{i}(u)^{-1}X^{\epsilon}_{n}(v)k^{\pm}_{i}(u)&=(\delta_{in}\frac{u_{\mp\epsilon}-v-2\hbar}{u_{\mp\epsilon}-v}+\delta_{i<n})X^{\epsilon}_{n}(v)
&\text{for }\g_{N}=\spa_{2n},\\ \label{Dkix-n}
k^{\pm}_{i}(u)^{-1}X^{\epsilon}_{n}(v)k^{\pm}_{i}(u)&=(\delta_{i>n-2}\frac{u_{\mp\epsilon}-v-\hbar}{u_{\mp\epsilon}-v}+\delta_{i<n-1})X^{\epsilon}_{n}(v)
&\text{for }\g_{N}=\oa_{2n}.
\end{align}
%for $\g_{N}=\oa_{2n+1}$,
%\begin{align}
%k^{\pm}_{i}(u)^{-1}X^{+}_{n}(v)k^{\pm}_{i}(u)=(\delta_{in}\frac{u_{\mp}-v-2}{u_{\mp}-v}+\delta_{i<n})X^{+}_{n}(v)\\
%\label{Ckix-n}
%k^{\pm}_{i}(u)^{-1}X^{-}_{n}(v)k^{\pm}_{i}(u)=(\delta_{in}\frac{u_{\pm}-v-2}{u_{\pm}-v}+\delta_{i<n})X^{-}_{n}(v)
%\end{align}
%for $\g_{N}=\spa_{2n}$,
%\begin{align}
%k^{\pm}_{i}(u)^{-1}X^{+}_{n}(v)k^{\pm}_{i}(u)=(\delta_{i>n-2}\frac{u_{\mp}-v-1}{u_{\mp}-v}+\delta_{i<n-1})X^{+}_{n}(v)\\
%\label{Dkix-n}
%k^{\pm}_{i}(u)^{-1}X^{-}_{n}(v)k^{\pm}_{i}(u)=(\delta_{i>n-2}\frac{u_{\pm}-v-1}{u_{\pm}-v}+\delta_{i<n-1})X^{-}_{n}(v)
%\end{align}
%for $\g_{N}=\oa_{2n}$,
Moreover in all three cases, we have that
\begin{align}\label{BCDxixn}
[X^{+}_{i}(u),X^{-}_{n}(v)]&=[X^{+}_{n}(u),X^{-}_{i}(v)]=0
\quad (i<n), \\ \label{BCDkn+1xi}
[k^{\pm}_{n+1}(u),X^{+}_{i}(v)]&=[k^{\pm}_{n+1}(u),X^{-}_{i}(v)]=0 \quad (i<n-1).
\end{align}
%Moreover, if $i<n-1$ then
%\beq\label{BCDkn+1xi}
%[k^{\pm}_{n+1}(u),X^{+}_{i}(v)]=[k^{\pm}_{n+1}(u),X^{-}_{i}(v)]=0.
%\eeq
Here the notation $\delta_{i<n}$ means $\delta_{i<n}=1$ if $i<n$, otherwise $\delta_{i<n}=0$.
\epr

\bpf The relations \eqref{Bkix+n}-\eqref{Dkix-n} for $i=n$ were already given in Props.~\ref{pro:rellowrankB}-\ref{pro:rellowrankD}.
Now let $\g_N=\oa_{2n}$. By the defining relations
\eqref{defrel++} and \eqref{defrel+-}, the subalgebra generated by the coefficients of the series $t^{\pm}_{ij}(u)$
with $i,j\in J=\{1,\dots,n-1, n+1\}$ and $c$
is isomorphic to $\DY_\hbar(\gl_n)$. Furthermore, Lemma~\ref{lem:lowrank D} implies that
\ben
e^{\pm}_{n\ts n+1}(u)=f^{\pm}_{n+1\ts n}(u)=0.
\een
So we have the following Gauss decomposition
\beql{gdj}
T^{\pm}_J(u)=F^{\pm}_J(u)\ts H^{\pm}_J(u)\ts E^{\pm}_J(u),
\eeq
where the subscript $J$ indicates the submatrices in \eqref{gd} with rows and columns
labelled by integers of $J$.  This means that the entries of the matrices
$F^{\pm}_J(u)$, $H^{\pm}_J(u)$ and $E^{\pm}_J(u)$ satisfy the type $A$ relations as given in
 Prop.~\ref{pro:Arelations}, and the case of type $D$ is done.
Now let $\g_N=\oa_{2n+1}$ or $\spa_{2n}$. Almost all the relations
\eqref{Bkix+n}-\eqref{Ckix-n} with $i<n$, as well as \eqref{BCDxixn}
and \eqref{BCDkn+1xi} with $i<n-1$ follow from Cor.~\ref{cor:commu}. For instance, to check the commutation relation between $k_{1}^{-}(u)$ and $f_{n}^{+}(v)$, we get
$$\frac{(u_{-}-v_{+})^{2}}{(u_{-}-v_{+})^{2}-\hbar^{2}}t_{11}^{-}(u)t^{[n-1]+}_{n+1~n}(v)=\frac{(u_{+}-v_{-})^{2}}{(u_{+}-v_{-})^{2}-\hbar^{2}}t^{[n-1]+}_{n+1~n}(v)t_{11}^{-}(u)$$
from Cor.~\ref{cor:commu}. But $t_{11}^{-}(u)=k_{1}^{-}(u)$ and $t^{[n-1]+}_{n+1~n}(v)=f_{n}^{+}(v)k_{n}^{+}(v)$. We also have
$$\frac{(u_{-}-v_{+})^{2}}{(u_{-}-v_{+})^{2}-\hbar^{2}}k_{1}^{-}(u)k_{n}^{+}(v)=\frac{(u_{+}-v_{-})^{2}}{(u_{+}-v_{-})^{2}-\hbar^{2}}k_{n}^{+}(v)k_{1}^{-}(u).$$
It follows that $k_{1}^{-}(u)f_{n}^{+}(v)=f_{n}^{+}(v)k_{1}^{-}(u).$ So we have $k_{1}^{-}(u)^{-1}X_{n}^{-}(v)k_{1}^{-}(u)=X_{n}^{-}(v).$
As for the cases $i=n-1$ of \eqref{BCDxixn}, see Prop.~\ref{prop:BCen-1fn}.
\epf

\bpr\label{prop:Dxn-2xn}
We have the relations in $\DX_\hbar(\oa_{2n})$ for $\epsilon=\pm$ or $\pm \hbar$:
\begin{align}\label{Dx+n-2x+n}
&(u-v-(\hbar+\epsilon)/2)X^{\epsilon}_{n-2}(u)X^{\epsilon}_{n}(v)=(u-v-(\hbar-\epsilon)/2)X^{\epsilon}_{n}(v)X^{\epsilon}_{n-2}(u). %\\
%\label{Dx-n-2x-n}
%&(u-v)X^{-}_{n-2}(u)X^{-}_{n}(v)=(u-v-1)X^{-}_{n}(v)X^{-}_{n-2}(u).
\end{align}
\epr

\bpf The Gauss generators occurring as the entries of $F^{\pm}_J(u)$, $H^{\pm}_J(u)$ and $E^{\pm}_J(u)$ in \eqref{gdj}
satisfy the Yangian relations of type $A$ (see the proof of Proposition~\ref{pro:relxnxi}), so
\eqref{Dx+n-2x+n} follows from the Drinfeld realization of the Yangian double $\DY_\hbar(\gl_n)$.
%decomposition holds for
%As we pointed out in the proof of Proposition~\ref{pro:relxnxi},
%the Gaussian
%generators which occur as the entries of the matrices
%$F^{\pm}_J(u)$, $H^{\pm}_J(u)$ and $E^{\pm}_J(u)$ in \eqref{gdj}
%satisfy the type $A$ relations as described
%in Proposition~\ref{pro:Arelations}. Therefore,
%\eqref{Dx+n-2x+n} and \eqref{Dx-n-2x-n}
%follow from the Drinfeld presentation of the double Yangian $\DY_\hbar(\gl_n)$.
\epf

\bpr\label{prop:BCDxixn} For $1\leq i\leq n-1$ and $(\alpha_{i},\alpha_{n})=0$, we have
\begin{align}\label{xixn=0}
[X^{\pm}_{i}(u),X^{\pm}_{n}(v)]=0.
\end{align}
\epr

\bpf
If $i\leqslant n-2$ in types $B$ and $C$ or $i\leqslant n-3$ in type $D$, the relation
\eqref{xixn=0} follows from Corollary~\ref{cor:commu}. If $i=n-1$ in type $D$,
\eqref{xixn=0} is given by Lemma~\ref{lem:lowrank D} and
Corollary~\ref{cor:guass-embed}.
\epf

\bpr\label{prop:BCDknkn+1}
For $\g_N=\oa_{2n+1}$ we have
\begin{align}
k_{n}^{+}(u)k_{n+1}^{-}(v)&=k_{n+1}^{-}(v)k_{n}^{+}(u), \\
\frac{(u_{+}-v_{-}-\hbar)(u_{-}-v_{+})}{(u_{-}-v_{+}+\hbar)(u_{+}-v_{-})}k_{n}^{-}(u)k_{n+1}^{+}(v)&=\frac{(u_{-}-v_{+}-\hbar)(u_{+}-v_{-})}{(u_{+}-v_{-}+\hbar)(u_{-}-v_{+})}k_{n+1}^{+}(v)k_{n}^{-}(u),\\
\frac{u_{+}-v_{-}\mp \hbar}{u_{+}-v_{-}}k_{n+1}^{\pm}(u)k_{n+1}^{\mp}(v)&=\frac{u_{-}-v_{+}\mp \hbar}{u_{-}-v_{+}}k_{n+1}^{\mp}(v)k_{n+1}^{\pm}(u).
\end{align}
For $\g_N=\spa_{2n}$ we have
\begin{align}
\frac{u_{-}-v_{+}-2\hbar}{u_{-}-v_{+}-\hbar}k_{n}^{+}(u)k_{n+1}^{-}(v)&=\frac{u_{+}-v_{-}-2\hbar}{u_{+}-v_{-}-\hbar}k_{n+1}^{-}(v)k_{n}^{+}(u), \\
\frac{(u_{+}-v_{-}-2\hbar)(u_{-}-v_{+})}{(u_{-}-v_{+}+\hbar)(u_{+}-v_{-})}k_{n}^{-}(u)k_{n+1}^{+}(v)&=\frac{(u_{-}-v_{+}-2\hbar)(u_{+}-v_{-})}{(u_{+}-v_{-}+\hbar)(u_{-}-v_{+})}k_{n+1}^{+}(v)k_{n}^{-}(u),\\
\frac{u_{+}-v_{-}\mp \hbar}{u_{+}-v_{-}}k_{n+1}^{\pm}(u)k_{n+1}^{\mp}(v)&=\frac{u_{-}-v_{+}\mp \hbar}{u_{-}-v_{+}}k_{n+1}^{\mp}(v)k_{n+1}^{\pm}(u).
\end{align}
For $\g_N=\oa_{2n}$ we have
\begin{align}
\frac{(u_{-}-v_{+})^{2}}{(u_{-}-v_{+}+\hbar)^{2}}k_{n}^{+}(u)k_{n+1}^{-}(v)&=\frac{(u_{+}-v_{-})^{2}}{(u_{+}-v_{-}+\hbar)^{2}}k_{n+1}^{-}(v)k_{n}^{+}(u), \\
\frac{u_{-}-v_{+}+\hbar}{u_{-}-v_{+}-\hbar}k_{n}^{-}(u)k_{n+1}^{+}(v)&=\frac{u_{+}-v_{-}+\hbar}{u_{+}-v_{-}-\hbar}k_{n+1}^{+}(v)k_{n}^{-}(u),\\
\frac{u_{-}-v_{+}\pm \hbar}{u_{-}-v_{+}}k_{n+1}^{\pm}(u)k_{n+1}^{\mp}(v)&=\frac{u_{+}-v_{-}\pm \hbar}{u_{+}-v_{-}}k_{n+1}^{\mp}(v)k_{n+1}^{\pm}(u).
\end{align}
\epr

\bpf For $\g_N=\oa_{2n+1}$, Proposition \ref{pro:Arelations} implies that
\begin{align*}
\frac{u_{-}-v_{+}-\hbar}{u_{-}-v_{+}}k^{+}_{n}(u)k^{-}_{n+1}(v)k^{-}_{n}(v)^{-1}&=\frac{u_{+}-v_{-}-\hbar}{u_{+}-v_{-}}k^{-}_{n+1}(v)k^{-}_{n}(v)^{-1}k^{+}_{n}(u),\\
\frac{u_{+}-v_{-}-\hbar}{u_{+}-v_{-}}k^{+}_{n}(u)k^{-}_{n}(v)&=\frac{u_{-}-v_{+}-\hbar}{u_{-}-v_{+}}k^{-}_{n}(v)k^{+}_{n}(u).
\end{align*}
It follows that $k^{+}_{n}(u)k^{-}_{n+1}(v)=k^{-}_{n+1}(v)k^{+}_{n}(u)$.
Similarly, we have
\begin{align*}
\frac{u_{+}-v_{-}-\hbar}{u_{+}-v_{-}}k^{-}_{n}(u)k^{+}_{n+1}(v)k^{+}_{n}(v)^{-1}&=\frac{u_{-}-v_{+}-\hbar}{u_{-}-v_{+}}k^{+}_{n+1}(v)k^{+}_{n}(v)^{-1}k^{-}_{n}(u),\\
\frac{u_{+}-v_{-}+\hbar}{u_{+}-v_{-}}k^{-}_{n}(u)k^{+}_{n}(v)&=\frac{u_{-}-v_{+}+\hbar}{u_{-}-v_{+}}k^{+}_{n}(v)k^{-}_{n}(u).
\end{align*}
It follows that \begin{displaymath}\frac{u_{+}-v_{-}-\hbar}{u_{+}-v_{-}}\frac{u_{-}-v_{+}}{u_{-}-v_{+}+\hbar}k^{-}_{n}(u)k^{+}_{n+1}(v)=\frac{u_{-}-v_{+}-\hbar}{u_{-}-v_{+}}\frac{u_{+}-v_{-}}{u_{+}-v_{-}+\hbar}k^{+}_{n+1}(v)k^{-}_{n}(u).\end{displaymath}
Proposition \ref{pro:Arelations} also implies that
\begin{align*}
\frac{u_{-}-v_{+}-\hbar}{u_{-}-v_{+}}k^{+}_{n+1}(u)k^{-}_{n+1}(v)k^{-}_{n}(v)^{-1}=\frac{u_{+}-v_{-}-\hbar}{u_{+}-v_{-}}k^{-}_{n+1}(v)k^{-}_{n}(v)^{-1}k^{+}_{n+1}(u).
\end{align*}
Therefore, we get \begin{displaymath}\frac{u_{+}-v_{-}-\hbar}{u_{+}-v_{-}}k^{+}_{n+1}(u)k^{-}_{n+1}(v)=\frac{u_{-}-v_{+}-\hbar}{u_{-}-v_{+}}k^{-}_{n+1}(v)k^{+}_{n+1}(u).\end{displaymath}
This completes the proof for type $B$. The other relations can be checked similarly.
\epf

\subsection{Drinfeld realization}
\label{subsec:tdp}

We now prove the Drinfeld realization for $\DX_\hbar(\mathfrak{g}_N)$. The notation still follows that of Section~\ref{sec:nd}.

The Drinfeld realization for the Yangian in type $B, C, D$ has been proved to be isomorphic to the $RTT$ presentation
\cite{ji:iso} via the Gauss decomposition of the generating matrix. The identification %between the two realization
can also be obtained
through Drinfeld's third presentation as recently given in \cite{GRW1}.
Our method follows the first method.
%recent work of \cite{ji:iso} for the isomorphism between the Drinfeld realization and $RTT$ presentation for
%the Yangian algebra in types $B, C$ and $D$.

\bthm\label{thm:dp}
The extended Yangian double $\DX_\hbar(\mathfrak{g}_N)$ is topologically generated by
the coefficients of the series
$k_i^{\pm}(u)$ ($1\leq i\leq n+1$),
$e_i^{\pm}(u)$ and $f_i^{\pm}(u)$ ($1\leq i\leq n$), and $c$
subject to the following relations, where the indices
run through all admissible values unless specified otherwise. The relations are
\begin{align}
[k_{i}^{\pm}(u),k_{j}^{\pm}(v)]=0.
\end{align}
For $i<j\leqslant n$
\begin{align}
&k_i^+(u)k_j^-(v)=k_j^-(v)k_i^+(u), \\
%\frac{(u_{+}-v_{-}+1)(u_{+}-v_{-}-1)}{(u_{+}-v_{-})^{2}}&k_{j}^{+}(u)k_{i}^{-}(v)=\frac{(u_{-}-v_{+}+1)(u_{-}-v_{+}-1)}{(u_{-}-v_{+})^{2}}k_{i}^{-}(v)k_{j}^{+}(u).
\frac{(u_{+}-v_{-})^2-\hbar^{2}}{(u_{+}-v_{-})^{2}}&k_{j}^{+}(u)k_{i}^{-}(v)=\frac{(u_{-}-v_{+})^2-\hbar^{2}}{(u_{-}-v_{+})^{2}}k_{i}^{-}(v)k_{j}^{+}(u).
\end{align}
For $i\leqslant n-1$ we have
\begin{align}
&k_i^+(u)k_{n+1}^-(v)=k_{n+1}^-(v)k_i^+(u), \\
\frac{(u_{+}-v_{-})^2-\hbar^{2}}{(u_{+}-v_{-})^{2}}&k_{n+1}^{+}(u)k_{i}^{-}(v)=\frac{(u_{-}-v_{+})^2-\hbar^{2}}{(u_{-}-v_{+})^{2}}k_{i}^{-}(v)k_{n+1}^{+}(u).
\end{align}
For $i\leqslant n$ we have
\begin{align}
\frac{u_{+}-v_{-}\mp \hbar}{u_{+}-v_{-}}k_{i}^{\pm}(u)k_{i}^{\mp}(v)&=\frac{u_{-}-v_{+}\mp \hbar}{u_{-}-v_{+}}k_{i}^{\mp}(v)k_{i}^{\pm}(v).
\end{align}
For $i\leqslant n-1, j\leqslant n, j\neq i,i+1$, and $\epsilon=\pm$ or $\pm\hbar$ we have
\begin{align}
k_{i}^{\pm}(u)^{-1}X_{i}^{\epsilon}(v)k_{i}^{\pm}(u)&=\frac{u_{\mp\epsilon}-v-\hbar}{u_{\mp\epsilon}-v}X_{i}^{\epsilon}(v), \\
%k_{i}^{\pm}(u)X_{i}^{-}(v)k_{i}^{\pm}(u)^{-1}&=\frac{u_{\pm}-v-1}{u_{\pm}-v}X_{i}^{-}(v), \\
k_{i+1}^{\pm}(u)^{-1}X_{i}^{\epsilon}(v)k_{i+1}^{\pm\epsilon}(u)&=\frac{u_{\mp\epsilon}-v+\hbar}{u_{\mp\epsilon}-v}X_{i}^{\epsilon}(v), \\
%k_{i+1}^{\pm}(u)X_{i}^{-}(v)k_{i+1}^{\pm}(u)^{-1}&=\frac{u_{\pm}-v-1}{u_{\pm}-v}X_{i}^{-}(v), \\
k_{j}^{\pm}(u)^{-\epsilon}X_{i}^{\epsilon}(v)k_{j}^{\pm}(u)^{\epsilon}&=X_{i}^{\epsilon}(v), \\
%k_{j}^{\pm}(u)X_{i}^{-}(v)k_{j}^{\pm}(u)^{-1}&=X_{i}^{-}(v), \\
(u-v\pm \hbar)X_{i}^{\pm}(u)X_{i}^{\pm}(v)&=(u-v\mp \hbar)X_{i}^{\pm}(v)X_{i}^{\pm}(u), \\
(u-v-(\hbar+\epsilon)/2)X_{i-1}^{\epsilon}(u)X_{i}^{\epsilon}(v)&=(u-v-(\hbar-\epsilon)/2)X_{i}^{\epsilon}(v)X_{i-1}^{\epsilon}(u). %\\
%(u-v)X_{i-1}^{-}(u)X_{i}^{-}(v)&=(u-v-1)X_{i}^{-}(v)X_{i-1}^{-}(u).
\end{align}
For $i,j\leqslant n-1$ we have
\begin{align}
&[X_{i}^{+}(u),X_{j}^{-}(v)]=\delta_{ij}\hbar\{\delta(u_{-}-v_{+})k_{i+1}^{-}(v_{+})k_{i}^{-}(v_{+})^{-1}-\delta(u_{+}-v_{-})k_{i+1}^{+}(u_{+})k_{i}^{+}(u_{+})^{-1}\}\\
&where ~\delta(u-v)=\Sigma_{k\in\mathbb{Z}}u^{-k-1}v^{k}.\nonumber
\end{align}
For $\g_N=\oa_{2n+1}$ we have
\begin{align}
k_{n}^{+}(u)k_{n+1}^{-}(v)&=k_{n+1}^{-}(v)k_{n}^{+}(u), \\
\frac{(u_{+}-v_{-}-\hbar)(u_{-}-v_{+})}{(u_{-}-v_{+}+\hbar)(u_{+}-v_{-})}k_{n}^{-}(u)k_{n+1}^{+}(v)&=\frac{(u_{-}-v_{+}-\hbar)(u_{+}-v_{-})}{(u_{+}-v_{-}+\hbar)(u_{-}-v_{+})}k_{n+1}^{+}(v)k_{n}^{-}(u),\\
\frac{u_{+}-v_{-}\mp \hbar}{u_{+}-v_{-}}k_{n+1}^{\pm}(u)k_{n+1}^{\mp}(v)&=\frac{u_{-}-v_{+}\mp \hbar}{u_{-}-v_{+}}k_{n+1}^{\mp}(v)k_{n+1}^{\pm}(u), \\
k_{i}^{\pm}(u)X_{n}^{+}(v)=X_{n}^{+}(v)k_{i}^{\pm}(u)&,\quad k_{i}^{\pm}(u)X_{n}^{-}(v)=X_{n}^{-}(v)k_{i}^{\pm}(u), (i\leq n-1) \\
k_{n}^{\pm}(u)^{-\epsilon}X_{n}^{\epsilon}(v)k_{n}^{\pm}(u)^{\epsilon}&=\frac{u_{\mp\epsilon}-v-\hbar}{u_{\mp\epsilon}-v}X_{n}^{\epsilon}(v), \\
%k_{n}^{\pm}(u)X_{n}^{-}(v)k_{n}^{\pm}(u)^{-1}&=\frac{u_{\pm}-v-1}{u_{\pm}-v}X_{n}^{-}(v),
k_{n+1}^{\pm}(u)^{-\epsilon}X_{j}^{\epsilon}(v)k_{n+1}^{\pm}(u)^{\epsilon}&=X_{j}^{\epsilon}(v),\quad (j<n)\\
%k_{n+1}^{\pm}(u)X_{j}^{-}(v)k_{n+1}^{\pm}(u)^{-1}&=X_{j}^{-}(v),\quad (j<n)\\
k_{n+1}^{\pm}(u)^{-\epsilon}X_{n}^{\epsilon}(v)k_{n+1}^{\pm}(u)^{\epsilon}&=\frac{(u_{\mp\epsilon}-v-\hbar)(u_{\mp\epsilon}-v+\frac{1}{2}\hbar)}
{(u_{\mp\epsilon}-v)(u_{\mp\epsilon}-v-\frac{1}{2}\hbar)}X_{n}^{\epsilon}(v), \\
%k_{n+1}^{\pm}(u)X_{n}^{-}(v)k_{n+1}^{\pm}(u)^{-1}&=\frac{(u_{\pm}-v-1)(u_{\pm}-v+\frac{1}{2})}{(u_{\pm}-v)(u_{\pm}-v-\frac{1}{2})}X_{n}^{-}(v),  \\
(u-v\pm\frac{1}{2}\hbar)X_{n}^{\pm}(u)X_{n}^{\pm}(v)&=(u-v\mp\frac{1}{2}\hbar)X_{n}^{\pm}(v)X_{n}^{\pm}(u), \\
%(u-v-\frac{1}{2})X_{n}^{-}(u)X_{n}^{-}(v)&=(u-v+\frac{1}{2})X_{n}^{-}(v)X_{n}^{-}(u), \\
[X_{n-1}^{\pm}(u),X_{n}^{\mp}(v)]&=[X_{n}^{\pm}(u),X_{n-1}^{\mp}(v)]=0, \\
[X_{i}^{\pm}(u),X_{n}^{\pm}(v)]&=0,\quad (i<n-1) \\
(u-v-(\hbar+\epsilon)/2)X_{n-1}^{\epsilon}(u)X_{n}^{\epsilon}(v)&=(u-v-(\hbar-\epsilon)/2)X_{n}^{\epsilon}(v)X_{n-1}^{\epsilon}(u), \\
%(u-v)X_{n-1}^{-}(u)X_{n}^{-}(v)&=(u-v-1)X_{n}^{-}(v)X_{n-1}^{-}(u), \\
[X_{n}^{+}(u),X_{n}^{-}(v)]&=\hbar\{\delta(u_{-}-v_{+})k_{n+1}^{-}(v_{+})k_{n}^{-}(v_{+})^{-1}\\
&-\delta(u_{+}-v_{-})k_{n+1}^{+}(u_{+})k_{n}^{+}(u_{+})^{-1}\}.\nonumber
\end{align}
For $\g_N=\spa_{2n}$ we have
\begin{align}
\frac{u_{-}-v_{+}-2\hbar}{u_{-}-v_{+}-\hbar}k_{n}^{+}(u)k_{n+1}^{-}(v)&=\frac{u_{+}-v_{-}-2\hbar}{u_{+}-v_{-}-\hbar}k_{n+1}^{-}(v)k_{n}^{+}(u), \\
\frac{(u_{+}-v_{-}-2\hbar)(u_{-}-v_{+})}{(u_{-}-v_{+}+\hbar)(u_{+}-v_{-})}k_{n}^{-}(u)k_{n+1}^{+}(v)&=\frac{(u_{-}-v_{+}-2\hbar)(u_{+}-v_{-})}{(u_{+}-v_{-}+\hbar)(u_{-}-v_{+})}k_{n+1}^{+}(v)k_{n}^{-}(u),\\
\frac{u_{+}-v_{-}\mp \hbar}{u_{+}-v_{-}}k_{n+1}^{\pm}(u)k_{n+1}^{\mp}(v)&=\frac{u_{-}-v_{+}\mp \hbar}{u_{-}-v_{+}}k_{n+1}^{\mp}(v)k_{n+1}^{\pm}(u), \\
k_{i}^{\pm}(u)X_{n}^{\epsilon}(v)&=X_{n}^{\epsilon}(v)k_{i}^{\pm}(u), \quad (i\leq n-1)
%\quad k_{i}^{\pm}(u)X_{n}^{-}(v)=X_{n}^{-}(v)k_{i}^{\pm}(u),
\\
k_{n}^{\pm}(u)^{-\epsilon}X_{n}^{\epsilon}(v)k_{n}^{\pm}(u)^{\epsilon}&=\frac{u_{\mp\epsilon}-v-2\hbar}{u_{\mp\epsilon}-v}X_{n}^{\epsilon}(v), \\
%k_{n}^{\pm}(u)X_{n}^{-}(v)k_{n}^{\pm}(u)^{-1}&=\frac{u_{\pm}-v-2}{u_{\pm}-v}X_{n}^{-}(v), \\
k_{n+1}^{\pm}(u)^{-\epsilon}X_{j}^{\epsilon}(v)k_{n+1}^{\pm}(u)^{\epsilon}&=X_{j}^{\epsilon}(v),\quad (j<n-1)\\
%k_{n+1}^{\pm}(u)X_{j}^{-}(v)k_{n+1}^{\pm}(u)^{-1}&=X_{j}^{-}(v),\quad (j<n-1)\\
k_{n+1}^{\pm}(u)^{-\epsilon}X_{n-1}^{\epsilon}(v)k_{n+1}^{\pm}(u)^{\epsilon}&=\frac{u_{\mp\epsilon}-v+\hbar}{u_{\mp\epsilon}-v+2\hbar}X_{n-1}^{\epsilon}(v), \\
%k_{n+1}^{\pm}(u)X_{n-1}^{-}(v)k_{n+1}^{\pm}(u)^{-1}&=\frac{u_{\pm}-v+1}{u_{\pm}-v+2}X_{n-1}^{-}(v),
k_{n+1}^{\pm}(u)^{-\epsilon}X_{n}^{\epsilon}(v)k_{n+1}^{\pm}(u)^{\epsilon}&=\frac{u_{\mp}-v+2\hbar}{u_{\mp}-v}X_{n}^{\epsilon}(v), \\
%k_{n+1}^{\pm}(u)X_{n}^{-}(v)k_{n+1}^{\pm}(u)^{-1}&=\frac{u_{\pm}-v+2}{u_{\pm}-v}X_{n}^{-}(v), \\
(u-v\pm 2\hbar)X_{n}^{\pm}(u)X_{n}^{\pm}(v)&=(u-v\mp 2\hbar)X_{n}^{\pm}(v)X_{n}^{\pm}(u),
%(u-v-2)X_{n}^{-}(u)X_{n}^{-}(v)&=(u-v+2)X_{n}^{-}(v)X_{n}^{-}(u), \\
\end{align}
\begin{align}
[X_{n-1}^{\pm}(u),X_{n}^{\mp}(v)]&=[X_{n}^{\pm}(u),X_{n-1}^{\mp}(v)]=0, \\
[X_{i}^{\pm}(u),X_{n}^{\pm}(v)]&=0,\quad (i<n-1) \\
(u-v-(\hbar+\epsilon))X_{n-1}^{\epsilon}(u)X_{n}^{\epsilon}(v)&=(u-v-(\hbar-\epsilon))X_{n}^{\epsilon}(v)X_{n-1}^{\epsilon}(u), \\
%(u-v)X_{n-1}^{-}(u)X_{n}^{-}(v)&=(u-v-2)X_{n}^{-}(v)X_{n-1}^{-}(u), \\
[X_{n}^{+}(u),X_{n}^{-}(v)]&=2\hbar\{\delta(u_{-}-v_{+})k_{n+1}^{-}(v_{+})k_{n}^{-}(v_{+})^{-1}\\
&-\delta(u_{+}-v_{-})k_{n+1}^{+}(u_{+})k_{n}^{+}(u_{+})^{-1}\}.\nonumber
\end{align}
For $\g_N=\oa_{2n}$ we have
\begin{align}
\frac{(u_{-}-v_{+})^{2}}{(u_{-}-v_{+}+\hbar)^{2}}k_{n}^{+}(u)k_{n+1}^{-}(v)&=\frac{(u_{+}-v_{-})^{2}}{(u_{+}-v_{-}+\hbar)^{2}}k_{n+1}^{-}(v)k_{n}^{+}(u), \\
\frac{u_{-}-v_{+}+\hbar}{u_{-}-v_{+}-\hbar}k_{n}^{-}(u)k_{n+1}^{+}(v)&=\frac{u_{+}-v_{-}+\hbar}{u_{+}-v_{-}-\hbar}k_{n+1}^{+}(v)k_{n}^{-}(u),\\
\frac{u_{-}-v_{+}\pm \hbar}{u_{-}-v_{+}}k_{n+1}^{\pm}(u)k_{n+1}^{\mp}(v)&=\frac{u_{+}-v_{-}\pm \hbar}{u_{+}-v_{-}}k_{n+1}^{\mp}(v)k_{n+1}^{\pm}(u), \\
k_{i}^{\pm}(u)X_{n}^{\epsilon}(v)&=X_{n}^{\epsilon}(v)k_{i}^{\pm}(u),   \quad (i\leq n-2) \\
%\quad k_{i}^{\pm}(u)X_{n}^{-}(v)=X_{n}^{-}(v)k_{i}^{\pm}(u), (i\leq n-2) \\
k_{n-1}^{\pm}(u)^{-\epsilon}X_{n}^{\epsilon}(v)k_{n-1}^{\pm}(u)^{\epsilon}&=\frac{u_{\mp\epsilon}-v-\hbar}{u_{\mp\epsilon}-v}X_{n}^{\epsilon}(v), \\
%k_{n-1}^{\pm}(u)X_{n}^{-}(v)k_{n-1}^{\pm}(u)^{-1}&=\frac{u_{\pm}-v-1}{u_{\pm}-v}X_{n}^{-}(v), \\
k_{n}^{\pm}(u)^{-\epsilon}X_{n}^{\epsilon}(v)k_{n}^{\pm}(u)^{\epsilon}&=\frac{u_{\mp\epsilon}-v-\hbar}{u_{\mp\epsilon}-v}X_{n}^{\epsilon}(v),
%k_{n}^{\pm}(u)X_{n}^{-}(v)k_{n}^{\pm}(u)^{-1}&=\frac{u_{\pm}-v-1}{u_{\pm}-v}X_{n}^{-}(v),  \\
\end{align}
\begin{align}
k_{n+1}^{\pm}(u)^{-\epsilon}X_{j}^{\epsilon}(v)k_{n+1}^{\pm}(u)^{\epsilon}&=X_{j}^{\epsilon}(v),\quad (j<n-1)\\
%k_{n+1}^{\pm}(u)X_{j}^{-}(v)k_{n+1}^{\pm}(u)^{-1}&=X_{j}^{-}(v),\quad (j<n-1)\\
k_{n+1}^{\pm}(u)^{-\epsilon}X_{n-1}^{\epsilon}(v)k_{n+1}^{\pm}(u)^{\epsilon}&=\frac{u_{\mp\epsilon}-v-\hbar}{u_{\mp\epsilon}-v}X_{n-1}^{\epsilon}(v), \\
%k_{n+1}^{\pm}(u)X_{n-1}^{-}(v)k_{n+1}^{\pm}(u)^{-1}&=\frac{u_{\pm}-v-1}{u_{\pm}-v}X_{n-1}^{-}(v),
k_{n+1}^{\pm}(u)^{-\epsilon}X_{n}^{\epsilon}(v)k_{n+1}^{\pm}(u)^{\epsilon}&=\frac{u_{\mp\epsilon}-v+\hbar}{u_{\mp\epsilon}-v}X_{n}^{\epsilon}(v), \\
%k_{n+1}^{\pm}(u)X_{n}^{-}(v)k_{n+1}^{\pm}(u)^{-1}&=\frac{u_{\pm}-v+1}{u_{\pm}-v}X_{n}^{-}(v),  \\
(u-v\pm \hbar)X_{n}^{\pm}(u)X_{n}^{\pm}(v)&=(u-v\mp \hbar)X_{n}^{\pm}(v)X_{n}^{\pm}(u), \\
%(u-v-1)X_{n}^{-}(u)X_{n}^{-}(v)&=(u-v+1)X_{n}^{-}(v)X_{n}^{-}(u), \\
[X_{n-1}^{\pm}(u),X_{n}^{\mp}(v)]&=[X_{n}^{\pm}(u),X_{n-1}^{\mp}(v)]=0, \\
X_{n-1}^{\pm}(u)X_{n}^{\pm}(v)&=X_{n}^{\pm}(v)X_{n-1}^{\pm}(u), \\
%X_{n-1}^{-}(u)X_{n}^{-}(v)&=X_{n}^{-}(v)X_{n-1}^{-}(u), \\
(u-v-(\hbar+\epsilon))X_{n-2}^{\epsilon}(u)X_{n}^{\epsilon}(v)&=(u-v-(\hbar-\epsilon))X_{n}^{\epsilon}(v)X_{n-2}^{\epsilon}(u), \\
%(u-v)X_{n-2}^{-}(u)X_{n}^{-}(v)&=(u-v-1)X_{n}^{-}(v)X_{n-2}^{-}(u), \\
[X_{i}^{\pm}(u),X_{n}^{\pm}(v)]&=0,\quad (i<n-2) \\
[X_{n}^{+}(u),X_{n}^{-}(v)]=\hbar\{\delta(u_{-}-v_{+})k_{n+1}^{-}(v_{+})k_{n-1}^{-}(v_{+})^{-1}&
-\delta(u_{+}-v_{-})k_{n+1}^{+}(u_{+})k_{n-1}^{+}(u_{+})^{-1}\}.
\end{align}
In all three cases we have
\begin{align}
[X_{i}^{+}(u),X_{n}^{-}(v)]&=[X_{n}^{+}(u),X_{i}^{-}(v)]=0.\quad (i<n)
\end{align}
and the Serre relations are
\begin{align}
\sum_{\sigma\in\mathfrak{S}_{m}}[X_{i}^{\pm}(u_{\sigma(1)}),[X_{i}^{\pm}(u_{\sigma(2)})\cdots,[X_{i}^{\pm}(u_{\sigma(m)}),X_{j}^{\pm}(v)]\cdots]=0,
~~ i\neq j,m=1-a_{ij}.
\end{align}
Here $A=(a_{ij})$ is the Cartan matrix of the Lie algebra $\oa_{2n+1},\spa_{2n},\oa_{2n}$ respectively.
\ethm

\bpf
Apart from the Serre relations, all the relations are satisfied in the algebra $\DX_\hbar(\g_{N})$ due to
Propositions~\ref{cor:commu}, \ref{pro:rellowrankB}, \ref{pro:rellowrankC}, \ref{pro:rellowrankD}, \ref{prop:BCDkn+1xn-1}, \ref{pro:relxnxi}, \ref{prop:Dxn-2xn},
\ref{prop:BCDxixn}, \ref{prop:BCDknkn+1} and Corollary~\ref{cor:Bxn-1xn}, \ref{cor:Cxn-1xn}.  To see the Serre relations, we need to write down the commutation relations between $X^{\pm}_{i}(u)$ and $X^{\pm}_{j}(v)$, which have been calculated in the previous section. Then it is not difficult to verify the Serre relations.

Now let $\wh \DX_\hbar(\g_{N})$ be the algebra
with generators and relations as in the statement of the theorem.
The above argument implies that there is a homomorphism
\beql{surjhom}
\wh \DX_\hbar(\g_{N})\to\DX_\hbar(\g_{N})
\eeq
which takes
the generators $k_{i}^{(r)}$,
$e_{i}^{(r)}$, $f_{i}^{(r)}$ and $c$ of $\wh \DX_\hbar(\g_{N})$ to the corresponding elements
of $\DX_\hbar(\g_{N})$. If we show that
this homomorphism is bijective, the proof is done. To prove the surjectivity we need a lemma.
\epf

\ble\label{lem:eijr}
Assume $r\in\mathbb{Z}^{\times}$. For the case $\DX_\hbar(\oa_{2n+1})$, if $1\leqslant i<j\leqslant n$, we have
\begin{alignat}{2}
e_{i\ts j+1}^{(r)}&=\big[e_{i\tss j}^{(r)},\ts e_j^{(1)}\big],\qquad
&f_{j+1\ts i}^{(r)}&=\big[f_j^{(1)},\ts f_{j\tss i}^{(r)}\big],
\non\\[0.3em]
e_{i\ts j'}^{(r)}&=-\big[e_{i\tss j'-1}^{(r)},\ts e_{j}^{(1)}\big],\qquad
&f_{j'\ts i}^{(r)}&=-\big[f_{j}^{(1)},\ts f_{j'-1\tss i}^{(r)}\big].
\non
\end{alignat}
For the case $\DX_\hbar(\spa_{2n})$, if $1\leqslant i<j\leqslant n-1$, we have
\begin{alignat}{2}
e_{i\ts j+1}^{(r)}&=\big[e_{i\tss j}^{(r)},\ts e_j^{(1)}\big],\qquad
&f_{j+1\ts i}^{(r)}&=\big[f_j^{(1)},\ts f_{j\tss i}^{(r)}\big],
\non\\[0.3em]
e_{i\ts j'}^{(r)}&=-\big[e_{i\tss j'-1}^{(r)},\ts e_{j}^{(1)}\big],\qquad
&f_{j'\ts i}^{(r)}&=-\big[f_{j}^{(1)},\ts f_{j'-1\tss i}^{(r)}\big].
\non
\end{alignat}
Furthermore, if $1\leqslant i\leqslant n-1$, we have
\begin{alignat}{2}
e_{i\ts n'}^{(r)}&=\frac{1}{2}\ts\big[e_{i\tss n'-1}^{(r)},\ts e_{n}^{(1)}\big],\qquad
&f_{n'\ts i}^{(r)}&=\frac{1}{2}\ts\big[f_{n}^{(1)},\ts f_{n'-1\tss i}^{(r)}\big],
\non\\[0.3em]
e_{i\ts i'}^{(r)}&=-\big[e_{i\tss i'-1}^{(r)},\ts e_{i}^{(1)}\big],\qquad
&f_{i'\ts i}^{(r)}&=-\big[f_{i}^{(1)},\ts f_{i'-1\tss i}^{(r)}\big].
\non
\end{alignat}
For the case $\DX_\hbar(\oa_{2n})$, if $1\leqslant i<j\leqslant n-1$, we have
\begin{alignat}{2}
e_{i\ts j+1}^{(r)}&=\big[e_{i\tss j}^{(r)},\ts e_j^{(1)}\big],\qquad
&f_{j+1\ts i}^{(r)}&=\big[f_j^{(1)},\ts f_{j\tss i}^{(r)}\big],
\non\\[0.3em]
e_{i\ts j'}^{(r)}&=-\big[e_{i\tss j'-1}^{(r)},\ts e_{j}^{(1)}\big],\qquad
&f_{j'\ts i}^{(r)}&=-\big[f_{j}^{(1)},\ts f_{j'-1\tss i}^{(r)}\big].
\non
\end{alignat}
And if $1\leqslant i\leqslant n-2$, we have
\ben
e_{i\ts n'}^{(r)}=\big[e_{i\tss n-1}^{(r)},\ts e_{n}^{(1)}\big],\qquad
f_{n'\ts i}^{(r)}=\big[f_{n}^{(1)},\ts f_{n-1\tss i}^{(r)}\big].
\een
\ele

\bpf
All relations can be checked with the use of the Gauss decomposition \eqref{gd}
and defining relations \eqref{defrel++}-\eqref{defrel+-}. For example, let $\g_N=\spa_{2n},r<0$. By taking the coefficients of $v^{-1}$
on both sides of \eqref{defrel+-}, we get
$[t^{-}_{1j}(u),t_{j\tss j+1}^{(1)}]=t^{-}_{1\tss j+1}(u)$ for $1<j\leqslant n-1$. The Gauss decomposition \eqref{gd} implies that
$k^{-}_{1}(u)[e^{-}_{1j}(u),e_{j}^{(1)}]=k^{-}_{1}(u)e^{-}_{1\tss j+1}(u)$. It follows that $[e^{-}_{1j}(u),e_{j}^{(1)}]=e^{-}_{1\tss j+1}(u)$. Hence,
$e_{1\ts j+1}^{(r)}=[e_{1\tss j}^{(r)},\ts e_j^{(1)}]$ for $r<0$.

The similar argument gives that
$e_{1\ts j'}^{(r)}=[e_{1\tss j'-1}^{(r)},\ts e_{(j+1)'\tss j'}^{(1)}]$ for $1\leqslant j\leqslant n-1$ and $r<0$.
By Proposition~\ref{prop:eiei'}, this equality can be written as
$e_{1\ts j'}^{(r)}=-[e_{1\tss j'-1}^{(r)},\ts e_{j}^{(1)}]$.
We can verify the remaining cases with $i=1$ in the same way.
The proof for arbitrary values of $i$ can be done by
applying Corollary~\ref{cor:guass-embed}.
\epf

Let $\wt\DX_\hbar(\g_N)$ be the subalgebra of $\DX_\hbar(\g_{N})$ generated by the coefficients of the series
$k^{\pm}_i(u)$ with $i=1,\dots,n+1$, and
$e^{\pm}_j(u)$, $f^{\pm}_j(u)$ with $j=1,\dots,n$ and $c$.
By Lemma~\ref{lem:eijr}, for $r\in\mathbb{Z}^{\times}$, $i,j$ satisfy the conditions $i<j$ and $i<j\pr$ in the orthogonal case,
$i\leqslant j\pr$ in the symplectic case,
all elements $e_{ij}^{(r)}$ and $f_{ji}^{(r)}$ belong to the subalgebra $\wt\DX_\hbar(\g_N)$. Therefore,
the Gauss decomposition \eqref{gd} implies that the coefficients
of the series $t^{\pm}_{ij}(u)$ also belong to the subalgebra $\wt\DX_\hbar(\g_N)$ if the indices $i,j$ satisfy the above conditions.
Moreover, Theorem~\ref{thm:Center} implies that the coefficients of the series $z^{\pm}_N(u)$
are also in $\wt\DX_\hbar(\g_N)$. Finally, taking the coefficients of $u^{r}$
for $r=-1,-2,\dots$ or $r=0,1,\dots$ in \eqref{zcenter} and using induction on $r$, we can prove
the coefficients of all series $t^{\pm}_{ij}(u)$ belong to $\wt\DX_\hbar(\g_N)$.
So $\wt\DX_\hbar(\g_N)=\DX_\hbar(\g_N)$.
This completes the proof of the surjectivity.

Next we show that this homomorphism is injective. Under the isomorphism \eqref{isomapping},
the images of the elements $e_{ij}^{(r)}$
and $f_{ji}^{(r)}$ in the $(r-1)$-th component of the graded
algebra $\gr \DX_\hbar(\g_{N})$ respectively correspond to
$F_{ij}\tss x^{r-1}$ and $F_{ji}\tss x^{r-1}$, the images of the elements $e_{ij}^{(-r)}$
and $f_{ji}^{(-r)}$ in the $(-r)$-th component of the graded
algebra $\gr \DX_\hbar(\g_{N})$ respectively correspond to
$F_{ij}\tss x^{-r}$ and $F_{ji}\tss x^{-r}$. Moreover,
the image of $h_{i}^{(r)}$ correspond to
$F_{ii}\tss x^{r-1}+\ze_r/2$ for $i=1,\dots, n$, while the image of $h_{i}^{(-r)}$ correspond to
$F_{ii}\tss x^{-r}+\vs_r/2$ for $i=1,\dots, n$. \eqref{znim} and Theorem~\ref{thm:Center} implies that
\ben
\bar h^{(r)}_{n+1}\mapsto
\begin{cases}
\ze_r/2
\qquad&\text{for}\quad \oa_{2n+1}\\[0.2em]
-F_{n\tss n}\tss x^{r-1}+\ze_r/2
\qquad&\text{for}\quad \spa_{2n}\\[0.2em]
-F_{n-1\tss n-1}\tss x^{r-1}-F_{n\tss n}\tss x^{r-1}+\ze_r/2
\qquad&\text{for}\quad \oa_{2n},
\end{cases}
\een
\ben
\bar h^{(-r)}_{n+1}\mapsto
\begin{cases}
\vs_r/2
\qquad&\text{for}\quad \oa_{2n+1}\\[0.2em]
-F_{n\tss n}\tss x^{-r}+\vs_r/2
\qquad&\text{for}\quad \spa_{2n}\\[0.2em]
-F_{n-1\tss n-1}\tss x^{-r}-F_{n\tss n}\tss x^{-r}+\vs_r/2
\qquad&\text{for}\quad \oa_{2n}.
\end{cases}
\een
By the Poincar\'e-Birkhoff--Witt theorem for $\U(\g_{N}[x,x^{-1}]\oplus \mathbb{C}K)$, we conclude that the ordered monomials in the set of generators
$k_{i}^{(r)},k_{i}^{(-r)}$
with $i=1,\dots,n+1,r\geq 1$,
and $e_{ij}^{(r)},e_{ij}^{(-r)},f_{ji}^{(r)},f_{ji}^{(-r)}$ with $r\geq 1$
and the conditions $i<j$ and $i<j\pr$ in the orthogonal case, $i\leqslant j\pr$ in the symplectic case, and $c$, are linearly independent in the algebra $\DX_\hbar(\g_{N})$.

Now let $r\in\mathbb{Z}^{\times}$, define elements $e_{ij}^{(r)}$ and $f_{ji}^{(r)}$ of $\wh \DX_\hbar(\g_{N})$
inductively as follows. For the case $\g_N=\oa_{2n+1}$, we set $e_{i\ts i+1}^{(r)}=e_{i}^{(r)}$ and
$f_{i+1\ts i}^{(r)}=f_{i}^{(r)}$, and
\begin{alignat}{2}\label{Beijpone}
e_{i\ts j+1}^{(r)}&=[e_{i\tss j}^{(r)},\ts e_j^{(1)}],\qquad
&f_{j+1\ts i}^{(r)}&=[f_j^{(1)},\ts f_{j\tss i}^{(r)}],\\
e_{i\ts j'}^{(r)}&=-[e_{i\tss j'-1}^{(r)},\ts e_{j}^{(1)}],\qquad
&f_{j'\ts i}^{(r)}&=-[f_{j}^{(1)},\ts f_{j'-1\tss i}^{(r)}]
\non
\end{alignat}
if $1\leqslant i<j\leqslant n$. For the case $\g_N=\spa_{2n}$,
we set $e_{i\ts i+1}^{(r)}=e_{i}^{(r)}$ and
$f_{i+1\ts i}^{(r)}=f_{i}^{(r)}$, and
\begin{alignat}{2}\label{Ceijpone}
e_{i\ts j+1}^{(r)}&=[e_{i\tss j}^{(r)},\ts e_j^{(1)}],\qquad
&f_{j+1\ts i}^{(r)}&=[f_j^{(1)},\ts f_{j\tss i}^{(r)}],\\
e_{i\ts j'}^{(r)}&=-[e_{i\tss j'-1}^{(r)},\ts e_{j}^{(1)}],\qquad
&f_{j'\ts i}^{(r)}&=-[f_{j}^{(1)},\ts f_{j'-1\tss i}^{(r)}]
\non
\end{alignat}
if $1\leqslant i< j\leqslant n-1$. Moreover, we set
$e_{n\ts n'}^{(r)}=e_n^{(r)}$ and $f_{n'\ts n}^{(r)}=f_n^{(r)}$, and
\begin{alignat}{2}\label{Ceijponen}
&e_{i\ts n'}^{(r)}=\frac{1}{2}\ts[e_{i\tss n'-1}^{(r)},\ts e_{n}^{(1)}],\qquad
f_{n'\ts i}^{(r)}=\frac{1}{2}\ts[f_{n}^{(1)},\ts f_{n'-1\tss i}^{(r)}],\\
&e_{i\ts i'}^{(r)}=-[e_{i\tss i'-1}^{(r)},\ts e_{i}^{(1)}],\qquad
f_{i'\ts i}^{(r)}=-[f_{i}^{(1)},\ts f_{i'-1\tss i}^{(r)}]
\non
\end{alignat}
if $1\leqslant i\leqslant n-1$. For the case $\g_N=\oa_{2n}$,
we set $e_{i\ts i+1}^{(r)}=e_{i}^{(r)}$ and
$f_{i+1\ts i}^{(r)}=f_{i}^{(r)}$, and
\begin{alignat}{2}\label{Deijpone}
e_{i\ts j+1}^{(r)}&=[e_{i\tss j}^{(r)},\ts e_j^{(1)}],\qquad
&f_{j+1\ts i}^{(r)}&=[f_j^{(1)},\ts f_{j\tss i}^{(r)}],\\
e_{i\ts j'}^{(r)}&=-[e_{i\tss j'-1}^{(r)},\ts e_{j}^{(1)}],\qquad
&f_{j'\ts i}^{(r)}&=-[f_{j}^{(1)},\ts f_{j'-1\tss i}^{(r)}]
\non
\end{alignat}
if $1\leqslant i< j\leqslant n-1$. Moreover, we set
$e_{n-1\ts n'}^{(r)}=e_n^{(r)}$ and $f_{n'\ts n-1}^{(r)}=f_n^{(r)}$, and
\beql{Deijponen}
e_{i\ts n'}^{(r)}=[e_{i\tss n-1}^{(r)},\ts e_{n}^{(1)}],\qquad
f_{n'\ts i}^{(r)}=[f_{n}^{(1)},\ts f_{n-1\tss i}^{(r)}]
\eeq
if $1\leqslant i\leqslant n-2$.

Lemma~\ref{lem:eijr} implies that the images of the elements
$e_{ij}^{(r)}$ and $f_{ji}^{(r)}$ of the algebra
$\wh \DX_\hbar(\g_{N})$ under the homomorphism \eqref{surjhom}
coincide with the elements of $\DX_\hbar(\g_{N})$ denoted by the same symbols.

Now we will show that \eqref{surjhom} is injective. First of all, we introduce some notations. Denote by $\wh \Ec$, $\wh \Fc$
and $\wh \Hc$ the subalgebras of $\wh \DX_\hbar(\g_{N})$ respectively
generated by all elements
of the form $e_{i}^{(r)}$, $f_{i}^{(r)}$ and $k_{i}^{(r)}.$
Denote by $\wh \Ec^{+}$, $\wh \Fc^{+}$
and $\wh \Hc^{+}$ the subalgebras of $\wh \DX_\hbar(\g_{N})$ respectively
generated by all elements
of the form $e_{i}^{(r)}$, $f_{i}^{(r)}$ and $k_{i}^{(r)}$ with $r>0.$
Denote by $\wh \Ec^{-}$, $\wh \Fc^{-}$
and $\wh \Hc^{-}$ the subalgebras of $\wh \DX_\hbar(\g_{N})$ respectively
generated by all elements
of the form $e_{i}^{(r)}$, $f_{i}^{(r)}$ and $k_{i}^{(r)}$ with $r<0.$
Define an ascending filtration
on $\wh \Ec^{-}$ by setting $\deg e_{i}^{(-r)}=-r, \deg \hbar=0$.
Denote by $\gr\wh \Ec^{-}$ the corresponding graded algebra.
Let $\eb_{ij}^{\tss(-r)}$ be the image of $e_{ij}^{(-r)}$ in the
$(-r)$-th component of the graded algebra $\gr\wh \Ec^{-}$. By using the skew-symmetry conditions
\beql{skew-symmetry}
\eb_{i\tss j}^{\tss(-r)}=-\theta_{ij}\eb_{j'\tss i'}^{\tss(-r)},
\eeq
we can extend the range of subscripts of
$\eb_{ij}^{\tss(-r)}$ to all values $1\leqslant i<j\leqslant 1'$. Note that the algebra $\wh \Ec^{+}$ is spanned
by the set of monomials in the elements $e_{ij}^{\tss(r)}$ taken in some fixed order; see \cite[Sec.~5.5]{ji:iso}.
Similarly, the desired spanning property of
the algebra $\wh \Ec^{-}$ follows from the relations
\beql{TypeBeijeklbar-}
[\eb_{i\tss j}^{\tss(-r)},\eb_{k\tss l}^{\tss(-s)}]=
\de_{k\tss j}\ts\eb_{i\tss l}^{\tss(-r-s)}-\de_{i\tss l}\ts\eb_{kj}^{\tss(-r-s)}
-\theta_{ij}\tss \de_{k\tss i'}\ts\eb_{j'\tss l}^{\tss(-r-s)}
+\theta_{ij}\tss \de_{j'\tss l}\ts\eb_{k\tss i'}^{\tss(-r-s)}.
\eeq
We can use the same method in \cite[Sec.~5.5]{ji:iso} to prove \eqref{TypeBeijeklbar-}. In addition, we can swap $e_{ij}^{\tss(r)}$ and $e_{ij}^{\tss(-r)}$ by defining relations in Theorem~\ref{thm:dp}. Therefore, the algebra $\wh \Ec$ is spanned by the ordered monomials in the elements $e_{ij}^{\tss(r)}$ and $e_{ij}^{\tss(-r)}$. The same is true for $\wh \Fc$. Observe that the ordered monomials in $k_i^{(r)}$ span $\wh\Hc$. Moreover, the defining relations of $\wh \DX_\hbar(\g_{N})$ implies that the multiplication map
\beq
\wh\Fc~\widetilde{\otimes}~\wh\Hc~\widetilde{\otimes}~\wh\Ec~\widetilde{\otimes}~\Ac c  \to \wh \DX_\hbar(\g_{N})
\eeq
is surjective. Therefore, if we let
the elements of $\wh\Fc$ precede the elements
of $\wh\Hc$, and the latter precede the elements of $\wh \Ec$, $c$ included in the ordering in an
arbitrary way, then the ordered monomials in the set of elements $k_i^{(r)},e_{ij}^{\tss(r)},f_{ji}^{\tss(r)}$ and $c$ with $r\in \mathbb{Z}^{\times}$
span $\wh \DX_\hbar(\g_{N})$. It follows that
\eqref{surjhom} is injective.

\section{Isomorphism theorem for the Drinfeld Yangian double}
\label{sec:isom}

In this section, we will prove the Main Theorem as stated in the Introduction. It follows from the definition of the series $H_i^{\pm}(u),E_i(u),F_i(u)$ that all the coefficients $h^{}_{i\tss r},~\xi_{i\tss r}^{+},~\xi_{i\tss r}^{-}$ $(i=1,\dots,n,~r\in\mathbb{Z})$ of $H_i^{\pm}(u),E_i(u),F_i(u)$ belong to the subalgebra $\DY_\hbar(\g_N)$.

\bpr\label{prop:YgN}
The subalgebra $\DY_\hbar(\g_N)$ of $\DX_\hbar(\g_N)$ is generated by
the coefficients $h^{}_{i\tss r},~\xi_{i\tss r}^{+},~\xi_{i\tss r}^{-},~c$ with
$i=1,\dots,n$ and $r\in\mathbb{Z}$.
\epr

\bpf
Let $\wt \DY_\hbar(\g_N)$ be the subalgebra of $\DX_\hbar(\g_N)$ generated by
the coefficients $h^{}_{i\tss r},~\xi_{i\tss r}^{+},~\xi_{i\tss r}^{-},~c$ with
$i=1,\dots,n$ and $r\in\mathbb{Z}$. By the definition of the series $H^{\pm}_i(u)$, we can express $k^{\pm}_{1}(u)k^{\pm}_{n+1}(u)^{-1}$ as a product of the
series of the form $H^{\pm}_i(u)^{-1}$ with some shifts of $u$ by constants. As a next step, Theorem~\ref{thm:Center} implies that
$k^{\pm}_{1}(u+\kappa\hbar)k^{\pm}_{n+1}(u)$ equals $z^{\pm}_N(u)$ times the product of the shifted series $H^{\pm}_i(u)^{-1}$. It follows that all coefficients of $k^{\pm}_1(u)$ belong to the subalgebra $\wt \DY_\hbar(\g_N)$. Hence, all coefficients of the series $k^{\pm}_i(u)$ with $i=1,\dots,n+1$ belong to the subalgebra $\wt \DY_\hbar(\g_N)$. Moreover, for any $i$, the elements $e_i^{(r)}$ and $f_i^{(r)}$ can be expressed as linear combinations of the $\xi_{i\tss s}^{+}$ and $\xi_{i\tss s}^{-}$, respectively. This means
the elements $e_i^{(r)}$ and $f_i^{(r)}$ belong to the subalgebra $\wt \DY_\hbar(\g_N)$. Theorem~\ref{thm:dp} says that the coefficients of the series
$k^{\pm}_i(u)$ with $i=1,\dots,n+1$, and $e^{\pm}_i(u),f^{\pm}_i(u)$ with $i=1,\dots,n$ and $c$ generate the algebra $\DX_\hbar(\g_N)$. So the elements of $\wt \DY_\hbar(\g_N)$ together with the coefficients $z^{(r)}_N$ of the series $z^{\pm}_N(u)$ generate the algebra $\DX_\hbar(\g_N)$. By the tensor decomposition \eqref{tensordecom} of the algebra $\DX_\hbar(\g_N)$, we conclude that $\wt \DY_\hbar(\g_N)=\DY_\hbar(\g_N)$.
\epf

We are in the position to prove the generators
$h^{}_{i\tss r}$, $\xi_{i\tss r}^{+}$ and
$\xi_{i\tss r}^{-}$ of the subalgebra $\DY_\hbar(\g_N)$ satisfy the
same relations as the Drinfeld Yangian double $\DY_\hbar^{D}(\g_N)$.

\bthm\label{thm:relinYgN}
The Yangian double $\DY_\hbar(\mathfrak{g}_N)$ is topologically generated by
the coefficients of the series
$H_i^{\pm}(u)$ with $i=1,\dots,n$, and
$E_i(u)$ and $F_i(u)$ with $i=1,\dots, n$ and $c$
subject only to the following sets of relations, where the indices
take all admissible values unless specified otherwise.
\begin{align}
[H_{i}^{\pm}(u),H_{j}^{\pm}(v)]=0,
\end{align}
\begin{align}\label{BCDhihj}
(u_{\pm}-v_{\mp}-B_{ij}\hbar)(u_{\mp}-v_{\pm}+B_{ij}\hbar)H_{i}^{\pm}(u)H_{j}^{\mp}(v)\\
=(u_{\pm}-v_{\mp}+B_{ij}\hbar)(u_{\mp}-v_{\pm}-B_{ij}\hbar)H_{j}^{\mp}(v)H_{i}^{\pm}(u),\nonumber
\end{align}
\begin{align}\label{BCDhiej}
H_{i}^{\pm}(u)^{-1}E_{j}(v)H_{i}^{\pm}(u)=\frac{u_{\mp}-v+B_{ij}\hbar}{u_{\mp}-v-B_{ij}\hbar}E_{j}(v),\\
\label{BCDhifj}
H_{i}^{\pm}(u)F_{j}(v)H_{i}^{\pm}(u)^{-1}=\frac{u_{\pm}-v+B_{ij}\hbar}{u_{\pm}-v-B_{ij}\hbar}F_{j}(v),
\end{align}
\begin{align}\label{BCDeiej}
(u-v+B_{ij}\hbar)E_{i}(u)E_{j}(v)=(u-v-B_{ij}\hbar)E_{j}(v)E_{i}(u),\\
\label{BCDfifj}
(u-v-B_{ij}\hbar)F_{i}(u)F_{j}(v)=(u-v+B_{ij}\hbar)F_{j}(v)F_{i}(u),
\end{align}
\begin{align}\label{Serreeiej}
\sum_{\sigma\in\mathfrak{S}_{m}}[E_{i}(u_{\sigma(1)}),[E_{i}(u_{\sigma(2)})\cdots,[E_{i}(u_{\sigma(m)}),E_{j}(v)]\cdots]=0, \\
\label{Serrefifj}
\sum_{\sigma\in\mathfrak{S}_{m}}[F_{i}(u_{\sigma(1)}),[F_{i}(u_{\sigma(2)})\cdots,[F_{i}(u_{\sigma(m)}),F_{j}(v)]\cdots]=0, \\
~~ i\neq j,m=1-a_{ij},\nonumber
\end{align}
\begin{align}\label{BCDeifj}
[E_{i}(u),F_{j}(v)]=\frac{1}{\hbar}\delta_{ij}\{\delta(u_{-}-v_{+})H_{i}^{-}(v_{+})-\delta(u_{+}-v_{-})H_{i}^{+}(u_{+})\}.
\end{align}
Here we set $B_{ij}=\frac{1}{2}(\al_i,\al_j)$.
\ethm

\bpf
It suffices to write the relations of Theorem~\ref{thm:dp} in terms of the series $H_{i}^{\pm}(u),E_{i}(u)$ and $F_{i}(u)$. All relations are immediate from the definition of the series $H_{i}^{\pm}(u),E_{i}(u)$ and $F_{i}(u)$. In the following we will demonstrate several relations
for various types to show the idea. %$B$, $C$ and $D$.
\paragraph{Type $B_n$.}
To verify the $i=j=n$ case of \eqref{BCDhiej}, we get $H^{\pm}_{n}(u)=k^{\pm}_n\big(u-(n-1)\hbar/2\big)^{-1}\ts k^{\pm}_{n+1}\big(u-(n-1)\hbar/2\big)$ and $E_{n}(u)=X_n^+\big(u-(n-1)\hbar/2\big)$ from definition. By Theorem~\ref{thm:dp}, we have
\begin{align*}
k_{n}^{\pm}(u)^{-1}X_{n}^{+}(v)k_{n}^{\pm}(u)&=\frac{u_{\mp}-v-\hbar}{u_{\mp}-v}X_{n}^{+}(v),\\
k_{n+1}^{\pm}(u)^{-1}X_{n}^{+}(v)k_{n+1}^{\pm}(u)&=\frac{(u_{\mp}-v-\hbar)(u_{\mp}-v+\frac{1}{2}\hbar)}{(u_{\mp}-v)(u_{\mp}-v-\frac{1}{2}\hbar)}X_{n}^{+}(v).
\end{align*}
It follows that $H_{n}^{\pm}(u)^{-1}E_{n}(v)H_{n}^{\pm}(u)=\frac{u_{\mp}-v+B_{nn}\hbar}{u_{\mp}-v-B_{nn}\hbar}E_{n}(v)$, where $B_{nn}=\frac{1}{2}$.
To verify the $i=j=n$ case of \eqref{BCDeifj}, we get $E_{n}(u)=X_n^+\big(u-(n-1)\hbar/2\big)$ and $F_{n}(u)=X_n^-\big(u-(n-1)\hbar/2\big)$ from definition. By Theorem~\ref{thm:dp}, we have
\begin{align*}
[X_{n}^{+}(u),X_{n}^{-}(v)]=\hbar\{\delta(u_{-}-v_{+})k_{n+1}^{-}(v_{+})k_{n}^{-}(v_{+})^{-1}
-\delta(u_{+}-v_{-})k_{n+1}^{+}(u_{+})k_{n}^{+}(u_{+})^{-1}\}.
\end{align*}
Since $H^{\pm}_{n}(u)=k^{\pm}_n\big(u-(n-1)\hbar/2\big)^{-1}\ts k^{\pm}_{n+1}\big(u-(n-1)\hbar/2\big)$, This yields $[E_{n}(u),F_{n}(v)]=\frac{1}{\hbar}\{\delta(u_{-}-v_{+})H_{n}^{-}(v_{+})-\delta(u_{+}-v_{-})H_{n}^{+}(u_{+})\}$.
By Proposition~\ref{pro:Arelations}, we have
\begin{align*}
\frac{u-v-\frac{\hbar}{2}-\frac{\hbar C}{4}}{u-v-\frac{\hbar C}{4}}K^{+}_{1}(2u+\hbar)&K^{-}_{2}(2v+\hbar)K^{-}_{1}(2v+\hbar)^{-1}\\
&=\frac{u-v-\frac{\hbar}{2}+\frac{\hbar C}{4}}{u-v+\frac{\hbar C}{4}}K^{-}_{2}(2v+\hbar)K^{-}_{1}(2v+\hbar)^{-1}K^{+}_{1}(2u+\hbar),\\
\frac{u-v-\hbar-\frac{\hbar C}{4}}{u-v-\frac{\hbar}{2}-\frac{\hbar C}{4}}K^{+}_{1}(2u)&K^{-}_{2}(2v+\hbar)K^{-}_{1}(2v+\hbar)^{-1}\\
&=\frac{u-v-\hbar+\frac{\hbar C}{4}}{u-v-\frac{\hbar}{2}+\frac{\hbar C}{4}}K^{-}_{2}(2v+\hbar)K^{-}_{1}(2v+\hbar)^{-1}K^{+}_{1}(2u).
\end{align*}
in $\DY_\hbar(\gl_2)$.
Apply Lemma~\ref{lem:lowrank B}, we get $$\frac{u_{-}-v_{+}-\hbar}{u_{-}-v_{+}}k^{+}_{n}(u)k^{-}_{n+1}(v)k^{-}_{n}(v)^{-1}=\frac{u_{+}-v_{-}-\hbar}{u_{+}-v_{-}}k^{-}_{n+1}(v)k^{-}_{n}(v)^{-1}k^{+}_{n}(u).$$ Corollary~\ref{cor:commu} implies $k^{+}_{n-1}(u)^{-1}k^{-}_{n+1}(v)=k^{-}_{n+1}(v)k^{+}_{n-1}(u)^{-1}$. Since $k^{+}_{n-1}(u)^{-1}k^{-}_{n}(v)^{-1}=k^{-}_{n}(v)^{-1}k^{+}_{n-1}(u)^{-1}$, then we have $k^{+}_{n-1}(u)^{-1}H^{-}_{n}(v)=H^{-}_{n}(v)k^{+}_{n-1}(u)^{-1}$.
It follows that $$(u_{-}-v_{+}-\frac{\hbar}{2})(u_{+}-v_{-}+\frac{\hbar}{2})H^{+}_{n-1}(u)H^{-}_{n}(v)=(u_{-}-v_{+}+\frac{\hbar}{2})(u_{+}-v_{-}-\frac{\hbar}{2})H^{-}_{n}(v)H^{+}_{n-1}(u),$$ which is just \eqref{BCDhihj} when $i=n-1,j=n$. All other remaining cases follow by similar
calculations.

\paragraph{Type $C_n$.}
To verify the $i=j=n$ case of \eqref{BCDhiej}, we get $H^{\pm}_{n}(u)=2\tss k^{\pm}_n\big(u-n\hbar/2\big)^{-1}\ts k^{\pm}_{n+1}\big(u-n\hbar/2\big)$ and $E_n(u)=X_n^+\big(u-n\hbar/2\big)$ from definition. By Theorem~\ref{thm:dp}, we have
\begin{align*}
k_{n}^{\pm}(u)^{-1}X_{n}^{+}(v)k_{n}^{\pm}(u)&=\frac{u_{\mp}-v-2\hbar}{u_{\mp}-v}X_{n}^{+}(v)\\
k_{n+1}^{\pm}(u)^{-1}X_{n}^{+}(v)k_{n+1}^{\pm}(u)&=\frac{u_{\mp}-v+2\hbar}{u_{\mp}-v}X_{n}^{+}(v)
\end{align*}
It follows that $H_{n}^{\pm}(u)^{-1}E_{n}(v)H_{n}^{\pm}(u)=\frac{u_{\mp}-v+B_{nn}\hbar}{u_{\mp}-v-B_{nn}\hbar}E_{n}(v)$, where $B_{nn}=2$.
To verify the $i=j=n$ case of \eqref{BCDeifj}, we get $E_{n}(u)=X_n^+\big(u-n\hbar/2\big)$ and $F_{n}(u)=X_n^-\big(u-n\hbar/2\big)$ from definition. By Theorem~\ref{thm:dp}, we have
\begin{align*}
[X_{n}^{+}(u),X_{n}^{-}(v)]=2\hbar\{\delta(u_{-}-v_{+})k_{n+1}^{-}(v_{+})k_{n}^{-}(v_{+})^{-1}
-\delta(u_{+}-v_{-})k_{n+1}^{+}(u_{+})k_{n}^{+}(u_{+})^{-1}\}.
\end{align*}
This yields $[E_{n}(u),F_{n}(v)]=\frac{1}{\hbar}\{\delta(u_{-}-v_{+})H_{n}^{-}(v_{+})-\delta(u_{+}-v_{-})H_{n}^{+}(u_{+})\}$.
Apply Lemma~\ref{lem:lowrank C} and Proposition~\ref{pro:Arelations}, we have
\begin{align*}
\frac{u_{-}-v_{+}-\hbar}{u_{-}-v_{+}+\hbar}k^{+}_{n}(u-\frac{n-2}{2}\hbar)H^{-}_{n}(v)=\frac{u_{+}-v_{-}-\hbar}{u_{+}-v_{-}+\hbar}H^{-}_{n}(v)k^{+}_{n}(u-\frac{n-2}{2}\hbar).
\end{align*}
Corollary~\ref{cor:commu} implies $k^{+}_{n-1}(u)^{-1}k^{-}_{n+1}(v)=k^{-}_{n+1}(v)k^{+}_{n-1}(u)^{-1}$. Since $k^{+}_{n-1}(u)^{-1}k^{-}_{n}(v)^{-1}=k^{-}_{n}(v)^{-1}k^{+}_{n-1}(u)^{-1}$, then we have $k^{+}_{n-1}(u)^{-1}H^{-}_{n}(v)=H^{-}_{n}(v)k^{+}_{n-1}(u)^{-1}$.
It follows that $$(u_{-}-v_{+}-\hbar)(u_{+}-v_{-}+\hbar)H^{+}_{n-1}(u)H^{-}_{n}(v)=(u_{-}-v_{+}+\hbar)(u_{+}-v_{-}-\hbar)H^{-}_{n}(v)H^{+}_{n-1}(u),$$ which is just \eqref{BCDhihj} when $i=n-1,j=n$. All other remaining cases follow by similar calculations.

\paragraph{Type $D_n$.}
To verify the $i=n-1,j=n$ case of \eqref{BCDhiej}, we get $H^{\pm}_{n-1}(u)=k^{\pm}_{n-1}\big(u-(n-2)\hbar/2\big)^{-1}\ts k^{\pm}_{n}\big(u-(n-2)\hbar/2\big)
$ and $E_n(u)=X_n^+\big(u-(n-2)\hbar/2\big)$ from definition.
By Theorem~\ref{thm:dp}, we have
\begin{align*}
k_{n}^{\pm}(u)^{-1}X_{n}^{+}(v)k_{n}^{\pm}(u)&=\frac{u_{\mp}-v-\hbar}{u_{\mp}-v}X_{n}^{+}(v),\\
k_{n-1}^{\pm}(u)^{-1}X_{n}^{+}(v)k_{n-1}^{\pm}(u)&=\frac{u_{\mp}-v-\hbar}{u_{\mp}-v}X_{n}^{+}(v),
\end{align*}
It follows that $H_{n-1}^{\pm}(u)^{-1}E_{n}(v)H_{n-1}^{\pm}(u)=E_{n}(v)$, which is just \eqref{BCDhiej} when $i=n-1,j=n$.
To verify the $i=j=n$ case of \eqref{BCDeifj}, we get $E_{n}(u)=X_n^+\big(u-(n-2)\hbar/2\big)$ and $F_{n}(u)=X_n^-\big(u-(n-2)\hbar/2\big)$ from definition. By Theorem~\ref{thm:dp}, we have
\begin{align*}
[X_{n}^{+}(u),X_{n}^{-}(v)]=\hbar\{\delta(u_{-}-v_{+})k_{n+1}^{-}(v_{+})k_{n-1}^{-}(v_{+})^{-1}
-\delta(u_{+}-v_{-})k_{n+1}^{+}(u_{+})k_{n-1}^{+}(u_{+})^{-1}\}.
\end{align*}
Since $H^{\pm}_n(u)=k^{\pm}_{n-1}\big(u-(n-2)\hbar/2\big)^{-1}\ts k^{\pm}_{n+1}\big(u-(n-2)\hbar/2\big)$, this yields $[E_{n}(u),F_{n}(v)]=\frac{1}{\hbar}\{\delta(u_{-}-v_{+})H_{n}^{-}(v_{+})-\delta(u_{+}-v_{-})H_{n}^{+}(u_{+})\}$.
Apply Lemma~\ref{lem:lowrank D} and Proposition~\ref{pro:Arelations}, we have
\begin{align*}
\frac{u_{-}-v_{+}-\hbar}{u_{-}-v_{+}}k^{+}_{n-1}(u)k^{-}_{n+1}(v)k^{-}_{n-1}(v)^{-1}&=\frac{u_{+}-v_{-}-\hbar}{u_{+}-v_{-}}k^{-}_{n+1}(v)k^{-}_{n-1}(v)^{-1}k^{+}_{n-1}(u),\\
\frac{u_{-}-v_{+}-\hbar}{u_{-}-v_{+}}k^{+}_{n}(u)k^{-}_{n+1}(v)k^{-}_{n-1}(v)^{-1}&=\frac{u_{+}-v_{-}-\hbar}{u_{+}-v_{-}}k^{-}_{n+1}(v)k^{-}_{n-1}(v)^{-1}k^{+}_{n}(u).
\end{align*}
Since $H^{+}_{n-1}(u)=k^{+}_{n-1}\big(u-(n-2)\hbar/2\big)^{-1}\ts k^{+}_{n}\big(u-(n-2)\hbar/2\big)$, then we get $H^{+}_{n-1}(u)H^{-}_{n}(v)=H^{-}_{n}(v)H^{+}_{n-1}(u)$,
which is just \eqref{BCDhihj} when $i=n-1,j=n$. All other remaining cases follow by similar calculations.

It follows from Propositions~\ref{prop:YgN} that the mapping
$\DY_\hbar^D(\g_{N})\to\DY_\hbar(\g_N)$ considered in the Main Theorem
is a surjective homomorphism. By the decomposition
\beql{decom}
\DY_\hbar(\g_N)= \Ec\widetilde{\otimes} (\DY_\hbar(\g_N)\cap\Hc)\widetilde{\otimes} \Fc\widetilde{\otimes}\Ac c,
\eeq
the corresponding arguments of the proof of Theorem \ref{thm:dp} gives that the mapping is also injective. This completes the proof of the Main Theorem.
\epf

\section{Bosonization of level 1 modules}
\label{sec:boson}

Here we construct level $1$ $\DY_\hbar(\g)$-module for $\g_N=\oa_{2n+1},\spa_{2n},\oa_{2n}$ in terms of bosons. Our ideal comes from Iohara's work \cite{io:br}. Introduce bosons $\{a_{i,k} \mid 1\leq i\leq n+4, k\in \mathbb{Z} \backslash \{0\}\}$ satisfying:
$$[a_{i,k},a_{j,l}]=k\delta_{ij}\delta_{k+l,0}.$$
Define the standard bilinear form $(\epsilon_{i},\epsilon_{j})=\delta_{ij}$. For $\g_N=\oa_{2n+1}$, set $\alpha_{j}=\epsilon_{j}-\epsilon_{j+1}$ for $j=1,\cdots,n-1$,$\alpha_{n}=\epsilon_{n},\beta_{1}=\epsilon_{n+1},\beta_{i}=\sqrt{2}\epsilon_{n+i}$ for $i=2,3,4$. We introduce the Fock space $\mathcal {F}=\Ac[a_{j,-k}(1\leq i\leq n+3,k\in \mathbb{Z}\backslash \{0\})]\widetilde{\bigotimes}\Ac[Q]$,
where $Q=\bigoplus_{i=1}^{n}\mathbb{Z}\alpha_{i}\oplus\bigoplus_{j=1}^{4}\mathbb{Z}\beta_{j}$, $\Ac[Q]$ is the group algebra of $Q$ over $\Ac$. Note that $\mathcal{F}$ is a topological tensor product. On this space, we define the action of the operators $a_{i,k},\partial\epsilon_{j},e^{\epsilon_{j}}(1\leq i\leq n+3,1\leq j\leq n+4)$ by
\begin{equation}
a_{i,k}\cdot f\otimes e^{\beta}=
\begin{cases}
a_{i,k}f\otimes e^{\beta}&\text{if $k<0$}\\
[a_{i,k},f]\otimes e^{\beta}&\text{if $k>0$}
\end{cases}
\end{equation}
\begin{align*}
\partial\epsilon_{j}\cdot f\otimes e^{\beta}=&(\epsilon_{j},\beta)f\otimes e^{\beta},\\
e^{\epsilon_{j}}\cdot f\otimes e^{\beta}=&f\otimes e^{\epsilon_{j}+\beta},
\end{align*}
for $f\otimes e^{\beta}\in \mathcal {F}$.

\bthm
The following assignment defines a $\DY_\hbar(\oa_{2n+1})$-module structure on $\mathcal {F}$.\\
For $1\leq j \leq n-2$,
\begin{align*}
E_{j}(u)\mapsto
&exp[\sum_{k>0}\frac{a_{j,-k}}{k}\{(u-\frac{1}{4}\hbar)^{k}+(u+\frac{3}{4}\hbar)^{k}\}-\sum_{k>0}\frac{a_{j+1,-k}+a_{j-1,-k}}{k}(u+\frac{1}{4}\hbar)^{k}]\\
&exp[-\sum_{k>0}\frac{a_{j,k}}{k}(u-\frac{1}{4}\hbar)^{-k}]e^{\alpha_{j}}[(-1)^{j-1}(u-\frac{1}{4}\hbar)]^{\partial\alpha_{j}},
\end{align*}
\begin{align*}
F_{j}(u)\mapsto
&exp[-\sum_{k>0}\frac{a_{j,-k}}{k}\{(u+\frac{1}{4}\hbar)^{k}+(u-\frac{3}{4}\hbar)^{k}\}+\sum_{k>0}\frac{a_{j+1,-k}+a_{j-1,-k}}{k}(u-\frac{1}{4}\hbar)^{k}]\\
&exp[\sum_{k>0}\frac{a_{j,k}}{k}(u+\frac{1}{4}\hbar)^{-k}]e^{-\alpha_{j}}[(-1)^{j-1}(u+\frac{1}{4}\hbar)]^{-\partial\alpha_{j}},
\end{align*}
\begin{align*}
H_{j}^{+}(u)\mapsto
exp[-\sum_{k>0}\frac{a_{j,k}}{k}\{(u-\frac{1}{2}\hbar)^{-k}-(u+\frac{1}{2}\hbar)^{-k}\}](\frac{u+\frac{1}{2}\hbar}{u-\frac{1}{2}\hbar})^{-\partial\alpha_{j}},
\end{align*}
\begin{align*}
H_{j}^{-}(u)\mapsto
exp[-\sum_{k>0}\frac{a_{j,-k}}{k}&\{(u-\hbar)^{k}-(u+\hbar)^{k}\}\\
&+\sum_{k>0}\frac{a_{j+1,-k}+a_{j-1,-k}}{k}\{(u-\frac{1}{2}\hbar)^{k}-(u+\frac{1}{2}\hbar)^{k}\}].
\end{align*}

\begin{align*}
E_{n-1}(u)\mapsto
&exp[\sum_{k>0}\frac{a_{n-1,-k}}{k}\{(u-\frac{1}{4}\hbar)^{k}+(u+\frac{3}{4}\hbar)^{k}\}-\sum_{k>0}\frac{a_{n-2,-k}+a_{n,-k}}{k}(u+\frac{1}{4}\hbar)^{k}]\\
&exp[-\sum_{k>0}\frac{a_{n-1,k}}{k}(u-\frac{1}{4}\hbar)^{-k}]e^{\alpha_{n-1}}[(-1)^{n-2}(u-\frac{1}{4}\hbar)]^{\partial\alpha_{n-1}+\partial\beta_{1}},
\end{align*}
\begin{align*}
F_{n-1}(u)\mapsto
&exp[-\sum_{k>0}\frac{a_{n-1,-k}}{k}\{(u+\frac{1}{4}\hbar)^{k}+(u-\frac{3}{4}\hbar)^{k}\}+\sum_{k>0}\frac{a_{n-2,-k}+a_{n,-k}}{k}(u-\frac{1}{4}\hbar)^{k}]\\
&exp[\sum_{k>0}\frac{a_{n-1,k}}{k}(u+\frac{1}{4}\hbar)^{-k}]e^{-\alpha_{n-1}}[(-1)^{n-2}(u+\frac{1}{4}\hbar)]^{-\partial\alpha_{n-1}-\partial\beta_{1}},
\end{align*}
\begin{align*}
E_{n}(u)\mapsto
&exp[-\sum_{k>0}\frac{a_{n-1,-k}}{k}(u+\frac{1}{4}\hbar)^{k}+\sum_{k>0}\frac{a_{n+1,-k}}{k}\{(u-\frac{1}{4}\hbar)^{k}+(u+\frac{3}{4}\hbar)^{k}\}\\
&+\sum_{k>0}\frac{a_{n+2,-k}}{k}\{(u-\frac{3}{4}\hbar)^{k}+(u+\frac{3}{4}\hbar)^{k}\}]exp[-\sum_{k>0}\frac{a_{n,k}+a_{n+1,k}}{k}(u-\frac{1}{4}\hbar)^{-k}\\
&-\sum_{k>0}\frac{a_{n+2,k}}{k}(u+\frac{1}{4}\hbar)^{-k}]e^{-\beta_{1}}[(-1)^{n-1}(u-\frac{1}{4}\hbar)]^{\partial\alpha_{n}}e^{\beta_{2}+\beta_{3}}(u-\frac{1}{4}\hbar)^{\partial\beta_{2}}(u+\frac{1}{4}\hbar)^{\partial\beta_{3}},
\end{align*}
\begin{align*}
F_{n}(u)\mapsto
&exp[\sum_{k>0}\frac{a_{n-1,-k}}{k}(u-\frac{1}{4}\hbar)^{k}-\sum_{k>0}\frac{a_{n+1,-k}}{k}\{(u+\frac{1}{4}\hbar)^{k}+(u-\frac{3}{4}\hbar)^{k}\}\\
&-\sum_{k>0}\frac{a_{n+3,-k}}{k}\{(u-\frac{3}{4}\hbar)^{k}+(u+\frac{3}{4}\hbar)^{k}\}]exp[\sum_{k>0}\frac{a_{n,k}+a_{n+1,k}}{k}(u+\frac{1}{4}\hbar)^{-k}\\
&+\sum_{k>0}\frac{a_{n+3,k}}{k}(u-\frac{1}{4}\hbar)^{-k}]e^{\beta_{1}}[(-1)^{n-1}(u+\frac{1}{4}\hbar)]^{-\partial\alpha_{n}}e^{-\beta_{2}+\beta_{4}}(u+\frac{1}{4}\hbar)^{-\partial\beta_{2}}(u-\frac{1}{4}\hbar)^{\partial\beta_{4}},
\end{align*}
\begin{align*}
H_{n-1}^{+}(u)\mapsto
exp[-\sum_{k>0}\frac{a_{n-1,k}}{k}\{(u-\frac{1}{2}\hbar)^{-k}-(u+\frac{1}{2}\hbar)^{-k}\}](\frac{u+\frac{1}{2}\hbar}{u-\frac{1}{2}\hbar})^{-\partial\alpha_{n-1}-\partial\beta_{1}},
\end{align*}
\begin{align*}
H_{n-1}^{-}(u)\mapsto
exp[-\sum_{k>0}\frac{a_{n-1,-k}}{k}&\{(u-\hbar)^{k}-(u+\hbar)^{k}\}\\
&+\sum_{k>0}\frac{a_{n-2,-k}+a_{n,-k}}{k}\{(u-\frac{1}{2}\hbar)^{k}-(u+\frac{1}{2}\hbar)^{k}\}],
\end{align*}
\begin{align*}
H_{n}^{+}(u)\mapsto
&exp[\sum_{k>0}\frac{a_{n+2,-k}}{k}\{(u-\hbar)^{k}+(u+\frac{1}{2}\hbar)^{k}\}-\sum_{k>0}\frac{a_{n+3,-k}}{k}\{(u+\hbar)^{k}+(u-\frac{1}{2}\hbar)^{k}\}]\\
&exp[\sum_{k>0}\frac{a_{n,k}+a_{n+1,k}}{k}\{(u+\frac{1}{2}\hbar)^{-k}-(u-\frac{1}{2}\hbar)^{-k}\}+\sum_{k>0}\frac{a_{n+3,k}-a_{n+2,k}}{k}u^{-k}]\\
&(\frac{u+\frac{1}{2}\hbar}{u-\frac{1}{2}\hbar})^{-\partial\alpha_{n}-\partial\beta_{2}}e^{\beta_{3}+\beta_{4}}u^{\partial\beta_{3}+\partial\beta_{4}},
\end{align*}
\begin{align*}
H_{n}^{-}(u)\mapsto
&exp[\sum_{k>0}\frac{a_{n-1,-k}}{k}\{(u-\frac{1}{2}\hbar)^{k}-(u+\frac{1}{2}\hbar)^{k}\}-\sum_{k>0}\frac{a_{n+1,-k}}{k}\{(u-\hbar)^{k}-(u+\hbar)^{k}\}\\
&+\sum_{k>0}\frac{a_{n+2,-k}}{k}\{(u-\frac{1}{2}\hbar)^{k}+(u+\hbar)^{k}\}-\sum_{k>0}\frac{a_{n+3,-k}}{k}\{(u+\frac{1}{2}\hbar)^{k}+(u-\hbar)^{k}\}]\\
&exp[-\sum_{k>0}\frac{a_{n+2,k}}{k}(u+\frac{1}{2}\hbar)^{-k}+\sum_{k>0}\frac{a_{n+3,k}}{k}(u-\frac{1}{2}\hbar)^{-k}]e^{\beta_{3}+\beta_{4}}(u+\frac{1}{2}\hbar)^{\partial\beta_{3}}(u-\frac{1}{2}\hbar)^{\partial\beta_{4}}.
\end{align*}

\ethm

\bpf
We need to verify the relations in Theorem~\ref{thm:relinYgN} with $c=1$. This is done by some routine calculations. For instance, to verify \eqref{BCDhiej} with $i=j=n$,
we define the normal ordering $:\cdot:$ of the fields by regarding $a_{j,k}(k<0),e^{\epsilon_{j}}$ as creation operators and $a_{j,k}(k>0),\partial_{\epsilon_{j}}$
as annihilation operators. The definition of $\partial\alpha$ and $e^{\beta}$ implies that $[\partial\alpha,e^{\beta}]=(\alpha,\beta)e^{\beta}$. Then we have
\begin{align*}
(\frac{u+\frac{1}{2}\hbar}{u-\frac{1}{2}\hbar})^{-\partial\beta_{2}}e^{\beta_{2}}&=e^{\beta_{2}}(\frac{u+\frac{1}{2}\hbar}{u-\frac{1}{2}\hbar})^{-\partial\beta_{2}}(\frac{u+\frac{1}{2}\hbar}{u-\frac{1}{2}\hbar})^{(-\beta_{2},\beta_{2})}\\
&=e^{\beta_{2}}(\frac{u+\frac{1}{2}\hbar}{u-\frac{1}{2}\hbar})^{-\partial\beta_{2}}(\frac{u+\frac{1}{2}\hbar}{u-\frac{1}{2}\hbar})^{-2}
\end{align*}
\begin{align*}
u^{\partial\beta_{3}}e^{\beta_{3}}=e^{\beta_{3}}u^{\partial\beta_{3}}u^{(\beta_{3},\beta_{3})}=e^{\beta_{3}}u^{\partial\beta_{3}}u^{2}
\end{align*}
Note that $exp(A)exp(B)=exp(B)exp(A)exp([A,B])$ if $[A,B]$ commutes with $A$ and $B$. Since $[a_{i,k},a_{i,-k}]=k$, using the expansion $\ln(1-u)=-\sum_{k>0}\frac{u^{k}}{k}$, we have the following relations
\begin{align*}
&exp[\sum_{k>0}\frac{a_{n+1,k}}{k}\{(u+\frac{1}{2}\hbar)^{-k}-(u-\frac{1}{2}\hbar)^{-k}\}]exp[\sum_{k>0}\frac{a_{n+1,-k}}{k}\{(v-\frac{1}{4}\hbar)^{k}+(v+\frac{3}{4}\hbar)^{k}\}]\\
=&exp[\sum_{k>0}\frac{a_{n+1,-k}}{k}\{(v-\frac{1}{4}\hbar)^{k}+(v+\frac{3}{4}\hbar)^{k}\}]exp[\sum_{k>0}\frac{a_{n+1,k}}{k}\{(u+\frac{1}{2}\hbar)^{-k}-(u-\frac{1}{2}\hbar)^{-k}\}]\\
&\cdot\frac{(u+\frac{1}{2}\hbar)^{2}(u-v-\frac{5}{4}\hbar)}{(u-\frac{1}{2}\hbar)^{2}(u-v+\frac{3}{4}\hbar)}
\end{align*}
\begin{align*}
&exp[-\sum_{k>0}\frac{a_{n+2,k}}{k}u^{-k}]exp[\sum_{k>0}\frac{a_{n+2,-k}}{k}\{(v-\frac{3}{4}\hbar)^{k}+(v+\frac{3}{4}\hbar)^{k}\}]\\
=&exp[\sum_{k>0}\frac{a_{n+2,-k}}{k}\{(v-\frac{3}{4}\hbar)^{k}+(v+\frac{3}{4}\hbar)^{k}\}]exp[-\sum_{k>0}\frac{a_{n+2,k}}{k}u^{-k}]\frac{(u-v+\frac{3}{4}\hbar)(u-v-\frac{3}{4}\hbar)}{u^{2}}
\end{align*}
So we obtain the following operator product expansion (OPE):
\begin{align*}
H^{+}_{n}(u)E_{n}(v)=:H^{+}_{n}(u)E_{n}(v):(u-v-\frac{5}{4}\hbar)(u-v-\frac{3}{4}\hbar).
\end{align*}
Similarly, we have $E_{n}(v)H^{+}_{n}(u)=:E_{n}(v)H^{+}_{n}(u):(u-v-\frac{5}{4}\hbar)(u-v+\frac{1}{4}\hbar)$. This implies $H^{+}_{n}(u)^{-1}E_{n}(v)H^{+}_{n}(u)=\frac{u-v+\frac{1}{4}\hbar}{u-v-\frac{3}{4}\hbar}E_{n}(v)$. To verify \eqref{BCDeiej} with $i=n-1,j=n$, we have
\begin{align*}
&exp[-\sum_{k>0}\frac{a_{n-1,k}}{k}(u-\frac{1}{4}\hbar)^{-k}]exp[-\sum_{k>0}\frac{a_{n-1,-k}}{k}(v+\frac{1}{4}\hbar)^{k}]\\
=&exp[-\sum_{k>0}\frac{a_{n-1,-k}}{k}(v+\frac{1}{4}\hbar)^{k}]exp[-\sum_{k>0}\frac{a_{n-1,k}}{k}(u-\frac{1}{4}\hbar)^{-k}]\frac{u-\frac{1}{4}\hbar}{u-v-\frac{1}{2}\hbar}
\end{align*}
Also, we have
\begin{align*}
[(-1)^{n-2}(u-\frac{1}{4}\hbar)]^{\partial\beta_{1}}e^{-\beta_{1}}=e^{-\beta_{1}}[(-1)^{n-2}(u-\frac{1}{4}\hbar)]^{\partial\beta_{1}}[(-1)^{n-2}(u-\frac{1}{4}\hbar)]^{-1}
\end{align*}
So we obtain the following operator product expansion (OPE):
\begin{align}\label{e:ope}
E_{n-1}(u)E_{n}(v)=:E_{n-1}(u)E_{n}(v):\frac{(-1)^{-n+2}}{u-v-\frac{1}{2}\hbar}.
\end{align}
Similarly, we have $E_{n}(v)E_{n-1}(u)=:E_{n}(v)E_{n-1}(u):\frac{(-1)^{-n+1}}{v-u-\frac{1}{2}\hbar}$. This implies $(u-v-\frac{1}{2}\hbar)E_{n-1}(u)E_{n}(v)=(u-v+\frac{1}{2}\hbar)E_{n}(v)E_{n-1}(u)$.
Note that the Serre relation \eqref{Serreeiej} can be essentially deduced from \eqref{BCDeiej} or more precisely the OPEs like \eqref{e:ope}.
For instance, to verify \eqref{Serreeiej} with $i=n-1,j=n$, we have
\begin{align*}
(u-v)[E_{n-1}(u),E_{n}(v)]=\frac{\hbar}{2}(E_{n-1}(u)E_{n}(v)+E_{n}(v)E_{n-1}(u)),\\
(u-v)[E_{n-1}(u),E_{n-1}(v)]=-\hbar(E_{n-1}(u)E_{n-1}(v)+E_{n-1}(v)E_{n-1}(u)),
\end{align*}
by \eqref{BCDeiej}. It follows that
\begin{align*}
&(u-v)(u-w)(v-w)([E_{n-1}(u),[E_{n-1}(v),E_{n}(w)]]+[E_{n-1}(v),[E_{n-1}(u),E_{n}(w)]])\\
=&-\frac{\hbar^2}{4}(u-v)([E_{n-1}(u),[E_{n-1}(v),E_{n}(w)]]+[E_{n-1}(v),[E_{n-1}(u),E_{n}(w)]]).
\end{align*}
Therefore, $[E_{n-1}(u),[E_{n-1}(v),E_{n}(w)]]+[E_{n-1}(v),[E_{n-1}(u),E_{n}(w)]]=0$ when $u, v\neq w$. It can be further
shown that $[E_{n-1}(u),[E_{n-1}(v),E_{n}(w)]]+[E_{n-1}(v),[E_{n-1}(u),E_{n}(w)]]=0$ by using the technique of vertex operators.
This completes the proof of all relations in Theorem~\ref{thm:relinYgN}.
\epf

For $\g_N=\spa_{2n}$, set $\alpha_{j}=\epsilon_{j}-\epsilon_{j+1}$ for $j=1,\cdots,n-1$,$\beta_{i}=\epsilon_{n+i}$ for $i=1,2$,$\beta_{i}=\sqrt{2}\epsilon_{n+i}$ for $i=3,4,5$. We introduce the Fock space $\mathcal {F}=\Ac[a_{j,-k}(1\leq i\leq n+4,k\in \mathbb{Z}\backslash \{0\})]\widetilde{\bigotimes}\Ac[Q]$,
where $Q=\bigoplus_{i=1}^{n-1}\mathbb{Z}\alpha_{i}\oplus\bigoplus_{j=1}^{5}\mathbb{Z}\beta_{j}$, $\Ac[Q]$ is the group algebra of $Q$ over $\Ac$. On this space, we define the action of the operators $a_{i,k},\partial\epsilon_{j},e^{\epsilon_{j}}(1\leq i\leq n+4,1\leq j\leq n+5)$ by
\begin{equation}
a_{i,k}\cdot f\otimes e^{\beta}=
\begin{cases}
a_{i,k}f\otimes e^{\beta}&\text{if $k<0$}\\
[a_{i,k},f]\otimes e^{\beta}&\text{if $k>0$}
\end{cases}
\end{equation}
\begin{align*}
\partial\epsilon_{j}\cdot f\otimes e^{\beta}=&(\epsilon_{j},\beta)f\otimes e^{\beta},\\
e^{\epsilon_{j}}\cdot f\otimes e^{\beta}=&f\otimes e^{\epsilon_{j}+\beta},
\end{align*}
for $f\otimes e^{\beta}\in \mathcal {F}$.

\bthm
The following assignment defines a $\DY_\hbar(\spa_{2n})$-module structure on $\mathcal {F}$.\\
For $1\leq j \leq n-2$,
\begin{align*}
E_{j}(u)\mapsto
&exp[\sum_{k>0}\frac{a_{j,-k}}{k}\{(u-\frac{1}{4}\hbar)^{k}+(u+\frac{3}{4}\hbar)^{k}\}-\sum_{k>0}\frac{a_{j+1,-k}+a_{j-1,-k}}{k}(u+\frac{1}{4}\hbar)^{k}]\\
&exp[-\sum_{k>0}\frac{a_{j,k}}{k}(u-\frac{1}{4}\hbar)^{-k}]e^{\alpha_{j}}[(-1)^{j-1}(u-\frac{1}{4}\hbar)]^{\partial\alpha_{j}},
\end{align*}
\begin{align*}
F_{j}(u)\mapsto
&exp[-\sum_{k>0}\frac{a_{j,-k}}{k}\{(u+\frac{1}{4}\hbar)^{k}+(u-\frac{3}{4}\hbar)^{k}\}+\sum_{k>0}\frac{a_{j+1,-k}+a_{j-1,-k}}{k}(u-\frac{1}{4}\hbar)^{k}]\\
&exp[\sum_{k>0}\frac{a_{j,k}}{k}(u+\frac{1}{4}\hbar)^{-k}]e^{-\alpha_{j}}[(-1)^{j-1}(u+\frac{1}{4}\hbar)]^{-\partial\alpha_{j}},
\end{align*}
\begin{align*}
H_{j}^{+}(u)\mapsto
exp[-\sum_{k>0}\frac{a_{j,k}}{k}\{(u-\frac{1}{2}\hbar)^{-k}-(u+\frac{1}{2}\hbar)^{-k}\}](\frac{u+\frac{1}{2}\hbar}{u-\frac{1}{2}\hbar})^{-\partial\alpha_{j}},
\end{align*}
\begin{align*}
H_{j}^{-}(u)\mapsto
exp[-\sum_{k>0}\frac{a_{j,-k}}{k}&\{(u-\hbar)^{k}-(u+\hbar)^{k}\}\\
&+\sum_{k>0}\frac{a_{j+1,-k}+a_{j-1,-k}}{k}\{(u-\frac{1}{2}\hbar)^{k}-(u+\frac{1}{2}\hbar)^{k}\}].
\end{align*}

\begin{align*}
E_{n-1}(u)\mapsto
&exp[\sum_{k>0}\frac{a_{n-1,-k}}{k}\{(u-\frac{1}{4}\hbar)^{k}+(u+\frac{3}{4}\hbar)^{k}\}-\sum_{k>0}\frac{a_{n-2,-k}}{k}(u+\frac{1}{4}\hbar)^{k}]\\
&exp[-\sum_{k>0}\frac{a_{n-1,k}}{k}(u-\frac{1}{4}\hbar)^{-k}]e^{\alpha_{n-1}}[(-1)^{n-2}(u-\frac{1}{4}\hbar)]^{\partial\alpha_{n-1}}\\
&exp[\sum_{k>0}\frac{a_{n,-k}}{k}(u+\frac{1}{2}\hbar)^{k}]exp[\sum_{k>0}\frac{a_{n+1,k}}{k}(u-\frac{1}{2}\hbar)^{-k}]
e^{\beta_{1}}[(-1)^{n-2}(u-\frac{1}{2}\hbar)]^{\partial\beta_{2}},
\end{align*}
\begin{align*}
F_{n-1}(u)\mapsto
&exp[-\sum_{k>0}\frac{a_{n-1,-k}}{k}\{(u+\frac{1}{4}\hbar)^{k}+(u-\frac{3}{4}\hbar)^{k}\}+\sum_{k>0}\frac{a_{n-2,-k}}{k}(u-\frac{1}{4}\hbar)^{k}]\\
&exp[\sum_{k>0}\frac{a_{n-1,k}}{k}(u+\frac{1}{4}\hbar)^{-k}]e^{-\alpha_{n-1}}[(-1)^{n-2}(u+\frac{1}{4}\hbar)]^{-\partial\alpha_{n-1}}\\
&exp[\sum_{k>0}\frac{a_{n,-k}}{k}(u-\frac{1}{2}\hbar)^{k}]exp[\sum_{k>0}\frac{a_{n+1,k}}{k}(u+\frac{1}{2}\hbar)^{-k}]
e^{-\beta_{1}}[(-1)^{n-2}(u+\frac{1}{2}\hbar)]^{-\partial\beta_{2}},
\end{align*}
\begin{align*}
E_{n}(u)\mapsto
&exp[\sum_{k>0}\frac{a_{n+1,-k}}{k}(u+\frac{1}{2}\hbar)^{k}]exp[\sum_{k>0}\frac{a_{n,k}}{k}(u-\frac{1}{2}\hbar)^{-k}]e^{-\beta_{2}}[(-1)^{n-1}(u-\frac{1}{2}\hbar)]^{-\partial\beta_{1}}\\
&exp[\sum_{k>0}\frac{a_{n+2,-k}}{k}\{(u-\frac{1}{4}\hbar)^{k}+(u+\frac{3}{4}\hbar)^{k}\}+\sum_{k>0}\frac{a_{n+3,-k}}{k}\{(u-\frac{3}{4}\hbar)^{k}+(u+\frac{9}{4}\hbar)^{k}\}]\\
&exp[-\sum_{k>0}\frac{a_{n+2,k}}{k}(u-\frac{1}{4}\hbar)^{-k}-\sum_{k>0}\frac{a_{n+3,k}}{k}(u+\frac{1}{4}\hbar)^{-k}]
e^{\beta_{3}+\beta_{4}}(u-\frac{1}{4}\hbar)^{\partial\beta_{3}}(u+\frac{1}{4}\hbar)^{\partial\beta_{4}},
\end{align*}
\begin{align*}
F_{n}(u)\mapsto
&exp[\sum_{k>0}\frac{a_{n+1,-k}}{k}(u-\frac{1}{2}\hbar)^{k}]exp[\sum_{k>0}\frac{a_{n,k}}{k}(u+\frac{1}{2}\hbar)^{-k}]e^{-\beta_{2}}[(-1)^{n-1}(u+\frac{1}{2}\hbar)]^{-\partial\beta_{1}}\\
&exp[-\sum_{k>0}\frac{a_{n+2,-k}}{k}\{(u+\frac{1}{4}\hbar)^{k}+(u-\frac{3}{4}\hbar)^{k}\}-\sum_{k>0}\frac{a_{n+4,-k}}{k}\{(u+\frac{3}{4}\hbar)^{k}+(u-\frac{9}{4}\hbar)^{k}\}]\\
&exp[\sum_{k>0}\frac{a_{n+2,k}}{k}(u+\frac{1}{4}\hbar)^{-k}+\sum_{k>0}\frac{a_{n+4,k}}{k}(u-\frac{1}{4}\hbar)^{-k}]
e^{-\beta_{3}+\beta_{5}}(u+\frac{1}{4}\hbar)^{-\partial\beta_{3}}(u-\frac{1}{4}\hbar)^{\partial\beta_{5}},
\end{align*}
\begin{align*}
H_{n-1}^{+}(u)\mapsto
&exp[-\sum_{k>0}\frac{a_{n-1,k}}{k}\{(u-\frac{1}{2}\hbar)^{-k}-(u+\frac{1}{2}\hbar)^{-k}\}](\frac{u+\frac{1}{2}\hbar}{u-\frac{1}{2}\hbar})^{-\partial\alpha_{n-1}}\\
&exp[\sum_{k>0}\frac{a_{n,-k}}{k}\{(u+\frac{1}{4}\hbar)^{k}+(u-\frac{1}{4}\hbar)^{k}\}]exp[\sum_{k>0}\frac{a_{n+1,k}}{k}\{(u-\frac{3}{4}\hbar)^{-k}+(u+\frac{3}{4}\hbar)^{-k}\}]\\
&(\frac{u+\frac{3}{4}\hbar}{u-\frac{3}{4}\hbar})^{-\partial\beta_{2}},
\end{align*}
\begin{align*}
H_{n-1}^{-}(u)\mapsto
&exp[-\sum_{k>0}\frac{a_{n-1,-k}}{k}\{(u-\hbar)^{k}-(u+\hbar)^{k}\}+\sum_{k>0}\frac{a_{n-2,-k}}{k}\{(u-\frac{1}{2}\hbar)^{k}-(u+\frac{1}{2}\hbar)^{k}\}]\\
&exp[\sum_{k>0}\frac{a_{n,-k}}{k}\{(u+\frac{3}{4}\hbar)^{k}+(u-\frac{3}{4}\hbar)^{k}\}]exp[\sum_{k>0}\frac{a_{n+1,k}}{k}\{(u-\frac{1}{4}\hbar)^{-k}+(u+\frac{1}{4}\hbar)^{-k}\}]\\
&(\frac{u+\frac{1}{4}\hbar}{u-\frac{1}{4}\hbar})^{-\partial\beta_{2}},
\end{align*}
\begin{align*}
H_{n}^{+}(u)\mapsto
&exp[\sum_{k>0}\frac{a_{n+1,-k}}{k}\{(u+\frac{1}{4}\hbar)^{k}+(u-\frac{1}{4}\hbar)^{k}\}]exp[\sum_{k>0}\frac{a_{n,k}}{k}\{(u-\frac{3}{4}\hbar)^{-k}+(u+\frac{3}{4}\hbar)^{-k}\}]\\
&e^{-2\beta_{2}}[(-1)^{n-1}(u-\frac{3}{4}\hbar)]^{-\partial\beta_{1}}[(-1)^{n-1}(u+\frac{3}{4}\hbar)]^{-\partial\beta_{1}}\\
&exp[\sum_{k>0}\frac{a_{n+3,-k}}{k}\{(u-\hbar)^{k}+(u+2\hbar)^{k}\}-\sum_{k>0}\frac{a_{n+4,-k}}{k}\{(u+\hbar)^{k}+(u-2\hbar)^{k}\}]\\
&exp[-\sum_{k>0}\frac{a_{n+2,k}}{k}\{(u-\frac{1}{2}\hbar)^{-k}-(u+\frac{1}{2}\hbar)^{-k}\}+\sum_{k>0}\frac{a_{n+4,k}-a_{n+3,k}}{k}u^{-k}]\\
&e^{\beta_{4}+\beta_{5}}(\frac{u+\frac{1}{2}\hbar}{u-\frac{1}{2}\hbar})^{-\partial\beta_{3}}u^{\partial\beta_{4}+\partial\beta_{5}},
\end{align*}
\begin{align*}
H_{n}^{-}(u)\mapsto
&exp[\sum_{k>0}\frac{a_{n+1,-k}}{k}\{(u+\frac{3}{4}\hbar)^{k}+(u-\frac{3}{4}\hbar)^{k}\}]exp[\sum_{k>0}\frac{a_{n,k}}{k}\{(u-\frac{1}{4}\hbar)^{-k}+(u+\frac{1}{4}\hbar)^{-k}\}]\\
&e^{-2\beta_{2}}[(-1)^{n-1}(u-\frac{1}{4}\hbar)]^{-\partial\beta_{1}}[(-1)^{n-1}(u+\frac{1}{4}\hbar)]^{-\partial\beta_{1}}\\
&exp[-\sum_{k>0}\frac{a_{n+2,-k}}{k}\{(u-\hbar)^{k}-(u+\hbar)^{k}\}+\sum_{k>0}\frac{a_{n+3,-k}}{k}\{(u-\frac{1}{2}\hbar)^{k}+(u+\frac{5}{2}\hbar)^{k}\}\\
&-\sum_{k>0}\frac{a_{n+4,-k}}{k}\{(u+\frac{1}{2}\hbar)^{k}+(u-\frac{5}{2}\hbar)^{k}\}]exp[-\sum_{k>0}\frac{a_{n+3,k}}{k}(u+\frac{1}{2}\hbar)^{-k}\\
&+\sum_{k>0}\frac{a_{n+4,k}}{k}(u-\frac{1}{2}\hbar)^{-k}]e^{\beta_{4}+\beta_{5}}(u+\frac{1}{2}\hbar)^{\partial\beta_{4}}(u-\frac{1}{2}\hbar)^{\partial\beta_{5}}.
\end{align*}

\ethm

\bpf
The relations can be checked similarly as type B.
\epf

For $\g_N=\oa_{2n}$, set $\alpha_{j}=\epsilon_{j}-\epsilon_{j+1}$ for $j=1,\cdots,n-1$, $\alpha_{n}=\epsilon_{n-1}+\epsilon_{n}$. We introduce the Fock space $\mathcal {F}=\Ac[a_{j,-k}(1\leq i\leq n,k\in \mathbb{Z}\backslash \{0\})]\widetilde{\bigotimes}\Ac[Q]$,
where $Q=\bigoplus_{i=1}^{n}\mathbb{Z}\alpha_{i}$, $\Ac[Q]$ is the group algebra of $Q$ over $\Ac$. On this space, we define the action of the operators $a_{j,k},\partial\epsilon_{j},e^{\epsilon_{j}}(1\leq j\leq n)$ by
\begin{equation}
a_{j,k}\cdot f\otimes e^{\beta}=
\begin{cases}
a_{j,k}f\otimes e^{\beta}&\text{if $k<0$}\\
[a_{j,k},f]\otimes e^{\beta}&\text{if $k>0$}
\end{cases}
\end{equation}
\begin{align*}
\partial\epsilon_{j}\cdot f\otimes e^{\beta}=&(\epsilon_{j},\beta)f\otimes e^{\beta},\\
e^{\epsilon_{j}}\cdot f\otimes e^{\beta}=&f\otimes e^{\epsilon_{j}+\beta},
\end{align*}
for $f\otimes e^{\beta}\in \mathcal {F}$.

\bthm\label{VertexD}
The following assignment defines a $\DY_\hbar(\oa_{2n})$-module structure on $\mathcal {F}$.\\
For $1\leq j \leq n-3$,
\begin{align*}
E_{j}(u)\mapsto
&exp[\sum_{k>0}\frac{a_{j,-k}}{k}\{(u-\frac{1}{4}\hbar)^{k}+(u+\frac{3}{4}\hbar)^{k}\}-\sum_{k>0}\frac{a_{j+1,-k}+a_{j-1,-k}}{k}(u+\frac{1}{4}\hbar)^{k}]\\
&exp[-\sum_{k>0}\frac{a_{j,k}}{k}(u-\frac{1}{4}\hbar)^{-k}]
e^{\alpha_{j}}[(-1)^{j-1}(u-\frac{1}{4}\hbar)]^{\partial\alpha_{j}},
\end{align*}
\begin{align*}
F_{j}(u)\mapsto
&exp[-\sum_{k>0}\frac{a_{j,-k}}{k}\{(u+\frac{1}{4}\hbar)^{k}+(u-\frac{3}{4}\hbar)^{k}\}+\sum_{k>0}\frac{a_{j+1,-k}+a_{j-1,-k}}{k}(u-\frac{1}{4}\hbar)^{k}]\\
&exp[\sum_{k>0}\frac{a_{j,k}}{k}(u+\frac{1}{4}\hbar)^{-k}]
e^{-\alpha_{j}}[(-1)^{j-1}(u+\frac{1}{4}\hbar)]^{-\partial\alpha_{j}},
\end{align*}
\begin{align*}
H_{j}^{+}(u)\mapsto
exp[-\sum_{k>0}\frac{a_{j,k}}{k}\{(u-\frac{1}{2}\hbar)^{-k}-(u+\frac{1}{2}\hbar)^{-k}\}](\frac{u+\frac{1}{2}\hbar}{u-\frac{1}{2}\hbar})^{-\partial\alpha_{j}},
\end{align*}
\begin{align*}
H_{j}^{-}(u)\mapsto
exp[-\sum_{k>0}\frac{a_{j,-k}}{k}&\{(u-\hbar)^{k}-(u+\hbar)^{k}\}\\
&+\sum_{k>0}\frac{a_{j+1,-k}+a_{j-1,-k}}{k}\{(u-\frac{1}{2}\hbar)^{k}-(u+\frac{1}{2}\hbar)^{k}\}].
\end{align*}

\begin{align*}
E_{n-2}(u)\mapsto
&exp[\sum_{k>0}\frac{a_{n-2,-k}}{k}\{(u-\frac{1}{4}\hbar)^{k}+(u+\frac{3}{4}\hbar)^{k}\}-\sum_{k>0}\frac{a_{n-1,-k}+a_{n-3,-k}+a_{n,-k}}{k}\\
&(u+\frac{1}{4}\hbar)^{k}]exp[-\sum_{k>0}\frac{a_{n-2,k}}{k}(u-\frac{1}{4}\hbar)^{-k}]
e^{\alpha_{n-2}}[(-1)^{n-3}(u-\frac{1}{4}\hbar)]^{\partial\alpha_{n-2}},
\end{align*}
\begin{align*}
F_{n-2}(u)\mapsto
&exp[-\sum_{k>0}\frac{a_{n-2,-k}}{k}\{(u+\frac{1}{4}\hbar)^{k}+(u-\frac{3}{4}\hbar)^{k}\}+\sum_{k>0}\frac{a_{n-1,-k}+a_{n-3,-k}+a_{n,-k}}{k}\\
&(u-\frac{1}{4}\hbar)^{k}]exp[\sum_{k>0}\frac{a_{n-2,k}}{k}(u+\frac{1}{4}\hbar)^{-k}]
e^{-\alpha_{n-2}}[(-1)^{n-3}(u+\frac{1}{4}\hbar)]^{-\partial\alpha_{n-2}},
\end{align*}
\begin{align*}
H_{n-2}^{+}(u)\mapsto
exp[-\sum_{k>0}\frac{a_{n-2,k}}{k}\{(u-\frac{1}{2}\hbar)^{-k}-(u+\frac{1}{2}\hbar)^{-k}\}](\frac{u+\frac{1}{2}\hbar}{u-\frac{1}{2}\hbar})^{-\partial\alpha_{n-2}},
\end{align*}
\begin{align*}
H_{n-2}^{-}(u)\mapsto
exp&[-\sum_{k>0}\frac{a_{n-2,-k}}{k}\{(u-\hbar)^{k}-(u+\hbar)^{k}\}+\sum_{k>0}\frac{a_{n-1,-k}+a_{n-3,-k}+a_{n,-k}}{k}\\
&\{(u-\frac{1}{2}\hbar)^{k}-(u+\frac{1}{2}\hbar)^{k}\}].
\end{align*}

\begin{align*}
E_{n-1}(u)\mapsto
&exp[\sum_{k>0}\frac{a_{n-1,-k}}{k}\{(u-\frac{1}{4}\hbar)^{k}+(u+\frac{3}{4}\hbar)^{k}\}-\sum_{k>0}\frac{a_{n-2,-k}}{k}(u+\frac{1}{4}\hbar)^{k}]\\
&exp[-\sum_{k>0}\frac{a_{n-1,k}}{k}(u-\frac{1}{4}\hbar)^{-k}]
e^{\alpha_{n-1}}[(-1)^{n-2}(u-\frac{1}{4}\hbar)]^{\partial\alpha_{n-1}},
\end{align*}
\begin{align*}
F_{n-1}(u)\mapsto
&exp[-\sum_{k>0}\frac{a_{n-1,-k}}{k}\{(u+\frac{1}{4}\hbar)^{k}+(u-\frac{3}{4}\hbar)^{k}\}+\sum_{k>0}\frac{a_{n-2,-k}}{k}(u-\frac{1}{4}\hbar)^{k}]\\
&exp[\sum_{k>0}\frac{a_{n-1,k}}{k}(u+\frac{1}{4}\hbar)^{-k}]
e^{-\alpha_{n-1}}[(-1)^{n-2}(u+\frac{1}{4}\hbar)]^{-\partial\alpha_{n-1}},
\end{align*}
\begin{align*}
H_{n-1}^{+}(u)\mapsto
exp[-\sum_{k>0}\frac{a_{n-1,k}}{k}\{(u-\frac{1}{2}\hbar)^{-k}-(u+\frac{1}{2}\hbar)^{-k}\}](\frac{u+\frac{1}{2}\hbar}{u-\frac{1}{2}\hbar})^{-\partial\alpha_{n-1}},
\end{align*}
\begin{align*}
H_{n-1}^{-}(u)\mapsto
exp[-\sum_{k>0}\frac{a_{n-1,-k}}{k}&\{(u-\hbar)^{k}-(u+\hbar)^{k}\}\\
&+\sum_{k>0}\frac{a_{n-2,-k}}{k}\{(u-\frac{1}{2}\hbar)^{k}-(u+\frac{1}{2}\hbar)^{k}\}].
\end{align*}

\begin{align*}
E_{n}(u)\mapsto
&exp[\sum_{k>0}\frac{a_{n,-k}}{k}\{(u-\frac{1}{4}\hbar)^{k}+(u+\frac{3}{4}\hbar)^{k}\}-\sum_{k>0}\frac{a_{n-2,-k}}{k}(u+\frac{1}{4}\hbar)^{k}]\\
&exp[-\sum_{k>0}\frac{a_{n,k}}{k}(u-\frac{1}{4}\hbar)^{-k}]
e^{\alpha_{n}}[(-1)^{N}(u-\frac{1}{4}\hbar)]^{\partial\alpha_{n}},
\end{align*}
\begin{align*}
F_{n}(u)\mapsto
&exp[-\sum_{k>0}\frac{a_{n,-k}}{k}\{(u+\frac{1}{4}\hbar)^{k}+(u-\frac{3}{4}\hbar)^{k}\}+\sum_{k>0}\frac{a_{n-2,-k}}{k}(u-\frac{1}{4}\hbar)^{k}]\\
&exp[\sum_{k>0}\frac{a_{n,k}}{k}(u+\frac{1}{4}\hbar)^{-k}]
e^{-\alpha_{n}}[(-1)^{n}(u+\frac{1}{4}\hbar)]^{-\partial\alpha_{n}},
\end{align*}
\begin{align*}
H_{n}^{+}(u)\mapsto
exp[-\sum_{k>0}\frac{a_{n,k}}{k}\{(u-\frac{1}{2}\hbar)^{-k}-(u+\frac{1}{2}\hbar)^{-k}\}](\frac{u+\frac{1}{2}\hbar}{u-\frac{1}{2}\hbar})^{-\partial\alpha_{n}},
\end{align*}
\begin{align*}
H_{n}^{-}(u)\mapsto
exp[-\sum_{k>0}\frac{a_{n,-k}}{k}\{(u-\hbar)^{k}-(u+\hbar)^{k}\}+\sum_{k>0}\frac{a_{n-2,-k}}{k}\{(u-\frac{1}{2}\hbar)^{k}-(u+\frac{1}{2}\hbar)^{k}\}].
\end{align*}

\ethm

\bpf
The relations can be checked similarly as type B.
\epf

We remark that the analogous constructions of our level one modules in the affine Lie algebras
are different from those in \cite{Be, FJ, JK, JKM}. Our vertex operator representation covers
both symmetric and non-symmetric classical cases in type $B, C, D$ for the double Yangian algebras.
We remark that the symmetric case--Theorem~\ref{VertexD} is a special case of \cite[Thm.~5.5]{GRW2}.
The only difference is that the operators are shifted to account for slightly different presentations.

\bigskip

\centerline{\bf Acknowledgments}

\medskip

The research is supported by
the National Natural Science Foundation of China (grant nos. 11531004 and 11701182), Simons Foundation
(grant no. 523868), and the Australian Research Council (grant no. DP180101825).


\bigskip





\begin{thebibliography}{99}



\bibitem{amr:otr}
{D. Arnaudon, A. Molev and E. Ragoucy},
{\it On the $R$-matrix realization of Yangians
and their representations},
Annales Henri Poincar\'e 7 (2006), 1269--1325.

\bibitem{Be}
D. Bernard, {\it Vertex operator representations of the quantum affine algebra $U_q(B^{(1)}_r)$}, Lett. Math. Phys. 17 (1989), 239--245.

\bibitem{bk:pp}
{J. Brundan and A. Kleshchev},
{\it Parabolic presentations of the Yangian $\Y(\gl_n)$},
Comm. Math. Phys. 254 %{\bf 254}
(2005), 191-220.

\bibitem{cp:gq}
{V. Chari and A. Pressley},
{\it A guide to quantum groups},
Cambridge University Press, Cambridge, 1994.

\bibitem{C} I. V. Cherednik, {\it On the method of constructing factorized S-matrices in elementary functions}, Theor. Math. Phys. 43 (1980), 117-119. (Russian)

\bibitem{DHHZ} X. M. Ding, B. Y. Hou, B. Yuan Hou, L. Zhao,
{\it Free boson realization of $DY_h(gl_N)_k$}, J. Math. Phys. 39 (1998), 2273. %arXiv:hep-th/9709016.

\bibitem{Dr} V. G. Drinfeld, {\it Hopf algebras and the quantum Yang-Baxter equation}, Sov. Math. Dokl. 32 (1985), 254-258.

\bibitem{D}
V. G. Drinfeld, {\it A new realization of Yangians and of quantum affine algebras}, (Russian) Dokl. Akad. Nauk SSSR 296 (1987), 13-17; transl.: Soviet Math. Dokl. 36 (1988), 212-216.

\bibitem{FM} E. Frenkel and E. Mukhin, {\it Combinatorics of q-characters of finite-dimensional representations of quantum affine algebras}, Comm. Math. Phys. 216 (2001), 23-57.

\bibitem{FR1} E. Frenkel and N. Y. Reshetikhin, {\it Quantum affine algebras and deformations of the Virasoro and $\mathcal W$-algebras}, Comm. Math. Phys. 178 (1996), 237-264.

\bibitem{FJ}
I. B. Frenkel and N. Jing, {\it Vertex representations of quantum affine algebras}, Proc. Nat'l. Acad. Sci. USA {85} (1988), 9373-9377.

\bibitem{FR2} I. B. Frenkel and N. Y. Reshetikhin, {\it Quantum affine algebras and holonomic difference equations}, Comm. Math. Phys. 146 (1992), 1-60.

\bibitem{gr:dm}
{I. M. Gelfand and V. S. Retakh}, {\it Determinants
of matrices over noncommutative rings}, {Funct. Anal. Appl.} 25 %{\bf 25}
(1991), 91-102.


\bibitem{GRW1} N. Guay, V. Regelskis and C. Wendlandt, {\it Equivalences between three presentations of orthogonal and symplectic Yangians},
Lett. Math. Phys. 109 (2019), 327-379. %arXiv:1706.05176.


\bibitem{GRW2} N. Guay, V. Regelskis and C. Wendlandt, {\it Vertex representations for Yangians of Kac-Moody algebras},
J. \'Ec. Polytech. Math. 6 (2019), 665-706. % (arXiv:1804.04081).


\bibitem{io:br}
{K. Iohara},
{\it Bosonic representations of Yangian double $\DY_{\hbar}(\g)$ with $\g=\gl_N, \mathfrak{sl}_N$},
J. Phys. A 29 (1996), %no. 15,
4593-4621.
%{\tt arXiv:9603033v1}.

\bibitem{Jb} M. Jimbo, {\it A q-difference analogue of $U(g)$ and the Yang-Baxter equation}, Lett. Math. Phys. 10 (1985), 63-69.

\bibitem{Jb1}
M. Jimbo, {\it Quantum R-matrix for the generalized Toda system}, Comm. Math. Phys. 102 (1986), 537-547.

\bibitem{JK}
N. Jing, Y. Koyama, {\it Vertex operators of admissible modules of $U_q(C^{(1)}_n)$}, J. Algebra 205 (1998), 294-316.

\bibitem{JKM}
N. Jing, Y. Koyama, K. C. Misra, {\it Bosonic realizations of $U_q(C^{(1)}_n)$}, J. Algebra 200 (1998), 155-172.



\bibitem{ji:cen}
{N. Jing, S. Kozic, A. Molev and F. Yang},
{\it Center of the quantum affine vertex algebra in type A}, J. Algebra 496 (2018), 138-186.
%{\tt arXiv:1603.00237v1}.


\bibitem{ji:iso}
{N. Jing, M. Liu and A. Molev},
{\it Isomorphism between the R-matrix and Drinfeld
presentations of Yangian in types B,C and D}, Comm. Math. Phys. 361 (2018), 827-872.
%{\tt arXiv:1705.08155v1}.

\bibitem{CKa} C. Kassel, {\em Quantum Groups. Vol. 155, GTM. New York: Springe-Verlag}, 1995.

\bibitem{KT}
S. M. Khoroshkin and V. N. Tolstoy, {\it Yangian double}, Lett. Math. Phys. 36 (1996), 373-402.
%the authors present the triangular decomposition of DY(g) and a factorization for the canonical pairing of the Yangian with its dual inside Y0(g). The theory of the Cartan-Weyl basis allows one to describe explicitly the universal R-matrix.

\bibitem{Ko}
S. Kozic, {\it Commutative operators for double Yangian $DY(sl_n)$}, Glas. Mat. Ser. III 53 %(73)
(2018), 97-113.

\bibitem{KR}
P. P. Kulish and N. Y. Reshetikhin, {\it The quantum linear problem for the sine-Gordon equation and higher representations}, Questions in quantum
field theory and statistical physics, 2. Zap. Nauch. Sem. LOMI 101 (1981), 101-110.
(Russian)


\bibitem{m:yc}
A. Molev,
{\it Yangians and classical Lie algebras}, Mathematical
Surveys and Monographs, 143. American Mathematical Society,
Providence, RI, 2007.


\bibitem{R} N. Y. Reshetikhin, {\it Quasitriangular Hopf algebras and invariants of links}, (Russian) Algebra i Analiz 1 (1989),% no. 2,
169-188; translation in Leningrad Math. J. 1 (1990), %no. 2,
491-513

\bibitem{RT} N. Y. Reshetikhin and V. G. Turaev, {\it Invariants of 3-manifolds via link polynomials and quantum groups}, Invent. Math. 103 (1991), 547-597.

\bibitem{Ta} V. O. Tarasov, {\it Structures of L-operators for the R-matrix of the XXZ-model}, Theor. Math. Phys. 61 (1984), 1065-1072.

\bibitem{T} A. Tsymbaliuk, {\it PBW theorems and shuffle realizations for $U_v(L{\mathfrak sl}_n)$, $U_{v1,v2}(L{\mathfrak sl}_n)$, $U_v(L{\mathfrak sl}(m|n))$}, arXiv:1808.09536.

\bibitem{XZ} Y. Xu and R. B. Zhang, {\em Drinfeld realisations and vertex operator representations of quantum affine superalgebras},
arXiv:1802.09702.

\bibitem{zz:rf}
{A. B. Zamolodchikov and Al. B. Zamolodchikov},
{\it Factorized $S$-matrices in two dimensions as the exact solutions
of certain relativistic quantum field models},
{Ann. Phys.} 120 %{\bf 120}
(1979), 253-291.

\bibitem{Z} H. F. Zhang, {\it Representations of quantum affine superalgebras}, Math. Z. 278 (2014), 663-703.



\end{thebibliography}
\end{document}